\documentclass[preprint,11pt,numbers,sort&compress]{elsarticle}

\usepackage{soul}
\usepackage{psfrag}
\usepackage{a4wide}
\usepackage{rotating}
\usepackage{color}
\usepackage{amssymb}
\usepackage{amsmath}
\usepackage{amsthm}
\usepackage{verbatim}
\usepackage{cancel}
\usepackage{nicematrix}
\usepackage{multirow}
\usepackage{tikz}

\newcommand{\f}[1]{\mathbf{#1}}
\newcommand{\ab}[1]{\boldsymbol{#1}}
\def\bfm#1{\boldsymbol{#1}}
\newcommand{\bb}[1]{\bfm{#1}}
\newcommand{\N}{\mathbb N}

\newcommand{\R}{\mathbb R}

\newcommand{\V}{\mathcal{V}}

\newcommand{\W}{\mathcal{W}}

\newcommand{\LL}{i_0}
\newcommand{\RR}{i_1}
\newcommand{\g}{\gamma}
\newcommand{\gC}{\gamma}

\newcommand{\Side}{\tau}
\newcommand{\dd}{\partial}

\DeclareMathOperator{\Rank}{rank}
\DeclareMathOperator{\Diag}{diag}

\DeclareMathOperator{\Span}{span}

\theoremstyle{definition}
\newtheorem{ex}{Example}

\newtheorem{rem}{Remark}

\newproof{pf}{proof}
\advance\textheight by 0.4cm
\advance\topmargin by -0.2cm

\definecolor{gold}{rgb}{1,0.7,0}
\definecolor{dred}{rgb}{0.92,0,0}
\definecolor{dgreen}{rgb}{0,0.6,0}

\begin{document}

\begin{frontmatter}

\title{An Isogeometric Tearing and Interconnecting (IETI) method for solving high order partial differential equations over planar multi-patch geometries}

\cortext[cor]{Corresponding author}

\author[vil]{Mario Kapl}
\ead{m.kapl@fh-kaernten.at}

\author[slo1]{Alja\v z Kosma\v c}
\ead{aljaz.kosmac@iam.upr.si}

\author[slo1]{Vito Vitrih\corref{cor}}
\ead{vito.vitrih@upr.si}

\address[vil]{Department of Engineering $\&$ IT, Carinthia University of Applied Sciences, Villach, Austria}

\address[slo1]{IAM and FAMNIT, University of Primorska, Koper, Slovenia}

\begin{abstract}
We present a novel method for solving high-order partial differential equations (PDEs) over planar multi-patch geometries demonstrated on the basis of the polyharmonic equation of order~$m$, $m \geq 1$, which is a particular linear elliptic PDE of order~$2m$. Our approach is based on the concept of Isogeometric Tearing and Interconnecting (IETI) \cite{KlPeSaJu12} and allows to couple the numerical solution of the PDE with $C^s$-smoothness,  $s \geq m-1$, across the edges of the multi-patch geometry. The proposed technique relies on the use of a particular class of multi-patch geometries, called bilinear-like $G^s$ multi-patch parameterizations~\cite{KaVi20b}, to represent the multi-patch domain. The coupling between the neighboring patches is done via the use of Lagrange multipliers and leads to a saddle point problem, which can be solved efficiently first by a small dual problem for a subset of the Lagrange multipliers followed by local, parallelizable problems on the single patches for the coefficients of the numerical solution.
Several numerical examples of solving the polyharmonic equation of order $m=2$ and $m=3$, i.e. the biharmonic and triharmonic equation,  respectively, over different multi-patch geometries are shown to demonstrate the potential of our IETI method for high-order problems.
\end{abstract}
\begin{keyword}
isogeometric analysis; 
Galerkin method; 
$C^s$-smoothness; tearing and interconnecting; multi-patch domain; polyharmonic equation
\MSC[2010] 65N30 \sep 65D17 \sep 68U07
\end{keyword}

\end{frontmatter}

\section{Introduction}

Multi-patch geometries with possibly extraordinary vertices, i.e. unstructured quadrilateral meshes, are a common and powerful tool in computer-aided geometric design~\cite{HoLa93,Fa97} to model complex physical domains. Solving high-order PDEs, such as the polyharmonic equation of order~$m$, $m \geq 1$, which is a particular linear elliptic PDE of order~$2m$, over these multi-patch geometries via the Galerkin method and the concept of isogeometric analysis~\cite{ANU:9260759,CottrellBook,HuCoBa04} requires at least the use of globally $C^s$-smooth spline spaces with $s \geq m-1$. The construction of such $C^s$-smooth isogeometric multi-patch spline spaces for the case $s \geq 1$ can be roughly classified into two approaches. While the first strategy directly constructs the $C^s$-smooth spline space by generating a basis of the space which can then be used for discretizing and solving the PDE, the second approach usually employs the ``standard'' isogeometric spline basis functions on the single patches which are in general discontinuous across the edges of the multi-patch domain and enforces the $C^s$-smoothness of the solution of the PDE e.g. by adapting the variational form of the considered problem.

Most of the existing techniques of the first approach have been studied for the case $s=1$ and construct exactly $C^1$-smooth isogeometric spline spaces. These different constructions depend on the used parameterization of multi-patch domain including the use of $C^1$-smooth multi-patch parameterizations with singularities at the extraordinary vertices, e.g.~\cite{NgPe16,ToSpHu17,WeLiQiHuZhCa22}, $C^1$-smooth multi-patch parameterizations with $G^1$-smooth caps in the vicinity of the extraordinary vertices, e.g.~\cite{KaPe17,KaPe18,WeFaLiWeCa23}, particular class of $G^1$-smooth multi-patch parameterizations called analysis-suitable $G^1$ multi-patch parameterizations, e.g.~\cite{CoSaTa16,FaJuKaTa22,KaSaTa17b}, general $G^1$-smooth multi-patch parameterizations, e.g.~\cite{ChAnRa18,ChAnRa19,MaMaMo2024,mourrain2015geometrically} as well as $G^1$-smooth multi-patch parameterizations with singularities at so-called scaling centers, e.g. \cite{ArReKlSi23,ReArSiKl23}. In particular for planar multi-patch domains, also methods have been developed for the case $s=2$, e.g.~\cite{KaVi17b,KaVi17c,KaVi19a,KaVi20}, and for an arbitrary~$s \geq 1$, e.g.~\cite{KaVi20b,KaKoVi24b,KaKoVi24c,KaKoVi24}, where these methods rely on the use of so-called bilinear-like $G^s$ multi-patch geometries~\cite{KaVi20b} to represent the multi-patch domain. In addition to exactly $C^s$-smooth spline spaces, approximately $C^1$-smooth spline spaces have been studied by constructing an approximate $C^1$-smooth basis via a null space computation~\cite{SaJu21}, an explicit representation~\cite{WeTa21,WeTa22} or subdivision based techniques, e.g.~\cite{CIOrSch00,RiAuFe16,ZhSaCi18}.    

Examples of the second strategy in the context of isogeometric analysis, which have been used to solve fourth-order PDEs over multi-patch domains, are the Nitsche's method, e.g. \cite{BeEvMcTa21,Guo2015881,Nguyen2014,WeTa22}, the penalty method, e.g. \cite{CoLoBu21,CoKiBu21,LeLiMaKiReGa20,DuRoSa17}, as well as the mortar method, e.g. \cite{BeLoSaTa23,Bouclier2017,DiSchWoHe19,SchDiWoKlHe19,MiZoScBoTh21}. While the Nitsche's method and the penalty methods add so-called coupling terms to the variational form of the problem to ensure the smoothness across the edges of the multi-patch domain, the mortar methods impose the smoothness conditions across the edges by means of Lagrange multipliers. However, all these methods have in common that the obtained numerical solution of the PDE is just approximately $C^1$ (or $C^s$-smooth). An alternative to the mentioned techniques is the concept of Isogeometric Tearing and Interconecting (IETI), firstly studied in \cite{KlPeSaJu12}, which is as the mortar method a Lagrangian multiplier based approach but which also allows the computation of exactly smooth numerical solutions of the PDE. 

The idea of the IETI-method has its origin in the Finite Element Tearing and Interconnecting (FETI) method, e.g.~\cite{Farhat1991,FarhatFETI-DP,Klawonn2000,Klawonn2006}, which is a domain decomposition technique that efficiently solves a linear elliptic PDE first by a small dual problem for the involved Lagrange multipliers and then by local, parallelizable problems on the single subdomains for the degrees of freedom of the numerical simulation. While the initial work~\cite{Farhat1991} requires the use of a pseudo-inverse to locally invert the stiffness matrices on the single subdomains, the techniques~\cite{FarhatFETI-DP,Klawonn2000,Klawonn2006} avoid this by introducing and solving an appropriate dual-primal formulation of the considered problem. In the last years, the concept of FETI has been extended to the framework of isogeometric analysis by developing several IETI-methods for solving PDEs over structured and unstructured quadrilateral meshes, e.g.~\cite{KlPeSaJu12,SoglTakacs_IETI_Elasticity,IETI_LowRank2024,HoLa17,SoTa23,MoSaSchTaTa23,WiZaScPa21}. However, the existing FETI and IETI-methods are mainly limited to second order PDEs or just to the case of structured quadrilateral meshes for PDEs of higher order.

In this paper, we will extend the concept of IETI method also to the solving of high-order PDEs, namely in our case to the polyharmonic equation of order $m$, $m\geq 1$, over planar multi-patch geometries with possibly extraordinary vertices, i.e. over unstructured planar quadrilateral meshes. Our proposed approach will be based on the parameterization of the considered multi-patch geometry by a bilinear-like $G^s$ multi-patch parameterization~\cite{KaVi20b} for $s \geq m-1$, and on the solving of a saddle point problem, where the required smoothness conditions of the numerical solution across the edges will be imposed by means of Lagrange multipliers. For an efficient solving of the saddle point problem, a dual-primal formulation of the problem will be presented, which will compute first a subset of the Lagrange multipliers via a small linear problem and then the remaining Lagrange multipliers as well as the degrees of freedom of the numerical solution via local linear problems for the single patches which can be solved in a parallel way. The power of our IETI method for high order problems will be demonstrated on the basis of several examples of solving the polyharmonic equation of order $m=2$ and $m=3$, which are also called the biharmonic and triharmonic equation, respectively. Thereby, all the numerical results will show optimal convergence behavior in the studied $H^m$-seminorm with respect to $h$-refinement. 

The remainder of the paper is organized as follows. In Section~\ref{sec:prel}, we will present some required preliminaries including the presentation of the used multi-patch setting, the class of bilinear-like $G^s$ multi-patch geometries as well as the concept of $C^s$-smooth isogeometric spline functions over planar bilinear-like $G^s$ multi-patch geometries. Section~\ref{sec:IETI} will first introduce the studied problem, namely the solving of the polyharmonic equation of order~$m$, $m \geq 1$, over planar bilinear-like $G^s$ multi-patch geometries, and will then describe the IETI-based method to efficiently solve this problem by means of a dual-primal formulation. Some technical details and proofs about the proposed IETI technique will be presented in the Appendix. In Section~\ref{section_numerical_examples}, we will show several numerical examples of solving the biharmonic and triharmonic equation, which are the polyharmonic equation of order $m=2$ and $m=3$, respectively. Finally, we conclude the paper in Section~\ref{sec:Conclusion}.

\section{Preliminaries} \label{sec:prel}

We will introduce the setting of our used planar multi-patch domains, 
describe a specific class of multi-patch geometries, called 
bilinear-like $G^s$-smooth multi-patch parameterizations~\cite{KaVi20b}, 
to 
represent these multi-patch domains, as well as recall the concept of $C^s$-smooth isogeometric multi-patch spline spaces.

\subsection{The planar multi-patch domain}

Let $\Omega \subset \R^2$ be an open domain whose closure~$\overline{\Omega}$ represents a planar multi-patch domain given as the disjoint union
\[
\overline{\Omega} = \bigcup_{i \in \mathcal{I}_{\Omega}} \Omega^{(i)}  \; \dot{\cup}  \bigcup_{i \in \mathcal{I}_{\Gamma}} \Gamma^{(i)} \; \dot{\cup} \bigcup_{i \in \mathcal{I}_{\Xi}} \bfm{\Xi}^{(i)}
\]
with open inner quadrilateral patches~$\Omega^{(i)}$, $i \in \mathcal{I}_{\Omega}^{I}$, and open boundary quadrilateral patches $\Omega^{(i)}$, $i \in \mathcal{I}_{\Omega}^{B}$, such that $\mathcal{I}_{\Omega} = \mathcal{I}_{\Omega}^I \dot{\cup} \mathcal{I}_{\Omega}^B$, with open inner edges~$\Gamma^{(i)}$, $i \in \mathcal{I}_\Gamma^{I}$, and open boundary edges~$\Gamma^{(i)}$, $i \in \mathcal{I}_{\Gamma}^{B}$, such that $\mathcal{I}_{\Gamma} = \mathcal{I}_{\Gamma}^I \dot{\cup} \mathcal{I}_{\Gamma}^B$, and with inner vertices $\bfm{\Xi}^{(i)}$, $ i \in \mathcal{I}_{\Xi}^I$, and boundary vertices $\bfm{\Xi}^{(i)}$, $ i \in \mathcal{I}_{\Xi}^B$, such that $\mathcal{I}_{\Xi} = \mathcal{I}_{\Xi}^I \dot{\cup} \mathcal{I}_{\Xi}^B$. The intersection $\overline{\Omega^{(i_0)}} \cap \overline{\Omega^{(i_1)}}$ of the closures of any two patches~$\Omega^{(i_0)}$ and $\Omega^{(i_1)}$, $i_0,i_1 \in \mathcal{I}_{\Omega}$, $i_0 \neq i_1$, is either empty, a vertex~$\bfm{\Xi}^{(i)}$ for an $i \in \mathcal{I}_{\Xi}$, or the closure~$\overline{\Gamma^{(i)}}$ of an inner edge~$\Gamma^{(i)}$ for an $i \in \mathcal{I}_{\Gamma}^I$. We will denote by $\nu_i$ the patch valency of a vertex~$\bfm{\Xi}^{(i)}$, $i \in \mathcal{I}_{\Xi}$, and by $\mathcal{I}_{\Xi}^{B,\nu \geq 2} \subset \mathcal{I}_{\Xi}^{B}$ the index set collecting all indices of boundary vertices~$\bfm{\Xi}^{(i)}$, $i \in \mathcal{I}^{B}_{\Xi}$, with a patch valency~$\nu_i \geq 2$.  

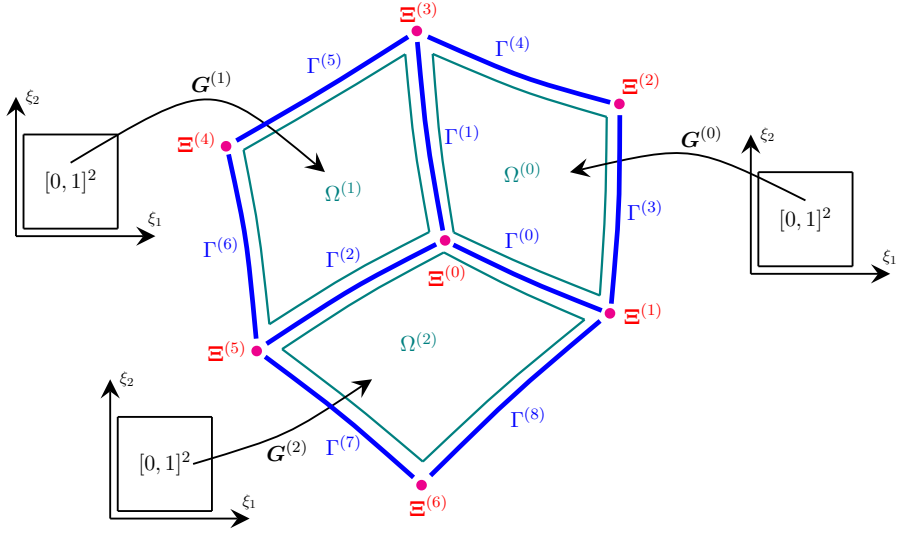
\begin{figure}
\begin{center}
\begin{tabular}{c}
\resizebox{0.8\textwidth}{!}{
 \begin{tikzpicture}
  \coordinate(A) at (0,0); \coordinate(B) at (1.8,2.5); \coordinate(C) at (0,4.2); \coordinate(D) at (-1.8,2.5);  \coordinate(E) at (-3.8,2.4); 
  \coordinate(F) at (-2.6,-0.4); \coordinate(G) at (-3,-2.9); \coordinate(H) at (-0.5,-2.4); \coordinate(I) at (2,-3.1); 
  \coordinate(J) at (2.5,-0.9); \coordinate(K) at (3.9,1.5);
  \draw[line width=0.12mm, teal] (0.0830276, -0.0202994) .. controls (0.00843625, 0.447714) .. (-0.0283338, 0.928053);
  \draw[line width=0.12mm,teal] (0.0830276, -0.0202994) .. controls (0.469181, -0.200027) .. (0.858765, -0.354949);
  \draw[line width=0.12mm,teal] (0.895323, 0.593403) .. controls (0.889936, 0.109089) .. (0.858765, -0.354949);
  \draw[line width=0.12mm,teal] (0.895323, 0.593403) .. controls (0.447681, 0.724473) .. (-0.0283338, 0.928053);
  \node[scale=0.3,teal] at (0.45,0.3) {$\Omega^{(0)}$};
  \draw[-stealth,line width=0.1mm] (1.95,0.15) .. controls (1.35,0.45) ..(0.7, 0.3);
  \node[scale=0.3] at (1.4,0.49) {$\ab{G}^{(0)}$};

  \draw[line width=0.12mm,teal] (-0.0463954, -0.0178091) .. controls (-0.464553, -0.232213) .. (-0.895251, -0.507768);
  \draw[line width=0.12mm,teal] (-0.0463954, -0.0178091) .. controls (-0.129422, 0.445037) .. (-0.174048, 0.926239);
  \draw[line width=0.12mm,teal] (-1.03215, 0.41779) .. controls (-0.590866, 0.670787) .. (-0.174048, 0.926239);
  \draw[line width=0.12mm,teal] (-1.03215, 0.41779) .. controls (-0.94911, -0.00646313) .. (-0.895251, -0.507768);
  \node[scale=0.3,teal] at (-0.5,0.2) {$\Omega^{(1)}$};
  \draw[-stealth,line width=0.1mm] (-1.95,0.35) .. controls (-1.2,0.8) .. (-0.6, 0.3);
  \node[scale=0.3] at (-1.2,0.78) {$\ab{G}^{(1)}$};

  \draw[line width=0.12mm,teal] (0.0309216, -0.126322) .. controls (0.402589, -0.321429) .. (0.777823, -0.483336);
  \draw[line width=0.12mm,teal] (0.0309216, -0.126322) .. controls (-0.392534, -0.352443) .. (-0.827679, -0.63994);
  \draw[line width=0.12mm,teal] (-0.0807774, -1.23732) .. controls (-0.452036, -0.923429) .. (-0.827679, -0.63994);
  \draw[line width=0.12mm,teal] (-0.0807774, -1.23732) .. controls (0.338466, -0.841568) .. (0.777823, -0.483336);
  \node[scale=0.3,teal] at (-0.1,-0.6) {$\Omega^{(2)}$};
  \draw[-stealth,line width=0.1mm] (-1.3,-1.25) .. controls (-0.8,-1.1) .. (-0.35, -0.8);
  \node[scale=0.3] at (-0.8,-1.18) {$\ab{G}^{(2)}$};

\draw[thick, magenta, fill=magenta] (0.0375, -0.0625) circle (0.14mm);

\draw[thick, magenta, fill=magenta] (-0.1125, 1.05) circle (0.14mm);

\draw[thick, magenta, fill=magenta] (0.9125, -0.45) circle (0.14mm);

\draw[thick, magenta, fill=magenta] (0.9625, 0.6625) circle (0.14mm);

\draw[thick, magenta, fill=magenta] (-1.125, 0.4375) circle (0.14mm);

\draw[thick, magenta, fill=magenta] (-0.9625, -0.65) circle (0.14mm);

\draw[thick, magenta, fill=magenta] (-0.0875, -1.3625) circle (0.14mm);

\draw[line width=0.25mm, blue] (0.02525, -0.00984375) .. controls (-0.0625, 0.478125) .. (-0.10975, 0.991406);

\draw[line width=0.25mm, blue] (0.08125, -0.0860312) .. controls (0.475, -0.278125) .. (0.86875, -0.434781);

\draw[line width=0.25mm, blue] (0.96475, 0.602719) .. controls (0.9625, 0.084375) .. (0.91975, -0.398531);

\draw[line width=0.25mm, blue] (0.9135, 0.670594) .. controls (0.45, 0.796875) .. (-0.054, 1.01934);

\draw[line width=0.25mm, blue] (-0.010125, -0.0829688) .. controls (-0.45, -0.309375) .. (-0.910125, -0.611719);

\draw[line width=0.25mm, blue] (-1.07022, 0.467531) .. controls (-0.596875, 0.740625) .. (-0.158969, 1.01878);

\draw[line width=0.25mm, blue] (-1.11153, 0.394406) .. controls (-1.01563, -0.046875) .. (-0.965281, -0.584344);

\draw[line width=0.25mm, blue] (-0.130063, -1.32153) .. controls (-0.51875, -0.978125) .. (-0.917563, -0.680281);

\draw[line width=0.25mm, blue] (-0.039875, -1.31153) .. controls (0.4, -0.878125) .. (0.860125, -0.490281);


\draw[line width=0.1mm] (-2.2, 0) -- (-1.7, 0);
\draw[line width=0.1mm] (-2.2, 0.5) -- (-1.7, 0.5);
\draw[line width=0.1mm] (-2.2, 0) -- (-2.2, 0.5);
\draw[line width=0.1mm] (-1.7, 0) -- (-1.7, 0.5);
\draw[-stealth,line width=0.1mm] (-2.24, -0.04) -- (-1.5, -0.04);
\draw[-stealth,line width=0.1mm] (-2.24, -0.04) -- (-2.24, 0.7);
\node[scale=0.3] at (-1.95,0.25) {$[0,1]^2$};
\node[scale=0.23] at (-1.5,0.05) {$\xi_1$};
\node[scale=0.23] at (-2.15,0.7) {$\xi_2$};


\draw[line width=0.1mm] (-1.7, -1.5) -- (-1.2, -1.5);
\draw[line width=0.1mm] (-1.7, -1) -- (-1.2, -1);
\draw[line width=0.1mm] (-1.7, -1.5) -- (-1.7, -1);
\draw[line width=0.1mm] (-1.2, -1.5) -- (-1.2, -1);
\draw[-stealth,line width=0.1mm] (-1.74, -1.54) -- (-1.0, -1.54);
\draw[-stealth,line width=0.1mm] (-1.74, -1.54) -- (-1.74, -0.8);
\node[scale=0.3] at (-1.45,-1.25) {$[0,1]^2$};
\node[scale=0.23] at (-1.0,-1.45) {$\xi_1$};
\node[scale=0.23] at (-1.65,-0.8) {$\xi_2$};


\draw[line width=0.1mm] (1.7, -0.2) -- (2.2, -0.2);
\draw[line width=0.1mm] (1.7, 0.3) -- (2.2, 0.3);
\draw[line width=0.1mm] (1.7, -0.2) -- (1.7, 0.3);
\draw[line width=0.1mm] (2.2, -0.2) -- (2.2, 0.3);
\draw[-stealth,line width=0.1mm] (1.66, -0.24) -- (2.4, -0.24);
\draw[-stealth,line width=0.1mm] (1.66, -0.24) -- (1.66, 0.5);
\node[scale=0.3] at (1.95,0.05) {$[0,1]^2$};
\node[scale=0.23] at (2.4,-0.15) {$\xi_1$};
\node[scale=0.23] at (1.75,0.5) {$\xi_2$};

\node[scale=0.3,red] at (0.05,-0.25) {$\bfm{\Xi}^{(0)}$};
\node[scale=0.3,red] at (1.09,-0.45) {$\bfm{\Xi}^{(1)}$};
\node[scale=0.3,red] at (1.08,0.77) {$\bfm{\Xi}^{(2)}$};
\node[scale=0.3,red] at (-0.1,1.15) {$\bfm{\Xi}^{(3)}$};
\node[scale=0.3,red] at (-1.28,0.46) {$\bfm{\Xi}^{(4)}$};
\node[scale=0.3,red] at (-1.12,-0.65) {$\bfm{\Xi}^{(5)}$};
\node[scale=0.3,red] at (-0.05,-1.47) {$\bfm{\Xi}^{(6)}$};

\node[scale=0.3,blue] at (0.45,-0.05) {$\Gamma^{(0)}$};
\node[scale=0.3,blue] at (0.13,0.5) {$\Gamma^{(1)}$};
\node[scale=0.3,blue] at (-0.5,-0.15) {$\Gamma^{(2)}$};
\node[scale=0.3,blue] at (1.1,0.1) {$\Gamma^{(3)}$};
\node[scale=0.3,blue] at (0.4,0.97) {$\Gamma^{(4)}$};
\node[scale=0.3,blue] at (-0.6,0.88) {$\Gamma^{(5)}$};
\node[scale=0.3,blue] at (-1.15,-0.1) {$\Gamma^{(6)}$};
\node[scale=0.3,blue] at (-0.5,-1.15) {$\Gamma^{(7)}$};
\node[scale=0.3,blue] at (0.48,-1) {$\Gamma^{(8)}$};

\end{tikzpicture}
 }
\end{tabular}
\end{center}
\caption{The schematic illustration of an example of a planar three patch domain~$\overline{\Omega}$ composed of the open patches~$\Omega^{(i)}$, $i\in \{0,1,2\}$, open edges~$\Gamma^{(i)}$, $i \in \{0, \ldots,8 \}$, and vertices~$\bfm{\Xi}^{(i)}$, $i \in \{0,\ldots,6 \}$, and with its corresponding geometry mappings~$\ab{G}^{(i)}$, $i \in \{0,1,2\}$.}
 \label{fig:three_patch_isogeometric_analysis}
\end{figure}

Each patch $\overline{\Omega^{(i)}}$ is parameterized by a bijective and regular geometry mapping (also called patch parameterization)~$\ab{G}^{(i)}$,
\begin{align*}
 \ab{G}^{(i)}: [0,1]^{2}  \rightarrow \R^{2}, \quad 
 \bb{\xi} =(\xi_1,\xi_2) \mapsto
 \ab{G}^{(i)}(\bb{\xi}) = \ab{G}^{(i)}(\xi_1,\xi_2), \quad i \in \mathcal{I}_{\Omega},
\end{align*}
such that $\overline{\Omega^{(i)}} = \ab{G}^{(i)}([0,1]^{2})$, cf.~Fig.~\ref{fig:three_patch_isogeometric_analysis}, with 
\[
\ab{G}^{(i)} \in \mathcal{S}_{h}^{\ab{p},\ab{r}}([0,1]^2) \times \mathcal{S}^{\ab{p},\ab{r}}_{h}([0,1]^2),
\]
where $\mathcal{S}_{h}^{\ab{p},\ab{r}}([0,1]^2)$ with $\ab{p}=(p,p)$ and $\ab{r}=(r,r)$ is the tensor-product spline space $\mathcal{S}_h^{p,r}([0,1]) \otimes \mathcal{S}_h^{p,r}([0,1])$ on the unit-square~$[0,1]^2$, obtained from the univariate spline space~$\mathcal{S}_h^{p,r}([0,1])$ on the unit interval~$[0,1]$ of degree~$p$, regularity~$r$ and mesh size~$h=\frac{1}{k+1}$ with $k$ different, uniformly distributed inner knots. We will denote by $\ab{G}$ the resulting multi-patch parameterization consisting of the single patch parameterizations~$\ab{G}^{(i)}$, $i \in \mathcal{I}_{\Omega}$, and we will further denote by $N_{j}^{p,r}$ and $N_{\ab{j}}^{\ab{p},\ab{r}}=N_{j_1}^{p,r}N_{j_2}^{p,r}$ with $j \in \mathcal{J}= \{0,1,\ldots,n-1\}$, $\ab{j}=(j_1,j_2) \in \bb{\mathcal{J}}=\mathcal{J} \times \mathcal{J}$, and $n= \dim \mathcal{S}_h^{p,r}([0,1]) = p+1+k(p-r)$, the B-splines of the spline spaces~$\mathcal{S}_h^{p,r}([0,1])$ and $\mathcal{S}_h^{\ab{p},\ab{r}}([0,1]^2)$, respectively. 

\subsection{The bilinear-like $G^s$ multi-patch parameterization}

In this work, we will represent the multi-patch parameterization~$\ab{G}$ of the multi-patch domain~$\overline{\Omega}$ by a specific multi-patch geometry called bilinear-like $G^s$ multi-patch parameterization~\cite{KaVi20b}. Let $s, r \geq 0$ with $r \geq s$. A $C^0$-smooth multi-patch parameterization~$\ab{G}$ is called bilinear-like $G^s$ if for each inner edge~$\Gamma^{(i)}$, $i \in \mathcal{I}^{I}_{\Gamma}$, with $\overline{\Gamma^{(i)}} = \overline{\Omega^{(i_0)}} \cap  \overline{\Omega^{(i_1)}}$, $i_0,i_1 \in \mathcal{I}_{\Omega}$, assuming without loss of generality that the two associated patch parameterizations~$\ab{G}^{(i_0)}$ and $\ab{G}^{(i_0)}$ are parameterized as
\begin{equation} \label{eq:standard_par}
\ab{G}^{(i_0)}(0,\xi) = \ab{G}^{(i_1)}(0,\xi), \quad \xi \in [0,1],
\end{equation}
there exist linear functions $\alpha^{(i,i_0)}:[0,1] \rightarrow \R$, $\alpha^{(i,i_1)}:[0,1] \rightarrow \R$, $\beta^{(i,i_0)}:[0,1] \rightarrow \R$ and $\beta^{(i,i_1)}: [0,1] \rightarrow \R$ such that
\begin{equation*}
\ab{G}^{(i,i_0)}_{\ell}(\xi) = \ab{G}^{(i,i_1)}_{\ell}(\xi) := \ab{G}^{(i)}_{\ell}(\xi), \quad \xi \in [0,1], \quad \ell =0,1,\ldots,s,  
\end{equation*}
with
\[
\ab{G}_{\ell}^{(i,\tau)}(\xi) = \left(\alpha^{(i,\tau)}(\xi) \right)^{-\ell} \partial_1^{\ell}\ab{G}^{(\tau)}(0,\xi) - \sum_{j=0}^{\ell-1} \binom{\ell}{j} \left( \frac{\beta^{(i,\tau)}(\xi)}{\alpha^{(i,\tau)}(\xi)}\right)^{\ell-j} \partial^{\ell -j}\ab{G}^{(i)}_{j}(\xi), \quad \tau\in \{i_0,i_1 \}.
\]
The class of bilinear-like $G^s$ multi-patch parameterizations is hence characterized by having linear connectivity functions $\alpha^{(i,i_0)}$, $\alpha^{(i,i_1)}$, $\beta^{(i,i_0)}$ and $\beta^{(i,i_1)}$ between the patches like bilinear multi-patch geometries, and contains {for any $s$} the bilinear multi-patch parameterizations as a subclass. Note that the class of bilinear-like $G^s$ multi-patch parameterizations is for $s=0$ just the class of $C^0$-smooth multi-patch parameterizations, and coincides for $s=1$ exactly with the class of analysis-suitable $G^1$-multi-patch parameterizations~\cite{CoSaTa16}. The use of bilinear-like $G^s$ multi-patch parameterizations allows to model multi-patch domains with curved edges and boundaries and enables further to solve high-order partial differential equations with $C^s$-smooth functions over the multi-patch domain with optimal convergence rates as numerically demonstrated e.g. for solving the biharmonic equation~\cite{KaSaTa19b,FaKaKoVi24} and the triharmonic equation~\cite{KaKoVi24b,KaVi19a} using the Galerkin approach and for solving the Poisson's equation~\cite{KaVi20} and the biharmonic equation~\cite{KaKoVi24,KaKoVi24c} via isogeometric collocation. For $s \geq 1$, the construction of bilinear-like $G^s$ multi-patch parameterizations was considered e.g. for~$s=1$ (i.e. for analysis-suitable $G^1$ multi-patch geometries) in \cite{KaSaTa17b,FaJuKaTa22,FaKaKoVi24}, for~$s=2$ in \cite{KaVi17c,KaSaTa19b,KaKoVi24b}, for~$s=3$ in \cite{KaVi20b} and for~$s=4$ in~\cite{KaKoVi24,KaKoVi24c}. 

\subsection{The space of $C^s$-smooth isogeometric spline functions} \label{subsec:Cs_space}

We recall the concept of $C^s$-smooth isogeometric spline functions over $G^s$-smooth bilinear-like multi-patch parameterizations which was first studied for an arbitrary smoothness~$s\geq 0$ in~\cite{KaVi20b}.

We denote by $\V_h(\overline{\Omega})$ the space of isogeometric spline functions over the multi-patch domain~$\overline{\Omega}$ (with respect to the bilinear-like $G^s$ multi-patch parameterization~$\ab{G}$ and mesh size~$h$) given as 
\begin{equation*}
\V_h(\overline{\Omega}) = \left\{ \phi \in L^2(\overline{\Omega}): \; \phi |_{\overline{\Omega}^{(i)}} \in {\mathcal{S}_{h}^{\ab{p},\ab{r}}([0,1]^{2})} \circ (\ab{G}^{(i)})^{-1} , \; 
i \in \mathcal{I}_{\Omega}   \right\} ,
\end{equation*}
and by $\V_h^s(\overline{\Omega})$ the corresponding space of $C^s$-smooth isogeometric spline functions, i.e
\begin{equation*}
\V_h^s(\overline{\Omega}) = \V_h(\overline{\Omega}) \cap C^s(\overline{\Omega}).
\end{equation*}
It was shown in~\cite[Theorem~2]{KaVi20b} that an isogeometric spline function $\phi \in \V_h(\overline{\Omega})$ belongs to the space~$\mathcal{V}_h^s(\overline{\Omega})$ if and only if for any inner edge~$\Gamma^{(i)}$, $i \in \mathcal{I}^{I}_{\Gamma}$, with $\overline{\Gamma^{(i)}} = \overline{\Omega^{(i_0)}} \cap  \overline{\Omega^{(i_1)}}$, $i_0,i_1 \in \mathcal{I}_{\Omega}$, assuming without loss of generality that the two associated patch parameterizations~$\ab{G}^{(i_0)}$ and $\ab{G}^{(i_1)}$ are parameterized as in~\eqref{eq:standard_par}, we have for all $\xi \in [0,1]$
\begin{equation}   \label{eq:gC}
 \g_\ell^{(i,\LL)}[\phi](\xi) = \g_\ell^{(i,\RR)}[\phi](\xi) =: \gC^{(i)}_\ell[\phi](\xi), \quad \ell =0,1,\ldots,s,
 \end{equation}
 with
 \begin{equation*}   \label{eq:gC2}
 \g_\ell^{(i,\Side)}[\phi](\xi) = \left(\alpha^{(i,\Side)}(\xi)\right)^{-\ell}\, \partial_1^\ell 
 \left(\phi \circ \ab{G}^{(\tau)}\right)(0,\xi) - \sum_{j=0}^{\ell-1} {\ell \choose j} 
 \left(\frac{\beta^{(i,\Side)}(\xi)}{\alpha^{(i,\Side)}(\xi)}\right)^{\ell-j}  \dd^{\ell-j} \gC_j^{(i)}[\phi](\xi),
 \end{equation*}
for $\Side\in \{\LL,\RR\}$. The spline functions $\gC_{\ell}^{(i)}[\phi]$ represent for $\ell=0$ the trace of the function $\phi$ along the inner edge~$\Gamma^{(i)}$ and for $\ell=1,\ldots,s$ specific derivatives of order~$\ell$ of the function~$\phi$ across the inner edge~$\Gamma^{(i)}$, cf.~\cite{CoSaTa16}. By means of~Eq.~\eqref{eq:gC}, the $C^s$-smooth spline space~$\V_h^s(\overline{\Omega})$ can be characterized as
\begin{equation} \label{eq:spaceVs}
\V_h^s(\overline{\Omega}) = \left\{ \phi \in \V_h(\overline{\Omega}): \;  \g_\ell^{(i,\LL)}[\phi](\xi) = \g_\ell^{(i,\RR)}[\phi](\xi), \; \xi \in [0,1], \, \ell=0,1,\ldots, s,\; i \in \mathcal{I}_{\Gamma}^I  \right\}.
\end{equation}
Since the dimension of $\V_h^s(\overline{\Omega})$ is highly dependent on the valencies of the inner vertices~$\bfm{\Xi}^{(i)}$, $ i \in \mathcal{I}_{\Xi}^I$, as well as on the parameterizations of the single patch parameterizations~$\ab{G}^{(i)}$, $i \in \mathcal{I}_{\Omega}$, cf.~\cite{KaVi20b,KaVi17b}, usually the subspace $\W_h^s(\overline{\Omega}) \subseteq \V_h^s(\overline{\Omega})$ given as
\begin{equation} \label{eq:spaceWs}
\W_h^s(\overline{\Omega}) = \left\{ \begin{array}{ll }\phi \in \V_h^s(\overline{\Omega}) : \; & \gC^{(i)}_\ell[\phi] \in \mathcal{S}_h^{p-\ell,r+s-\ell}([0,1]) , \, \ell=0,1,\ldots, s,  \, i \in \mathcal{I}_{\Gamma}^{I}, \mbox{ and } \\
&\phi \in C^{2s}(\bfm{\Xi}^{(i)}), i \in \mathcal{I}_{\Xi}^{I} \end{array} \right\} 
\end{equation}
is considered instead. To get an $h$-refinable space~$\W_h^s(\overline{\Omega})$, we have to require $p \geq 2s+1$, $s \leq r \leq p-(s+1)$, and $h \leq \frac{p-r-s}{3s-r+1}$, cf.~\cite{KaVi20b}.  The $C^s$-smooth spline space~$\W_h^s(\overline{\Omega})$ maintains the optimal approximation properties of the entire $C^s$-smooth spline space~$\V_h^s(\overline{\Omega})$, as numerically shown in \cite{KaVi20b}, and has a dimension which is independent of the valencies of the inner vertices~$\bfm{\Xi}^{(i)}$, $ i \in \mathcal{I}_{\Xi}^I$, and of the single patch parameterizations~$\ab{G}^{(i)}$, $i \in \mathcal{I}_{\Omega}$, cf.~\cite{KaVi20b}.

\section{The Isogeometric Tearing and Interconnecting (IETI) method} \label{sec:IETI}

We will introduce a IETI-based method for solving high order PDEs, given by the polyharmonic equations of order~$m$, $m \geq 1$, over planar multi-patch domains with bilinear-like $G^s$ multi-patch parameterizations with 
 $s \geq m-1$. 
For this purpose, we will first state the considered problem and then present a general framework of the method for an arbitrary $m$ and $s$ satisfying $s \geq m-1 \geq 0$.      

\subsection{The problem statement}

Let $m \geq 1$, 
$s \geq m-1$, and let $\Omega \subset \R^2$ be an open domain whose closure~$\overline{\Omega}$ is a planar multi-patch domain with a bilinear-like $G^s$ multi-patch parameterization~$\ab{G}$, cf. Section~\ref{sec:prel}. We consider as model problem the polyharmonic equation of order $m$ given by
\begin{align} \label{eq:polyharmonic}
(-\triangle)^{m} u (\bfm{x}) & =   f(\bfm{x}), \quad 
\ab{x} \in \Omega,  \nonumber \\[-0.6cm]
& \\
\mathcal{B}^{\ell}[u](\bfm{x}) & = g_{\ell}(\bfm{x}), \; \; \, 
\ab{x} \in \partial \Omega, \quad \ell = 0,1,\ldots,m-1, \nonumber
\end{align}
where
\[
\mathcal{B}^{\ell}[u](\ab{x}) = \begin{cases}
\nabla^{\ell}  u (\bfm{x}) & \mbox{if $\ell$ is even},\\
\left( \nabla^{\ell}  u (\bfm{x}) \right)^T \bfm{n} &  \mbox{if $\ell$ is odd},
\end{cases}  
\]
$f: \Omega \rightarrow \R$, $g_{\ell}: \partial \Omega \rightarrow \R$, $\ell=0,1,\ldots,m-1$, are sufficiently smooth functions, and $\bfm{n}$ 
is the outward unit normal vector at the boundary~$\partial \Omega$. 
The polyharmonic equation~\eqref{eq:polyharmonic} is called the Poisson's equation for $m=1$, the biharmonic equation for $m=2$ and the triharmonic equation for $m=3$. The weak form of problem~\eqref{eq:polyharmonic} is to find $u\in H_g^{m}(\Omega)$ such that
\begin{equation} \label{eq:weak_polyharmonic}
a(u,v) = F(v), \quad \mbox{for all } v \in H^{m}_0(\Omega),
\end{equation}
where $a$ is the bilinear form 
 \[
a(u,v) = \sum_{i \in \mathcal{I}_\Omega} a^{(i)}(u,v), \quad a^{(i)}(u,v) = \int_{\Omega^{(i)}} \left( \nabla^{m}u(\bfm{x}) \right)^T \nabla^{m}v(\bfm{x})  d \Omega^{(i)},   
 \]
 $F$ is the linear functional
\[
F(v) = \sum_{i \in \mathcal{I}_\Omega}F^{(i)}(v), \quad F^{(i)}(v)= \int_{\Omega^{(i)}} f(\bfm{x}) v(\bfm{x}) d\Omega^{(i)},
\]
and the spaces $H_g^{m}(\Omega)$ and $H_0^{m}(\Omega)$ are given as
\begin{equation*}
H_g^{m}(\Omega) = \left\{ \phi \in H^{m}(\Omega): 
\mathcal{B}^{\ell}[\phi](\bfm{x})=g_{\ell}, \, \ell=0,1,\ldots,m-1,
\mbox{ for }\bfm{x} \in \partial \Omega \right\}
\end{equation*}
and
\begin{equation*}
H_0^{m}(\Omega) = \left\{ \phi \in H^{m}(\Omega): 
\mathcal{B}^{\ell}[\phi](\bfm{x})=0, \, \ell=0,1,\ldots,m-1,
\mbox{ for }\bfm{x} \in \partial \Omega \right\},
\end{equation*}
respectively. Problem~\eqref{eq:weak_polyharmonic} is equivalent to computing $u \in H_g^{m}(\Omega)$ via the minimization problem
\begin{equation} \label{eq:minproblem}
u=\underset{v \in  H_g^{m}(\Omega)}{\arg \min} \frac{1}{2} a(v,v)-F(v).
\end{equation}

We denote by $\W_{g,h}^s(\overline{\Omega})$ the subspace of $\W_{h}^{s}(\overline{\Omega})$, whose functions~$\phi$ fulfill the boundary conditions $\mathcal{B}^{\ell}[\phi](\bfm{x})=g_{\ell}(\bfm{x}), \, \ell=0,\ldots,m-1$, for $\ab{x} \in \partial \Omega$, i.e.
\[
\W^s_{g,h}(\overline{\Omega}) =  \left\{ \phi \in \W^s_{h}(\overline{\Omega}): 
\mathcal{B}^{\ell}[\phi](\bfm{x})=g_{\ell}(\bfm{x}), \, \ell=0,1,\ldots,m-1, \mbox{ for }\bfm{x} \in \partial \Omega
\right\} .
\]
Note that if there does not exist a function~$\phi \in \W^{s}_h(\overline{\Omega})$ such that $\mathcal{B}^{\ell}[\phi](\bfm{x})=g_{\ell}(\bfm{x})$, $\ell=0,\ldots,m-1$, for $\bfm{x} \in \partial \Omega$, then the boundary data~$g_{\ell}$, $\ell=0,\ldots,m-1$, is just projected to $\W^{s}_h(\overline{\Omega})$.
By means of Eqs.~\eqref{eq:spaceVs} and \eqref{eq:spaceWs}, the space~$\W_{g,h}^s(\overline{\Omega})$ can be described as
\begin{equation} \label{eq:char_Wsgh}
\W_{g,h}^s(\overline{\Omega}) = \left\{ \begin{array}{ll} \hspace{-0.1cm} \phi \in \V_h(\overline{\Omega}): & \hspace{-0.05cm} \g_\ell^{(i,\LL)}[\phi](\xi) = \g_\ell^{(i,\RR)}[\phi](\xi), \; \xi \in [0,1], \, \ell=0,\ldots, s,\; i \in \mathcal{I}_{\Gamma}^I, \mbox{ and}  \\ & \gC^{(i)}_\ell[\phi] \in \mathcal{S}_h^{p-\ell,r+s-\ell}([0,1]) , \, \ell=0,1,\ldots, s,  \, i \in \mathcal{I}_{\Gamma}^{I}, \mbox{ and } \\
& \phi \in C^{2s}(\bfm{\Xi}^{(i)}), i \in \mathcal{I}_{\Xi}^{I}, \mbox{ and } \\
& \mathcal{B}^{\ell}[\phi](\bfm{x})=g_{\ell}(\bfm{x}), \, \ell=0,1,\ldots,m-1, \mbox{ for } \bfm{x} \in \partial \Omega   \end{array}   \right\}.
\end{equation}
By applying Galerkin projection using the finite dimensional subspace~$\W_{g,h}^s(\overline{\Omega}) \subseteq H_{g}^{m}(\Omega)$, problem~\eqref{eq:minproblem} transforms to finding $u_h \in \W^{s}_{g,h}(\overline{\Omega})$ via 
\begin{equation} \label{eq:minproblem2} 
u_h=\underset{v_h \in \W^{s}_{g,h}(\overline{\Omega})}{\arg \min} \frac{1}{2} a(v_h,v_h)-F(v_h).
\end{equation}
Using characterization~\eqref{eq:char_Wsgh} for the subspace~$\W^{s}_{g,h}(\overline{\Omega})$, problem~\eqref{eq:minproblem2} is equivalent to computing $u_h \in \W^{s}_{g,h}(\overline{\Omega})$ by solving the minimization problem
\begin{equation} \label{eq:minproblem3} 
u_h=\underset{v_h \in \V_h(\overline{\Omega})}{\arg \min} \frac{1}{2} a(v_h,v_h)-F(v_h)
\end{equation}
subject to
\begin{align}
 & \g_\ell^{(i,\LL)}[v_h](\xi) = \g_\ell^{(i,\RR)}[v_h](\xi), \; \xi \in [0,1], \, \ell=0,1,\ldots, s,\; i \in \mathcal{I}_{\Gamma}^I, \mbox{ and } \label{eq:constraintsI1} \\  &\gC^{(i)}_\ell[v_h] \in \mathcal{S}_h^{p-\ell,r+s-\ell}([0,1]) , \, \ell=0,1,\ldots, s,  \, i \in \mathcal{I}_{\Gamma}^{I}, \mbox{ and } \label{eq:constraintsI2} \\
 &v_h \in C^{2s}(\bfm{\Xi}^{(i)}), i \in \mathcal{I}_{\Xi}^{I}, \mbox{ and } \label{eq:constraintsV}\\
 &\mathcal{B}^{\ell}[v_h](\bfm{x})=g_{\ell}(\bfm{x}), \, \ell=0,1,\ldots,m-1, \mbox{ for }\ab{x} \in \partial \Omega. \label{eq:constraintsB}
\end{align}
The isogeometric spline space~$\V_h(\overline{\Omega})$ can be expressed as the direct sum of smaller subspaces corresponding to the single patches~$\Omega^{(i)}$, $i \in \mathcal{I}_\Omega$, i.e., 
\begin{equation*} \label{eq:V}
    \V_h(\overline{\Omega}) =  \bigoplus_{i \in \mathcal{I}_\Omega} \V_h^{(i)}(\overline{\Omega})
\end{equation*}
with
\begin{equation*}
\V_{h}^{(i)}(\overline{\Omega}) =  \Span \left\{ \phi_{\ab{j}}^{(i)} |\;  \ab{j} \in \bb{\mathcal{J}}  \right\},
\end{equation*}
where the isogeometric spline functions~$\phi^{(i)}_{\ab{j}}$, $j\in \bb{\mathcal{J}}$, are defined as
\begin{equation*} 
{\phi}^{(i)}_{\ab{j}}(\bfm{x})  = 
\begin{cases}
   (N_{\ab{j}}^{\ab{p},\ab{r}}\circ (\ab{G}^{(i)})^{-1})(\bfm{x}) ;
\mbox{ if }\f \, \bfm{x} \in \overline{\Omega^{(i)}},
\\ 0 \quad \mbox{ if }\f \, \bfm{x} \in \overline{\Omega} \backslash \overline{\Omega^{(i)}}.
\end{cases} 
\end{equation*}
For each $i \in \mathcal{I}_{\Omega}$, the functions $\phi^{(i)}_{\ab{j}}$, $\ab{j}\in \bb{\mathcal{J}}$, form a basis of $\V_h^{(i)}(\overline{\Omega})$, and hence, all functions $\phi^{(i)}_{\ab{j}}$, $\ab{j}\in \bb{\mathcal{J}}$, $i \in \mathcal{I}_{\Omega}$, build a basis of $\V_h(\overline{\Omega})$. Therefore, each function $u_h \in \V_h(\overline{\Omega})$ can be represented as
\begin{equation} \label{eq:spline_representation}
u_h(\ab{x}) =\sum_{i \in \mathcal{I}_{\Omega}} \sum_{\ab{j} \in \bb{\mathcal{J}}} u_{\ab{j}}^{(i)} \phi_{\ab{j}}^{(i)}(\ab{x}) = \sum_{i \in \mathcal{I}_{\Omega}} \left(\ab{u}^{(i)}\right)^T \ab{\phi}^{(i)}(\ab{x}) = 
\ab{u}^T \ab{\phi}(\ab{x})
\end{equation}
with coefficients~$u_{\ab{j}}^{(i)} \in \R$, where
\[
\ab{u}^{(i)} = [u_{\ab{j}}^{(i)}]_{\ab{j} \in \bb{\mathcal{J}}} \quad \mbox{and} \quad
\ab{\phi}^{(i)} = [\phi_{\ab{j}}^{(i)}]_{\ab{j} \in \bb{\mathcal{J}}} 
\]
are the column vectors of the coefficients and of the basis functions corresponding to the patch~$\Omega^{(i)}$, and where
\[
\ab{u}  = [\ab{u}^{(i)}]_{i \in \mathcal{I}_{\Omega}} \quad \mbox{and} \quad
\ab{\phi} = [\ab{\phi}^{(i)}]_{i \in \mathcal{I}_{\Omega}}
\]
are the column vectors of all coefficients and of all basis functions. The constraints~\eqref{eq:constraintsI1}--\eqref{eq:constraintsB} for a function~$v_h \in \V_h(\overline{\Omega})$ with the spline representation~$v_h(\ab{x})=\ab{v}^T \ab{\phi}(\ab{x})$ are linear with respect to the coefficients~$\ab{v}$, and can be expressed as
\[
\ab{C} \ab{v} = \ab{c}, 
\]
with a matrix~$\ab{C} \in \R^{q \times n^2|\mathcal{I}_{\Omega}|}$ and a column vector~$\ab{c} \in \R^{n^2|\mathcal{I}_{\Omega}|}$ for some~$q \in \N$, assuming without loss of generality that all linear equations $\ab{C} \ab{v} = \ab{c}$ are linearly independent, which means that $\Rank (\ab{C}) = q$. Then, the minimization problem~\eqref{eq:minproblem3} with constraints~\eqref{eq:constraintsI1}--\eqref{eq:constraintsB} becomes
\begin{equation*} 
u_h=\underset{\\ \ab{C} \ab{v}=\ab{c}}{\underset{v_h \in \V_h(\overline{\Omega})}{\arg \min}} \frac{1}{2} a(v_h,v_h)-F(v_h),
\end{equation*}
which is further equivalent to the following saddle point problem: Find $u_h \in \V_h(\overline{\Omega})$ with the spline representation~$u_h(\ab{x})=\ab{u}^T \ab{\phi}(\ab{x})$ as given in~\eqref{eq:spline_representation} and Lagrange multipliers~$\ab{\lambda} \in \R^{q}$, such that
\begin{equation} \label{eq:large_problem}
\left( \begin{array}{cc}
\ab{K} & \ab{C}^T \\
\ab{C} & \ab{0} 
\end{array} \right)
\left( 
\begin{array}{c}
\ab{u} \\ \ab{\lambda}
\end{array}
\right) = 
\left(
\begin{array}{c}
\ab{f} \\ \ab{c}
\end{array}
\right),
\end{equation}
where $\ab{K} \in \R^{n^2|\mathcal{I}_{\Omega}| \times n^2|\mathcal{I}_{\Omega}|} $ is the global stiffness matrix of the form
\[
\ab{K} = \Diag (\{ \ab{K}^{(i)}\}_{i \in \mathcal{I}_{\Omega}})
\]
with the ``patch-local'' stiffness matrices
\[
\ab{K}^{(i)} = [a^{(i)}(\phi^{(i)}_{\ab{j}_1},\phi^{(i)}_{\ab{j}_2})]_{\ab{j}_1 \in \bb{\mathcal{J}};\ab{j}_2 \in \bb{\mathcal{J}}},
\]
and $\ab{f}$ is the global load vector of the form
\[
\ab{f}= [\ab{f}^{(i)}]_{i \in \mathcal{I}_\Omega}
\]
with the``patch-local'' load vectors
\[
\ab{f}^{(i)} = [F^{(i)}(\phi^{(i)}_{\ab{j}})]_{\ab{j} \in \bb{\mathcal{J}}}.
\]

In the isogeometric formulation, the elements of the stiffness matrices $\ab{K}^{(i)}$ for the polyharmonic equation of order $m$ can be computed as 
\begin{equation*} 
a^{(i)}(\phi^{(i)}_{\ab{j}_1},\phi^{(i)}_{\ab{j}_2}) = 
\begin{cases}
\displaystyle \int_{[0,1]^2} \psi_{\ab{j}_1}^{(i,m)}(\ab{G}^{(i)}(\ab{\xi}))^T \left( N^{(i)}(\ab{\xi}) \, \psi_{\ab{j}_2}^{(i,m)}(\ab{G}^{(i)}(\ab{\xi}))\right) \, d\ab{\xi}, \qquad \mbox{$m$ is odd}, \\ 
\displaystyle \int_{[0,1]^2} | \det J \ab{G}^{(i)} (\ab{\xi})| \, \psi_{\ab{j}_1}^{(i,m)}(\ab{G}^{(i)}(\ab{\xi})) \,  \psi_{\ab{j}_2}^{(i,m)}(\ab{G}^{(i)}(\ab{\xi}))  \, d\ab{\xi},  \quad \mbox{$m$ is even},
\end{cases}
\end{equation*} 
where $\psi_{\ab{j}}^{(i,\ell)}(\ab{G}^{(i)}(\ab{\xi})) $ are vector ($\ell$ is odd) or scalar ($\ell$ is even)  functions, defined recursively as 
\begin{equation*} 
\psi_{\ab{j}}^{(i,\ell)}(\ab{G}^{(i)}(\ab{\xi}))  = 
\begin{cases}
\nabla \psi_{\ab{j}}^{(i,\ell-1)}(\ab{G}^{(i)}(\ab{\xi})) , \qquad \mbox{$\ell$ is odd},\\[0.3cm]
\displaystyle \frac{1}{| \det J \ab{G}^{(i)} (\ab{\xi})|} \nabla \cdot \left(N^{(i)}(\ab{\xi}) \, \psi_{\ab{j}}^{(i,\ell-1)}(\ab{G}^{(i)}(\ab{\xi}))   \right) , \quad \mbox{$\ell$ is even},\\[0.5cm]
\phi_{\ab{j}}^{(i)}(\ab{G}^{(i)}(\ab{\xi})) , \qquad \mbox{$\ell=0$},
\end{cases}
\end{equation*} 
with
\begin{equation*} 
N^{(i)}(\ab{\xi}) = | \det J \ab{G}^{(i)} (\ab{\xi})| \left( J \ab{G}^{(i)} (\ab{\xi})\right)^{-T} \left( J \ab{G}^{(i)} (\ab{\xi})\right)^{-1}, 
\end{equation*} 
and $J \ab{G}^{(i)}$ being the Jacobian matrix of geometry mapping $\ab{G}^{(i)}$.
Moreover, the elements of load vectors $\ab{f}^{(i)}$ can be computed as
\begin{equation*} 
F^{(i)}(\phi^{(i)}_{\ab{j}}) = \int_{[0,1]^2} f(\ab{G}^{(i)}(\ab{\xi})) \phi^{(i)}_{\ab{j}} (\ab{G}^{(i)}(\ab{\xi})) | \det J \ab{G}^{(i)} (\ab{\xi})| d \ab{\xi}.
\end{equation*} 

The idea of a IETI-based method is to solve instead of the large problem~\eqref{eq:large_problem} for the coefficients~$\ab{u}$ and the Lagrange multipliers~$\ab{\lambda}$ a smaller dual problem just for the Lagrange multipliers~$\ab{\lambda}$. If the matrix~$\ab{K}$ would be invertible, then the first row of~\eqref{eq:large_problem} would lead to 
\begin{equation} \label{eq:first_line}
\ab{u} = \ab{K}^{-1} \left(\ab{f}-\ab{C}^T \ab{\lambda}\right),
\end{equation}
and by inserting~\eqref{eq:first_line} into the second row
of~\eqref{eq:large_problem}, we would obtain simply the following dual problem: Find $\ab{\lambda} \in \R^q$ such that
\begin{equation} \label{eq:dual_problem}
\ab{C} \ab{K} \ab{C}^T \ab{\lambda} = \ab{C} \ab{K} \ab{f} - \ab{c}.
\end{equation}
However, the single matrices~$\ab{K}^{(i)}$, $i \in \mathcal{I}_{\Omega}$, and hence also the matrix~$\ab{K}$ are in general not invertible, cf.~\cite{FarhatFETI-DP, Farhat1991, KlPeSaJu12, SoglTakacs_IETI_Elasticity}). Therefore, we will introduce in the next subsection an adaption of the dual problem~\eqref{eq:dual_problem}, more precisely a dual-primal formulation, which will allow to solve the saddle point problem~\eqref{eq:large_problem} in an efficient way. 

\subsection{The general framework}
\label{sub:general}

The development of the general framework for solving the saddle point problem~\eqref{eq:large_problem} will be based on the introduction of a dual-primal formulation of the problem, cf.~\cite{FarhatFETI-DP, KlPeSaJu12, SoglTakacs_IETI_Elasticity}.
For this, we will start with the decomposition of the vectors of coefficients $\ab{u}^{(i)}$, $i \in \mathcal{I}_{\Omega}$, into two subvectors $\ab{u}_P^{(i)}$ and $\ab{u}_R^{(i)}$, of the load vectors~$\ab{f}^{(i)}$, $i \in \mathcal{I}_{\Omega}$, into two subvectors~$\ab{f}_P^{(i)}$ and $\ab{f}_R^{(i)}$, and consequently of the matrices $\ab{K}^{(i)}$, $i \in \mathcal{I}_{\Omega}$, into 
\begin{equation*} 
\ab{K}^{(i)} = \left( \begin{array}{cc}
\ab{K}_{PP}^{(i)} & {\ab{K}_{RP}^{(i)}}^T \\
{\ab{K}_{RP}^{(i)}} & \ab{K}_{RR}^{(i)}
\end{array} \right),
\end{equation*}
where $P$ and $R$ stands for primal and remaining degrees of freedom. Since the rank deficiency of the matrices $\ab{K}^{(i)}$ depends on the spline degree $\ab{p}$ and the order $m$ of the polyharmonic equation, as numerically shown for $m=2,3$ and $m \leq p \leq 9$ on several bilinear multi-patch domains, and can be up to $2p+m$, cf. \ref{sec:AppendixA},  
the subvectors $\ab{u}_P^{(i)}$ should be large enough to imply invertible matrices $\ab{K}_{RR}^{(i)}$. In order to further simplify the construction, the subvectors $\ab{u}_P^{(i)}$ will contain those coefficients $u^{(i)}_{\ab{j}}$, $\ab{j} \in \bb{\mathcal{J}}$, which correspond to the $s+1$ outer rings of the coefficients (i.e. the boundary ring and the first $s$ neighboring rings of the coefficients) associated with the patch~$\Omega^{(i)}$. 
All remaining coefficients associated with the patch~$\Omega^{(i)}$ will be collected in the subvectors $\ab{u}_R^{(i)}$. Therefrom, we can write the load vector~$\ab{f}$ as 
\[
\ab{f}=\left(\begin{array}{c}
\ab{f}_P\\
\ab{f}_R
\end{array}
\right)
\]
with $\ab{f}_P=[\ab{f}_P^{(i)}]_{i \in \mathcal{I}_\Omega}$ and $\ab{f}_R=[\ab{f}_R^{(i)}]_{i \in \mathcal{I}_\Omega}$, and the stiffness matrix $\ab{K}$ as
\begin{equation*} 
\ab{K}^{} = \left( \begin{array}{cc}
\ab{K}_{PP}^{} & \ab{K}_{RP}^{T} \\
\ab{K}_{RP} & \ab{K}_{RR}^{}
\end{array} \right)
\end{equation*}
with
$\ab{K}_{PP} = \Diag (\{ \ab{K}_{PP}^{(i)}\}_{i \in \mathcal{I}_{\Omega}})$,  
$\ab{K}_{RP} = \Diag (\{ \ab{K}_{RP}^{(i)}\}_{i \in \mathcal{I}_{\Omega}})$ and
$\ab{K}_{RR} = \Diag (\{ \ab{K}_{RR}^{(i)}\}_{i \in \mathcal{I}_{\Omega}})$.
Furthermore, we will decompose the matrix of conditions $\ab{C}$ into three matrices: $\ab{C}_B$, $\ab{C}_{\Xi}$, and $\ab{C}_{\Gamma}$. Matrix
$\ab{C}_B$ will be a partial permutation matrix representing all the imposed boundary conditions~\eqref{eq:constraintsB} as $\ab{C}_B \, \ab{u}_P = \ab{g}$, matrix $\ab{C}_{\Xi}$ will represent the {$C^{2s}$-}smoothness conditions~\eqref{eq:constraintsV} at all inner vertices $\bfm{\Xi}^{(i)}$, $ i \in \mathcal{I}_{\Xi}^I$, and matrix $\ab{C}_{\Gamma}$ the {$C^s$-}smoothness conditions~\eqref{eq:constraintsI1} and \eqref{eq:constraintsI2} at all inner edges~$\Gamma^{(i)}$, $i \in \mathcal{I}_\Gamma^{I}$. Note that all these conditions only affect the subvector $\ab{u}_P$, so the system \eqref{eq:large_problem} can be written as
\begin{equation} \label{eq:decomposed_problem}
\left( \begin{array}{ccccc}
\ab{K}_{PP} & \ab{K}_{RP}^T & \ab{C}_B^T & \ab{C}_{\Xi}^T & \ab{C}_{\Gamma}^T \\
\ab{K}_{RP} & \ab{K}_{RR} & \ab{0} & \ab{0} & \ab{0} \\
\ab{C}_B & \ab{0} & \ab{0} & \ab{0} & \ab{0} \\
\ab{C}_{\Xi} & \ab{0} & \ab{0} & \ab{0} & \ab{0} \\
\ab{C}_{\Gamma} & \ab{0} & \ab{0} & \ab{0} & \ab{0} 
\end{array} \right)
\left( 
\begin{array}{c}
\ab{u}_P \\ 
\ab{u}_R \\ 
\ab{\lambda}_B\\
\ab{\lambda}_{\Xi}\\
\ab{\lambda}_{\Gamma}
\end{array}
\right) = 
\left(
\begin{array}{c}
\ab{f}_P \\ 
\ab{f}_R \\
\ab{g} \\
\ab{0}\\
\ab{0}
\end{array}
\right).
\end{equation} 
Let us additionally split and reorder the coefficients of $\ab{u}_P$ into $\ab{u}_P = \left(\ab{u}_B^T, \, \ab{u}_F^T\right)^T$, where $\ab{u}_B$ collects the coefficients of $\ab{u}_P$ which are determined by the boundary conditions~\eqref{eq:constraintsB},  
and $\ab{u}_F$ are the remaining coefficients of $\ab{u}_{P}$. 
We can reorder coefficients in $\ab{u}_B$ in such a way that matrix $\ab{C}_B$ becomes $\ab{C}_B = \left( \begin{array}{cc}
\ab{I} & \ab{0}
\end{array} \right)$
 and consequently 
 $\ab{u}_B = \ab{g}$. 
Before explaining a possible way of computing $\ab{u}_B$ from boundary conditions, let us present the  simplification of our saddle point problem \eqref{eq:decomposed_problem}, once the values of $\ab{u}_B$ are known. 
By writing 
$$
\ab{K}_{PP} = \left( \begin{array}{cc}
\ab{K}_{BB}^{} & \ab{K}_{FB}^T \\
\ab{K}_{FB} & \ab{K}_{FF}^{}
\end{array} \right),\quad
\ab{K}_{RP} = \left( \begin{array}{cc}
\ab{K}_{RB} &
\ab{K}_{RF}
\end{array} \right),
$$
and
$$
\ab{C}_{\Xi} = 
\left( \begin{array}{cc}
\ab{0} & \ab{C}_{\Xi,F} 
\end{array} \right), \quad
\ab{C}_{\Gamma} = \left(\begin{array}{cc}
\ab{C}_{\Gamma,B} \;\, \ab{C}_{\Gamma,F} 
\end{array}\right), 
$$
and using some trivial manipulations, we get to the following equivalent saddle point problem
\begin{equation} \label{eq:decomposed_problem_NoCB}
\left( \begin{array}{cccc}
\overline{\ab{K}}_{PP} & \overline{\ab{K}}_{RP}^T &  {\ab{C}}_{\Xi}^T & \overline{\ab{C}}_{\Gamma}^T \\
\overline{\ab{K}}_{RP} & \ab{K}_{RR} & \ab{0} & \ab{0}  \\
{\ab{C}}_{\Xi} & \ab{0} & \ab{0} & \ab{0}  \\
\overline{\ab{C}}_{\Gamma} & \ab{0} & \ab{0} & \ab{0} 
\end{array} \right)
\left( 
\begin{array}{c}
\ab{u}_P \\ 
\ab{u}_R \\ 
\ab{\lambda}_{\Xi}\\
\ab{\lambda}_{\Gamma}
\end{array}
\right) = 
\left(
\begin{array}{c}
\overline{\ab{f}}_P \\ 
\overline{\ab{f}}_R \\
\ab{0}\\
\overline{\ab{g}}
\end{array}
\right),
\end{equation}
where
$$
\overline{\ab{K}}_{PP} = \left( \begin{array}{cc}
\ab{I} & \ab{0} \\
\ab{0} & \ab{K}_{FF}
\end{array} \right),\quad
\overline{\ab{K}}_{RP} = \left( \begin{array}{cc}
\ab{0} &
\ab{K}_{RF}
\end{array} \right),
\quad
\overline{\ab{C}}_{\Gamma} = \left( \begin{array}{cc}
\ab{0} & \ab{C}_{\Gamma,F} 
\end{array} \right),
$$
and
$$
\overline{\ab{f}}_{P} = \left( \begin{array}{c}
\ab{g} \\
\ab{f}_{F} - \ab{K}_{FB} \,\ab{g}
\end{array} \right),\quad
\overline{\ab{f}}_{R} = {\ab{f}}_{R} - \ab{K}_{RB} \,\ab{g},\quad
\overline{\ab{g}}= - \ab{C}_{\Gamma,B} \,\ab{g}.
$$
Below, it remains to show how to compute the coefficients $\ab{u}_B$ from the given boundary conditions \eqref{eq:constraintsB}, how to efficiently solve the saddle point problem \eqref{eq:decomposed_problem_NoCB}, and how to construct the matrices $\ab{C}_{\Xi}$ and $\ab{C}_{\Gamma}$. 

\paragraph{Computing the coefficients of $\ab{u}_B$}

We compute the coefficients that correspond to the subvector $\ab{u}_B$ from the boundary conditions \eqref{eq:constraintsB} by solving a IETI type $L^2$-fitting saddle point problem
\begin{equation} \label{eq:Boundary_IETI_problem}
\left( \begin{array}{cc}
\ab{M} & \ab{B}^T \\
\ab{B} & \ab{0} 
\end{array} \right)
\left( 
\begin{array}{c}
\ab{u}_B \\ \ab{\mu}
\end{array}
\right) = 
\left(
\begin{array}{c}
\ab{b} \\ \ab{0}
\end{array}
\right),
\end{equation}
where $\ab{M}$ is the mass matrix of the form
$
\ab{M} = \Diag (\{ \ab{M}^{(i)}\}_{i \in \mathcal{I}^{B}_{\Omega}})
$ 
with
\begin{align*}
\ab{M}^{(i)} & = \left[ \int_{\Omega^{(i)}}\phi^{(i)}_{\ab{j}_1} \phi^{(i)}_{\ab{j}_2} d\Omega^{(i)} \right]_{\ab{j}_1, \ab{j}_2 \in {\bb{\mathcal{J}}^{(i)}_{\partial}}}  \\
& = 
\left[ \int_{[0,1]^2} | \det J \ab{G}^{(i)} (\ab{\xi})| \, \phi_{\ab{j}_1}^{(i)}(\ab{G}^{(i)}(\ab{\xi})) \,  \phi_{\ab{j}_2}^{(i)}(\ab{G}^{(i)}(\ab{\xi}))  \, d\ab{\xi}
\right]_{\ab{j}_1, \ab{j}_2 \in {\bb{\mathcal{J}}^{(i)}_{\partial}}}
,
\end{align*}
where $\bb{\mathcal{J}}^{(i)}_{\partial} \subset \bb{\mathcal{J}}$ collects the indices~$\ab{j} \in \bb{\mathcal{J}}$ of the basis functions~$\phi_{\ab{j}}^{(i)}$ of the boundary patch~$\Omega^{(i)}$, $i \in \mathcal{I}_{\Omega}^B$, which correspond to the
$m$ outer rings of indices $\ab{j} \in \bb{\mathcal{J}}$ that are associated with $\partial \Omega$.
Moreover, $\ab{b}$ is the load vector of the form
$
\ab{b}= [\ab{b}^{(i)}]_{i \in \mathcal{I}^{B}_{\Omega}}
$
with 
$
\ab{b}^{(i)} = [F(\phi^{(i)}_{\ab{j}})]_{\ab{j} \in {\bb{\mathcal{J}}^{(i)}_{\partial}}}.
$
{The matrix~$\ab{B}$ is of the form 
\[
\ab{B}=[\ab{B}^{(i)}]_{i \in \mathcal{I}_{\Xi}^{B,\nu \geq 2}},
\]
where each submatrix~$\ab{B}^{(i)}$ represents in the rows for the boundary vertex~$\bfm{\Xi}^{(i)}$ with patch valency~$\nu_i \geq 2$, i.e.~$i \in \mathcal{I}_{\Xi}^{B,\nu \geq 2}$, certain smoothness conditions on the coefficient vector~$\ab{u}_{B}$, that are required to ensure that the resulting solution~$u_h$ can satisfy the $C^s$-smoothness conditions~\eqref{eq:constraintsI1} and \eqref{eq:constraintsI2} across the inner edges~$\overline{\Gamma^{(i_j)}}$ containing $\bfm{\Xi}^{(i)}$. 

Let $\bfm{\Xi}^{(i)}$, $i \in \mathcal{I}_{\Xi}^{B,\nu \geq 2}$, be a boundary vertex with a patch valency~$\nu_i \geq 2$, assuming that the neighboring patches~$\Omega^{(i_0)}, \ldots, \Omega^{(i_{\nu_i-1})}$ with $\cap_{j=0}^{\nu_i-1} \overline{\Omega^{(i_j)}}=\bfm{\Xi}^{(i)}$ are labeled in counterclockwise order, and that the corresponding inner edges~$\Gamma^{(i_j)}$, $j=0,\ldots, \nu_i-2$, are given as $\overline{\Gamma^{(i_j)}} = \overline{\Omega^{(i_j)}} \cap \overline{\Omega^{(i_{j+1})}}$. 
We compute now the matrix
~$\ab{B}^{(i)}$ as follows. We first construct a matrix~
$\widetilde{\ab{B}}^{(i)}$,
which collects in its rows for the inner edges~$\Gamma^{(i_j)}$, $j=0,\ldots,\nu_i-2$, assuming that the inner edge~$\overline{\Gamma^{(i_j)}}$ is parameterized as in~\eqref{eq:gC}, i.e. as $\ab{G}^{(i_j)}(0,\xi) = \ab{G}^{(i_{j+1})}(0,\xi)$, $\xi \in [0,1]$, and assuming that the boundary vertex~$\bfm{\Xi}^{(i)}$ is given as $\bfm{\Xi}^{(i)} = \ab{G}^{(i_j)}(0,0) = \ab{G}^{(i_{j+1})}(0,0)$, the smoothness conditions
\begin{align} \label{eq:boundary_vertex}
 \partial^w \g_\ell^{(i_j,i_j)}[u_h](0) & = \partial^w \g_\ell^{(i_j,i_{j+1})}[u_h](0), \quad
w=0, \ldots, m+s-1-\ell, 
\quad \ell=0, \ldots, s,
\end{align}
for $j=0,\ldots,\nu_i-2$,
with respect to the coefficient vector~$\left(\widetilde{\ab{u}}_F^T, \, \ab{u}_B^T\right)^T$, where $\widetilde{\ab{u}}_F$ collects the coefficients of $\ab{u}_F$ which are involved in the smoothness conditions~\eqref{eq:boundary_vertex}.
We now compute the null space 
$\ab{N}^{(i)}$ of matrix 
$\left(\overline{\ab{B}}^{(i)}\right)^T$, where 
$\overline{\ab{B}}^{(i)}$ denotes the matrix obtained from the matrix 
$\widetilde{\ab{B}}^{(i)}$
by putting all the columns which correspond to the coefficient vector ${\ab{u}}_B$ to zero. 
The matrix~
$\ab{B}^{(i)}$
is finally obtained as
\begin{equation}  \label{eq:boundary_vertexOutside}
\left( \begin{array}{cc}
\ab{0} & 
\ab{B}^{(i)}
\end{array} \right)  = 
(\ab{N}^{(i)})^T \widetilde{\ab{B}}^{(i)},
\end{equation}
where the zero matrix represents the obtained zero columns with respect to the coefficient vector~$\widetilde{\ab{u}}_F$.
This ensures that the matrix~$\ab{B}$  represents in its rows the desired smoothness conditions with respect to the coefficient vector~$\ab{u}_{B}$ and that the matrix has full row rank. 

Let us now solve the linear system \eqref{eq:Boundary_IETI_problem}. Note that matrix $\ab{M}$ is invertible and we can efficiently solve the system by first solving the symmetric positive definite system 
$ \left(\ab{B} \ab{M}^{-1} \ab{B}^T \right) \ab{\mu}$ $ = \ab{B} \ab{M}^{-1} \ab{b} $ 
for the Lagrange multipliers $\ab{\mu}$, and then by computing $\ab{u}_B$ as
$$
\ab{u}_B = \ab{M}^{-1} \left( \ab{b} - \ab{B}^T \ab{\mu}  \right).
$$

\paragraph{Efficient solving of the saddle point problem  \eqref{eq:decomposed_problem_NoCB}}
Since the stiffness (sub)matrix $\ab{K}_{RR}$ is invertible, the second row of \eqref{eq:decomposed_problem_NoCB} implies 
\begin{equation} \label{eq:solution_u_R}
\ab{u}_R = \ab{K}_{RR}^{-1} \left( 
\overline{\ab{f}}_R - \overline{\ab{K}}_{RP} {\ab{u}}_P
\right).
\end{equation}
Note that $\ab{K}_{RR}$ is block-diagonal, so computing $\ab{K}_{RR}^{-1} = \Diag (\{ (\ab{K}_{RR}^{(i)})^{-1}\}_{i \in \mathcal{I}_{\Omega}})$ corresponds to solving local problems on each patch $\Omega^{(i)}$ by
(sparse) LU factorization.
Inserting \eqref{eq:solution_u_R} in the remaining row equations of \eqref{eq:decomposed_problem_NoCB} implies 
\begin{equation} \label{eq:decomposed_problem_Schur}
\left( \begin{array}{ccc}
\ab{S}_P  &  {\ab{C}}_{\Xi}^T & \overline{\ab{C}}_{\Gamma}^T \\
{\ab{C}}_{\Xi} & \ab{0} & \ab{0}  \\
\overline{\ab{C}}_{\Gamma} & \ab{0} & \ab{0} 
\end{array} \right)
\left( 
\begin{array}{c}
\ab{u}_P \\  
\ab{\lambda}_{\Xi}\\
\ab{\lambda}_{\Gamma}
\end{array}
\right) = 
\left(
\begin{array}{c}
\overline{\ab{f}}_P - \overline{\ab{K}}_{RP}^T \ab{K}_{RR}^{-1} \overline{\ab{f}}_R \\
\ab{0}\\
\overline{\ab{g}}
\end{array}
\right),
\end{equation}
where
\begin{equation}
\label{eq:Schur_complement}
\ab{S}_P := \overline{\ab{K}}_{PP} - \overline{\ab{K}}_{RP}^T \ab{K}_{RR}^{-1} \overline{\ab{K}}_{RP} = 
\left( \begin{array}{cc}
\ab{I} & \ab{0} \\
\ab{0} & \ab{S}_F
\end{array} \right)\;\,\mbox{ with }\;\,\ab{S}_{F} := \ab{K}_{FF} -{\ab{K}}_{RF}^T \ab{K}_{RR}^{-1} {\ab{K}}_{RF} 
\end{equation}
is the Schur complement of the block $\ab{K}_{RR}$ in matrix $\ab{K}$ having again block structure, since 
$
\ab{S}_{F} = \Diag (\{ \ab{S}_{F}^{(i)}\}_{i \in \mathcal{I}_{\Omega}}).
$
However, $\ab{S}_P$ is only symmetric positive semidefinite and therefore not invertible in general, so we can not express $\ab{u}_P$ directly as we could do in the case of $\ab{u}_R$ in \eqref{eq:solution_u_R}.
Therefore we rewrite \eqref{eq:decomposed_problem_Schur} as
\begin{equation} \label{eq:decomposed_problem_Schur_withS}
\left( \begin{array}{cc}
\ab{S} &  \widetilde{\ab{C}}_{\Gamma}^T \\
\widetilde{\ab{C}}_{\Gamma} & \ab{0} 
\end{array} \right)
\left( 
\begin{array}{c}
\ab{u}_{P,\Xi} \\  
\ab{\lambda}_{\Gamma}
\end{array}
\right) = 
\left(
\begin{array}{c}
\widetilde{\ab{f}}_P\\
\overline{\ab{g}}
\end{array}
\right),
\end{equation}
where
\begin{equation}
\label{eq:decomposed_problem_Schur_withS_Notation}
\ab{u}_{P,\Xi} = \left(\begin{array}{c}
\ab{u}_P\\
\ab{\lambda}_{\Xi}
\end{array}
\right), \;\;
\ab{S} = \left( \begin{array}{cc}
\ab{S}_P &  {\ab{C}}_{\Xi}^T \\
{\ab{C}}_{\Xi} & \ab{0} 
\end{array} \right), \;\;
\widetilde{\ab{C}}_{\Gamma} = \left( \begin{array}{cc}
\overline{\ab{C}}_{\Gamma} \; \, \ab{0} \end{array}\right), \;\; 
\widetilde{\ab{f}}_P = \left(\begin{array}{c}
\overline{\ab{f}}_P - \overline{\ab{K}}_{RP}^T \ab{K}_{RR}^{-1} \overline{\ab{f}}_R \\
\ab{0}
\end{array}
\right).
\end{equation}
Since $\ab{S}_P$ is symmetric positive semidefinite, ${\ab{C}}_{\Xi}$ has full row rank and $\ker(\ab{S}_P) \cap \ker({\ab{C}}_{\Xi}) = \{\ab{0}\}$, the matrix $\ab{S}$ is invertible. Note that $\ab{S}^{-1}$ can again be computed efficiently in a parallel way, see~\ref{sec:AppendixB}. This now implies 
\begin{equation*} \label{eq:solution_u_P}
\ab{u}_{P,\Xi} = \ab{S}^{-1} \left( 
\widetilde{\ab{f}}_P - \widetilde{\ab{C}}_{\Gamma}^T \ab{\lambda}_{\Gamma}
\right).
\end{equation*}
We are now left with the symmetric positive definite (cf.~\ref{sec:AppendixC}) linear system for $\ab{\lambda}_{\Gamma}$ of the form
\begin{equation} \label{eq:solution_LambdaGamma}
\left( \widetilde{\ab{C}}_{\Gamma} \ab{S}^{-1} \widetilde{\ab{C}}_{\Gamma}^T \right) \ab{\lambda}_{\Gamma} = \widetilde{\ab{C}}_{\Gamma} \ab{S}^{-1} \widetilde{\ab{f}}_P -\overline{\ab{g}},
\end{equation}
which can be solved by using e.g.~the conjugate gradient method.

\begin{rem}
We could also reverse the role of the matrices $\ab{C}_{\Gamma}$ and $\ab{C}_{\Xi}$ to change \eqref{eq:decomposed_problem_Schur_withS} and \eqref{eq:decomposed_problem_Schur_withS_Notation} into
\begin{equation*} 
\left( \begin{array}{cc}
\ab{S} &  \widetilde{\ab{C}}_{\Xi}^T \\
\widetilde{\ab{C}}_{\Xi} & \ab{0} 
\end{array} \right)
\left( 
\begin{array}{c}
\ab{u}_{P,\Gamma} \\  
\ab{\lambda}_{\Xi}
\end{array}
\right) = 
\left(
\begin{array}{c}
\widetilde{\ab{f}}_P\\
\ab{0}
\end{array}
\right)
\end{equation*}
with
\begin{equation*}
\ab{u}_{P,\Gamma} = \left(\begin{array}{c}
\ab{u}_P\\
\ab{\lambda}_{\Gamma}
\end{array}
\right), \;\;
\ab{S} = \left( \begin{array}{cc}
\ab{S}_P &  \overline{{\ab{C}}}_{\Gamma}^T \\
\overline{{\ab{C}}}_{\Gamma} & \ab{0} 
\end{array} \right), \;\;
\widetilde{\ab{C}}_{\Xi} = \left( \begin{array}{cc}
{\ab{C}}_{\Xi} \; \, \ab{0} \end{array}\right), \;\; 
\widetilde{\ab{f}}_P = \left(\begin{array}{c}
\overline{\ab{f}}_P - \overline{\ab{K}}_{RP}^T \ab{K}_{RR}^{-1} \overline{\ab{f}}_R \\
\overline{\ab{g}}
\end{array}
\right).
\end{equation*}
Then, we would first solve for $\ab{\lambda}_{\Xi}$ the symmetric positive definite linear system
\begin{equation*}
\left( \widetilde{\ab{C}}_{\Xi} \ab{S}^{-1} \widetilde{\ab{C}}_{\Xi}^T \right) \ab{\lambda}_{\Xi} = \widetilde{\ab{C}}_{\Xi} \ab{S}^{-1} \widetilde{\ab{f}}_P,
\end{equation*}
and then compute $\ab{u}_{P,\Gamma}$ via
\begin{equation*} 
\ab{u}_{P,\Gamma} = \ab{S}^{-1} \left( 
\widetilde{\ab{f}}_P - \widetilde{\ab{C}}_{\Xi}^T \ab{\lambda}_{\Xi}
\right).
\end{equation*}
\end{rem}

\paragraph{Constructing matrix $\ab{C}_{\Xi}$}

Recall that the matrix $\ab{C}_{\Xi}$ is of the form $\ab{C}_{\Xi} =  \left( \begin{array}{cc}
\ab{0}\;\; \ab{C}_{\Xi,F} 
\end{array} \right)$. The matrix~$\ab{C}_{\Xi,F}$ is further given by
\begin{equation} \label{eq:C_XiF}
\ab{C}_{\Xi,F}=[\ab{C}^{(i)}_{\Xi,F}]_{i \in \mathcal{I}_{\Xi}^{I}},
\end{equation}
where each submatrix~$\ab{C}^{(i)}_{\Xi,F}$ represents in the rows for the inner vertex~$\bfm{\Xi}^{(i)}$, $i \in \mathcal{I}_{\Xi}^{I}$, the $C^{2s}$-smoothness conditions~\eqref{eq:constraintsV} at the inner vertex~$\bfm{\Xi}^{(i)}$ with respect to the coefficient vector~$\ab{u}_{F}$. Let $\bfm{\Xi}^{(i)}$, $i \in \mathcal{I}_{\Xi}^{I}$, be an inner vertex assuming that the neighboring patches~$\Omega^{(i_0)}, \ldots, \Omega^{(i_{\nu_i-1})}$ with $\bfm{\Xi}^{(i)}=\cap_{\ell=0}^{\nu_i-1} \overline{\Omega^{(i_{\ell})}}$ are labeled in counterclockwise order. Then, the matrix~$\ab{C}^{(i)}_{\Xi,F}$ collects in its rows the $(\nu_i-1) \binom{2s+2}{2}$
linearly independent $C^{2s}$-smoothness conditions at 
the inner vertex $\ab{\Xi}^{(i)}$ given by
\begin{equation} \label{eq:inner_vertex}
  \left(\ab{u}^{(i_{\ell})}\right)^T \partial_1^{j_1} \partial_2^{j_2} \ab{\phi}^{(i_\ell)}(\ab{\Xi}^{(i)})  =  
\left(\ab{u}^{(i_{\ell+1})}\right)^T  \partial_1^{j_1} \partial_2^{j_2} \ab{\phi}^{(i_{\ell+1})}(\ab{\Xi}^{(i)}) , \quad j_1 + j_2 \leq 2s, \quad  \ell=0,\ldots,\nu_i-2.
\end{equation}
Note that for $\ell=1,\ldots,\nu_i-1$, 
the $\ell$-th block row of the matrix~$\ab{C}^{(i)}_{\Xi,F}$ 
consists of two nonzero matrices
\begin{equation}  \label{eq:C_XiF_blocks}
\ab{C}^{(i,i_{\ell-1})}_{\Xi,F} \in \R^{\binom{2s+2}{2} \times |\ab{u}_F^{(i_{\ell-1})}|} \quad {\rm and} \quad 
- \ab{C}^{(i,i_{\ell})}_{\Xi,F} \in \R^{\binom{2s+2}{2} \times |\ab{u}_F^{(i_\ell)}|} 
\end{equation}
corresponding to patches $\Omega^{(i_{\ell-1})}$ and $\Omega^{(i_{\ell})}$.
In addition, note that for the conditions~\eqref{eq:inner_vertex} also the geometry mappings~$\ab{G}^{(i_\ell)}$, $\ell=0,\ldots, \nu_i-1$, are involved.

\paragraph{Constructing matrix $\ab{C}_{\Gamma}$}

It remains to explain the construction of the matrix $\ab{C}_{\Gamma}$ which possesses the form 
\[
[\ab{C}_{\Gamma}^{(i)}]_{i \in \mathcal{I}_{\Gamma}^{I}},
\]
where each submatrix $\ab{C}_{\Gamma}^{(i)}$ represents in the rows for the inner edge~$\Gamma^{(i)}$, $i \in \mathcal{I}_\Gamma^I$, 
the $C^s$-smoothness conditions~\eqref{eq:constraintsI1} and \eqref{eq:constraintsI2} with respect to the coefficient vector~$\ab{u}_{P}$, which have to be linearly independent with the smoothness conditions \eqref{eq:constraintsV} and/or \eqref{eq:constraintsB} at the two endpoints of $\Gamma^{(i)}$, which have been already part of the matrices $\ab{C}_{\Xi}$ and/or $\ab{C}_B$. 
Let $\Gamma^{(i)}$, $i \in \mathcal{I}_{\Gamma}^{I}$, be an inner edge assuming that the two neighboring patches~$\Omega^{(i_0)}$ and $\Omega^{(i_1)}$, $i_0,i_1 \in \mathcal{I}_{\Omega}$, with
$\overline{\Gamma^{(i)}} = \overline{\Omega^{(i_0)}} \cap  \overline{\Omega^{(i_1)}}$, 
are parameterized as in 
\eqref{eq:standard_par}, and let $\xi^{p,r}_j$, $j=0,\ldots,n-1$, denote the Greville abscissae for the space $\mathcal{S}_h^{p,r}([0,1])$. {To construct the corresponding matrix~$\ab{C}_{\Gamma}^{(i)}$, we have to distinguish two cases.}

Assume first that $\overline{\Gamma^{(i)}} \cap \partial \Omega = \emptyset$, i.e.~both end points of $\overline{\Gamma^{(i)}}$ are inner vertices. Then these linearly independent smoothness conditions~\eqref{eq:constraintsI1} and \eqref{eq:constraintsI2} can be expressed as
\begin{equation*}
\label{eq:interfaceSmoothnessConditions_innerVertex}
    \begin{array}{ll}
    \g_\ell^{(i,\LL)}[u_h](\xi^{p,r}_{j_\ell}) = \gC^{(i)}_\ell[u_h](\xi^{p,r}_{j_\ell}), &  j_\ell = 
    0,\ldots,n-1
    \\
    \g_\ell^{(i,\RR)}[u_h](\xi^{p,r}_{j_\ell}) = \gC^{(i)}_\ell[u_h](\xi^{p,r}_{j_\ell}), & 
    j_\ell = 2s+1-\ell, \ldots, n-(2s-\ell) 
    \end{array},
    \quad \ell=0,1,\ldots, s,
\end{equation*}
where the spline functions $\gC_\ell^{(i)}[u_h] \in \mathcal{S}_h^{p-\ell,r+s-\ell}([0,1])$, $\ell=0,\ldots,s$, possess the spline representations
\[
\gC_\ell^{(i)}[u_h](\xi) = \sum_{j=0}^{n_{\ell}-1} d_{j,\ell}^{(i)} N_{j}^{p-\ell,r+s-\ell} 
\]
with the coefficients~$d_{j,\ell}^{(i)}$ and $n_{\ell}= \dim \mathcal{S}_h^{p-\ell,r+s-\ell}([0,1]) = p+1-\ell+k(p-r-s)$.  Additionally, let us collect these coefficients as $\ab{d}_{\ell}^{(i)} = \big[d_{j,\ell}^{(i)}\big]_{j=0}^{n_{\ell}-1}$ and $\ab{d}^{(i)} = \big[\ab{d}_\ell^{(i)}\big]_{\ell=0}^s$.
 In this way, we obtain linearly independent conditions which involve coefficients $\ab{u}_P^{(\LL)}$ and $\ab{u}_P^{(\RR)}$ for the patches $\Omega^{(\LL)}$ and $\Omega^{(\RR)}$ 
and coefficients~$\ab{d}^{(i)}$ of the spline functions $\gC_\ell^{(i)}[u_h]$,  
$\ell=0,\ldots,s$.
We can collect all the smoothness conditions for the inner edge $\Gamma^{(i)}$ in matrix form as
\begin{equation}  \label{eq:ui0ui1di_coefficients}
\left( \begin{array}{ccc}
\ab{C}^{(i,\LL)} &  \ab{0} & \ab{D}^{(i,\LL)} \\
\ab{0} & \ab{C}^{(i,\RR)} & \ab{D}^{(i,\RR)}
\end{array} \right)
\left( 
\begin{array}{c}
\ab{u}_P^{(\LL)} \\  
\ab{u}_P^{(\RR)}\\
\ab{d}^{(i)}
\end{array}
\right) = 
\left(
\begin{array}{c}
\ab{0}\\
\ab{0}
\end{array}
\right).    
\end{equation}
By denoting 
$
\ab{u}_P^{(\LL,\RR)} = \left(\begin{array}{c} \ab{u}_P^{(\LL)}\\ \ab{u}_P^{(\RR)} \end{array}\right)
$
and appropriately reordering conditions in \eqref{eq:ui0ui1di_coefficients}, we get a linear system 
\begin{equation}  
\label{eq:systemC1C2D1D2}
\left( \begin{array}{cc}
\ab{C}^{(i)}_{1} &  \ab{D}^{(i)}_1 \\
\ab{C}^{(i)}_2 & \ab{D}^{(i)}_2
\end{array} \right)
\left( 
\begin{array}{c}
\ab{u}_P^{(\LL,\RR)} \\  
\ab{d}^{(i)}
\end{array}
\right) = 
\left(
\begin{array}{c}
\ab{0}\\
\ab{0}
\end{array}
\right),    
\end{equation} 
where $\ab{D}^{(i)}_2$ is a square and invertible matrix. From the second row 
of \eqref{eq:systemC1C2D1D2} we now get
$\ab{d}^{(i)} = - (\ab{D}^{(i)}_2)^{-1} \ab{C}^{(i)}_2 \, \ab{u}_P^{(\LL,\RR)}$, which inserted into the first row 
of \eqref{eq:systemC1C2D1D2} implies conditions only on coefficients $\ab{u}_P^{(\LL)}$ and $\ab{u}_P^{(\RR)}$ as
$$
\left(\ab{C}_1^{(i)} -\ab{D}_1^{(i)} (\ab{D}^{(i)}_2)^{-1} \ab{C}^{(i)}_2 \right) \ab{u}_P^{(\LL,\RR)} = \ab{0}.
$$
Matrix $\ab{C}_1^{(i)} -\ab{D}_1^{(i)} (\ab{D}^{(i)}_2)^{-1} \ab{C}_2^{(i)}$, extended with zero columns at the positions corresponding to coefficients in $\ab{u}_P$ which are not in $\ab{u}_P^{(\LL)}$ and $\ab{u}_P^{(\RR)}$, represents now the matrix $\ab{C}^{(i)}_\Gamma$.

Assume now that one (or both) end points of $\overline{\Gamma^{(i)}}$ are boundary vertices of $\Omega$. 
Let us 
consider in detail the case that only one of the vertices, say $\ab{\Xi}^{(j)}$, $j \in \mathcal{I}_{\Xi}^B$, is a boundary vertex and that it equals $\ab{\Xi}^{(j)}=\ab{G}^{(i_0)}(0,0) = \ab{G}^{(i_1)}(0,0)$. The other cases would work similarly. {We construct first an initial matrix~$\ab{C}_{\Gamma}^{(i)}$ as described in the previous paragraph for the case that both end points of the edge would be inner vertices of~$\Omega$.} {Note that} there might be several edges $\Gamma^{(\ell)}$, $\ell \in \mathcal{I}_\Gamma^I$, having $\ab{\Xi}^{(j)}$ as one of their endpoints, i.~e.~$\nu_j >2$. Only in the case that $\Gamma^{(i)}$ is the edge with the smallest index of these edges $\Gamma^{(\ell)}$, we 
add now smoothness conditions representing the direct complement of conditions \eqref{eq:boundary_vertexOutside} with respect to conditions \eqref{eq:boundary_vertex} to the {matrix~$\ab{C}^{(i)}_{\Gamma}$}.  
{In case} that $\Gamma^{(i)}$ would not be the edge with the smallest index of these edges $\Gamma^{(\ell)}$, then these additional smoothness {conditions} have to be added to the {matrix~$\ab{C}_{\Gamma}^{(\ell)}$} with the smallest index $\ell$.

One possible strategy now to find this direct complement would be to take {the} conditions~\eqref{eq:boundary_vertex} starting with the condition having $w=m+s-1$ for each of the involved edges $\Gamma^{(\ell)}$, then the two conditions with $w=m+s-2$
for each of the edges $\Gamma^{(\ell)}$, etc., 
and adding these conditions one by one to the set of conditions \eqref{eq:boundary_vertexOutside} and keeping each newly added one only in the case if it is linearly independent with conditions \eqref{eq:boundary_vertexOutside} together with the previously added ones. We continue with 
this procedure until we get many enough (i.e.~the number of conditions \eqref{eq:boundary_vertex} minus the number of conditions \eqref{eq:boundary_vertexOutside}) added conditions, which now represent the direct complement to conditions \eqref{eq:boundary_vertexOutside} with respect to conditions \eqref{eq:boundary_vertex}. These are the required additional smoothness conditions that we have to add (extended with zeros at the positions corresponding
to the remaining coefficients in $\ab{u}_P$) to the matrix $\ab{C}_\Gamma^{(i)}$.

\section{Numerical examples} \label{section_numerical_examples}

In this section we will present some illustrative examples for solving the polyharmonic equation of order $m=2$ (the biharmonic equation) and $m=3$ (the triharmonic equation) over 
bilinear-like $G^{s}$ multi-patch geometries, $s\geq m-1$. To demonstrate the potential of our proposed IETI-based method, we will analyze the convergence behavior under $h$-refinement for the numerical solution $u_h$ obtained by solving the saddle point problem \eqref{eq:large_problem}.  
In all the examples below, the functions $f$ and $g_\ell$, $\ell=0,\ldots,m-1$, which define the right-side of the polyharmonic equation~\eqref{eq:polyharmonic} as well the boundary conditions, will be derived from the exact solution 
\begin{equation}  \label{eq:exactSolution}
 u(x_1,x_2)= \cos\left({x_1}\right) \sin\left({x_2}\right). 
\end{equation}
The resulting approximations $u_h$ will be compared with the exact solution \eqref{eq:exactSolution} by computing the relative errors with respect to the 
equivalent of the $H^m$-seminorm. Specifically, for the biharmonic ($m=2$) and triharmonic ($m=3$) equations, we will compute the relative errors in the $H^2$-seminorm and $H^3$-seminorm, respectively, as follows: 
\begin{equation} \label{eq:eqiuv2seminorms}
 \frac{\| \Delta u- \Delta u_h\|_{L^2}}{\| \Delta u \|_{L^2}} \quad {\rm and} \quad
\frac{ \| \nabla (\Delta u) - \nabla(\Delta u_h)\|_{L^2}}{\| \nabla(\Delta u )\|_{L^2}}.
\end{equation}
For brevity, we will refer to these equivalents of the seminorms as $H^2$ and $H^3$-seminorm, respectively. 

We will consider four different planar multi-patch domains (Domain~A--D) in our examples:
\begin{itemize}
\item \textbf{Domain~A:} This is the bilinear 11-patch domain~$\overline{\Omega}$, shown in Fig.~\ref{fig:domains} (top left), representing a hammer model. The parameterizations of the individual patches $\overline{\Omega^{(i)}}$ are defined by the vertices provided in the figure.
\item \textbf{Domain~B:} This is the bilinear 12-patch domain~~$\overline{\Omega}$, visualized in Fig. 2 (top right). Analogous to Domain A, the parameterizations of the individual patches $\overline{\Omega^{(i)}}$ are specified by the vertices shown in the figure.
\item \textbf{Domain~C:} This is the bilinear-like $G^1$ thirteen-patch controller domain~$\overline{\Omega}$, presented in Fig.~\ref{fig:domains} (bottom left). The individual patches are parameterized by bicubic geometry mappings from the space $\mathcal{S}_{1/3}^{\ab{3}, \ab{1}}([0,1]^2) \times \mathcal{S}_{1/3}^{\ab{3}, \ab{1}}([0,1]^2)$, as described in \cite{FaKaKoVi24}.
\item \textbf{Domain~D:} This is the bilinear-like $G^2$ five-patch Screw knob star domain~$\overline{\Omega}$, shown in Fig.~\ref{fig:domains} (bottom right). 
The individual patches are parameterized by biquintic geometry mappings  
from the space $\mathcal{S}_{1/4}^{\ab{5}, \ab{4}}([0,1]^2) \times \mathcal{S}_{1/4}^{\ab{5}, \ab{4}}([0,1]^2)$, as detailed in \cite{KaKoVi24c}. 
\end{itemize}
\begin{figure}[htb!]
    \centering
    \begin{tabular}{cc}
  \hspace{-0.8cm} \includegraphics[scale=0.14]{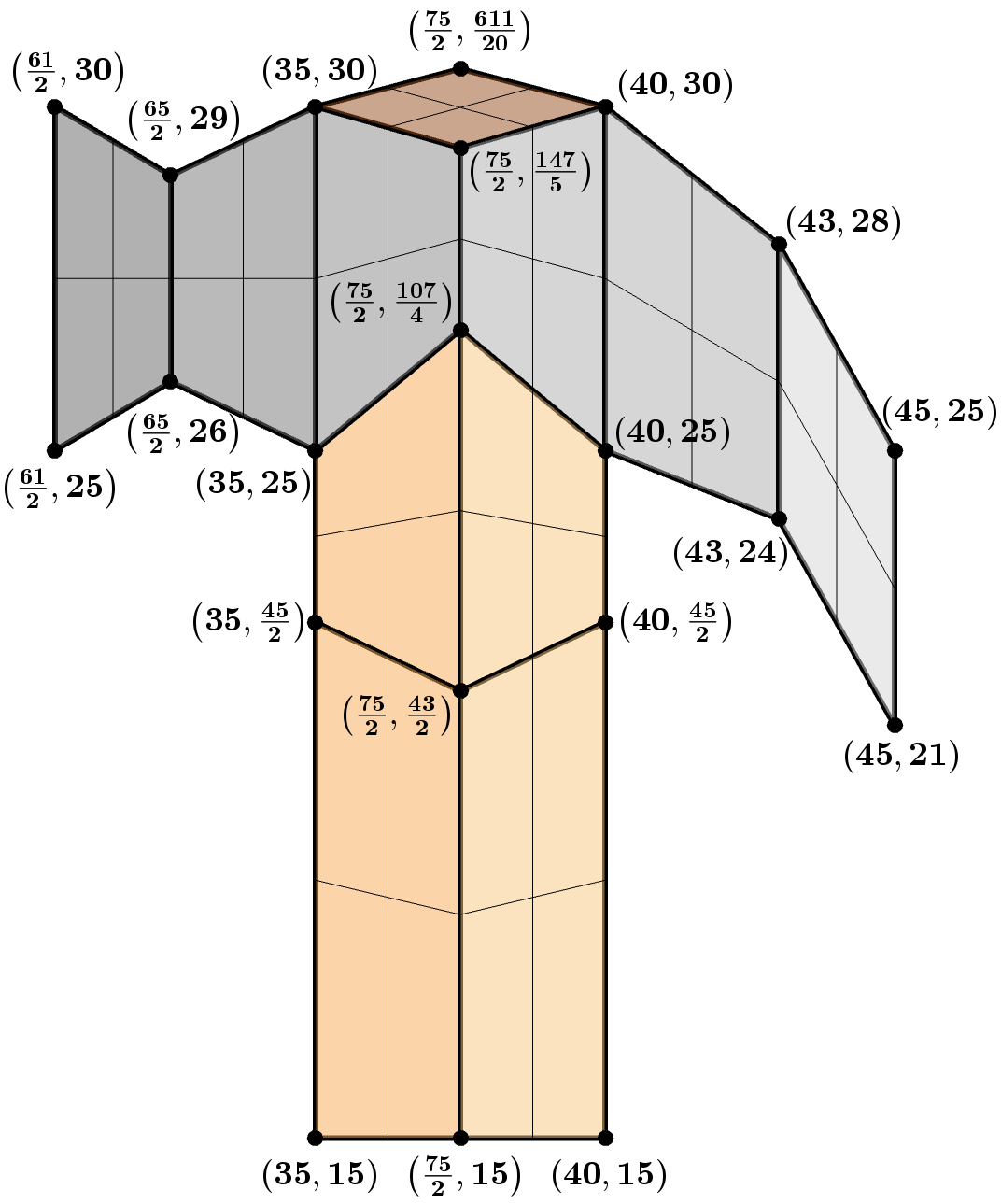}   & \hspace{-1.3cm}
  \includegraphics[scale=0.14]{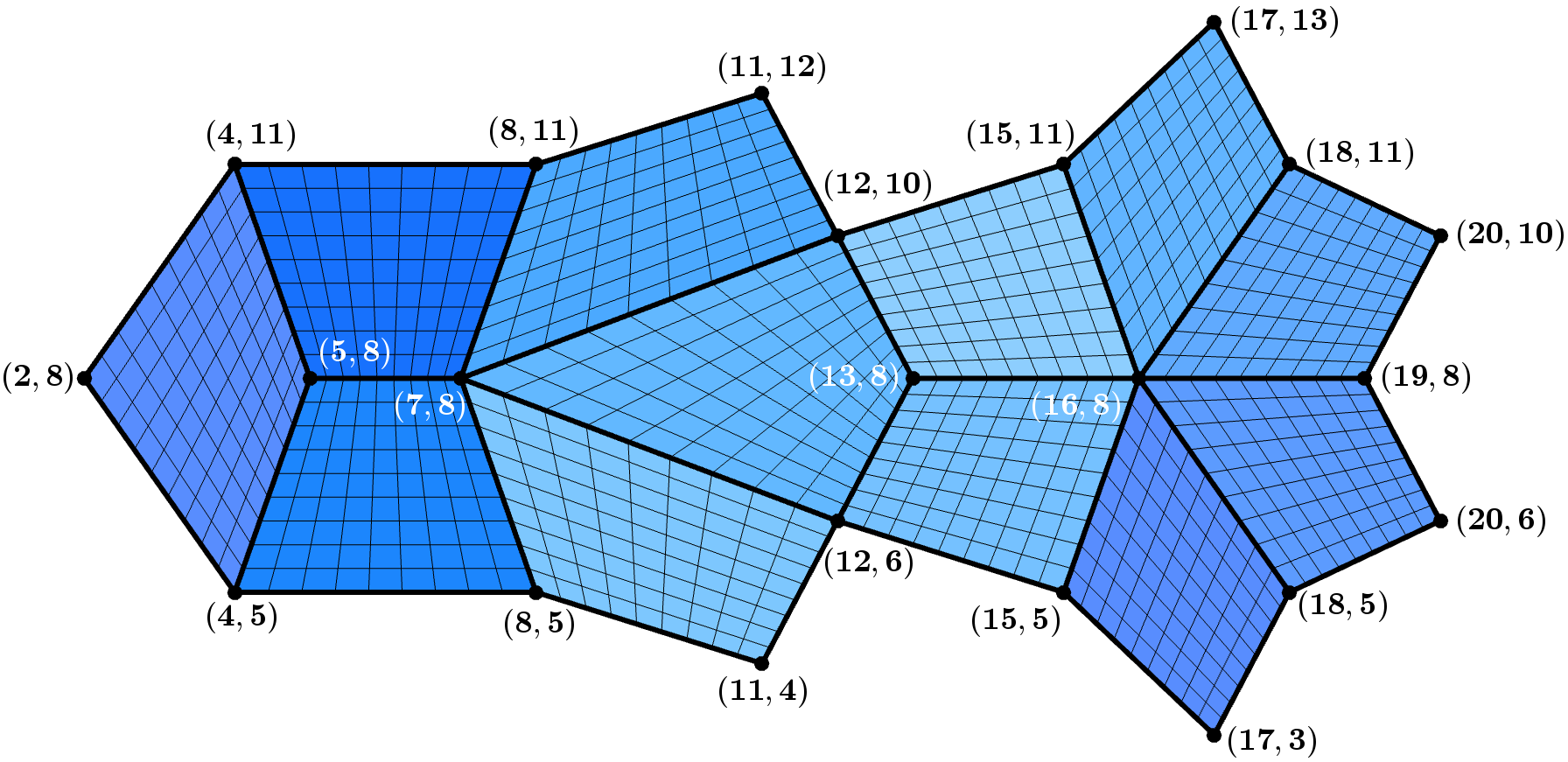} \\[0.2cm]
\hspace{-1.1cm} Domain A  & Domain B \\[0.8cm]
  \includegraphics[scale=0.14]{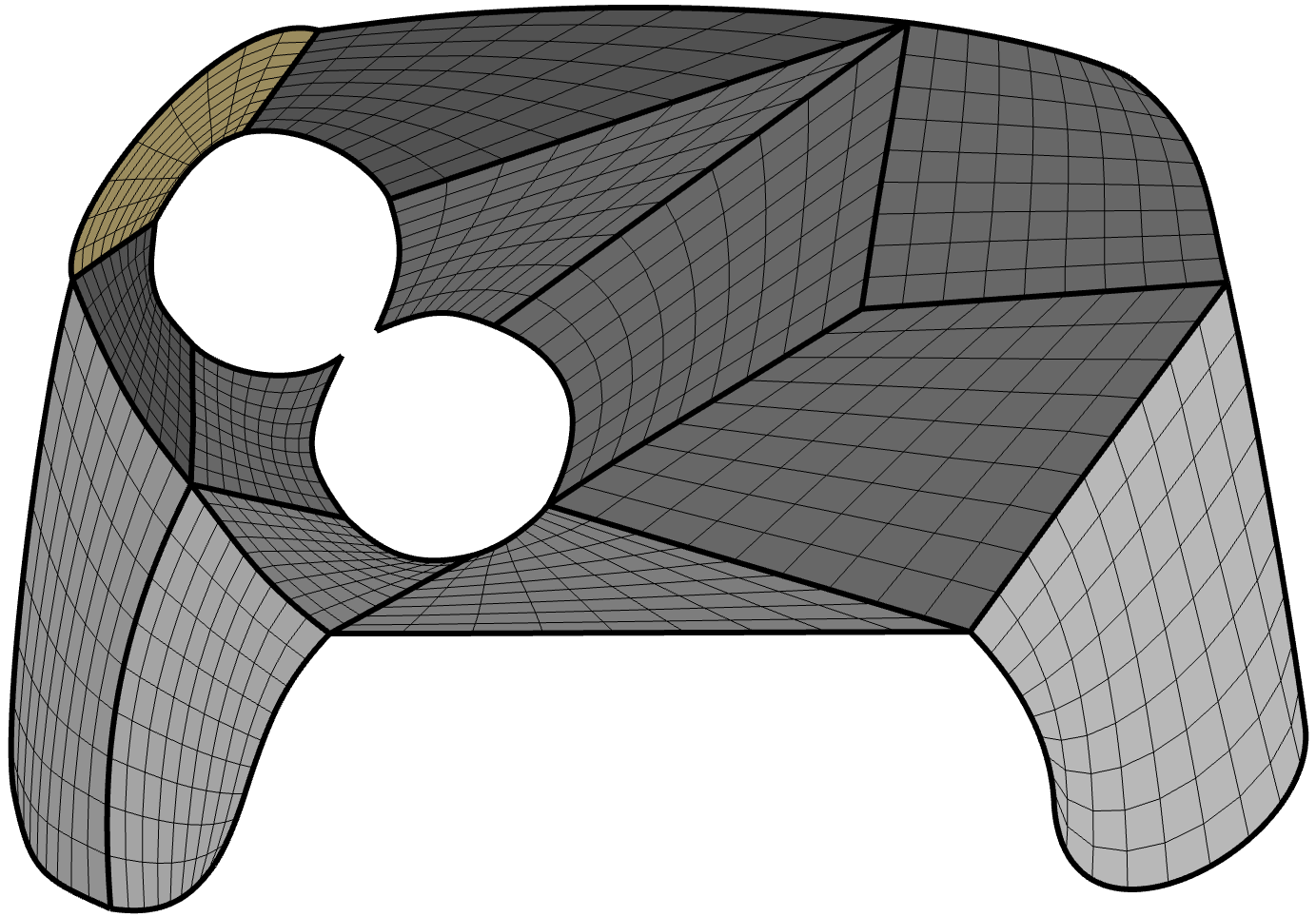}  \hskip3em &
  \includegraphics[scale=0.14]{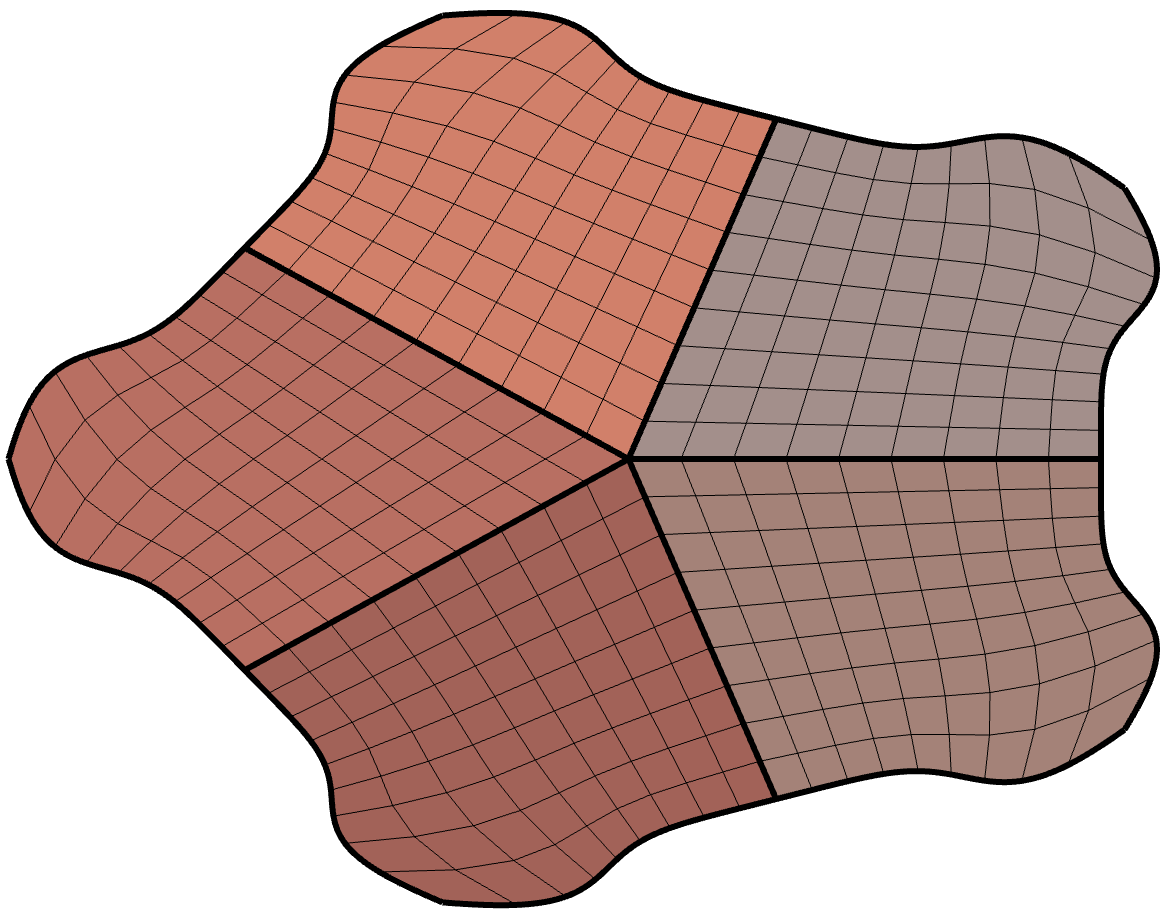}\\[0.4cm]
Domain C  & $\quad$ Domain D \\[0.4cm]
  \end{tabular}
\caption{The two bilinear multi-patch domains (Domain A (top left) and Domain B (top right)), the bilinear-like $G^1$ thirteen-patch controller domain (Domain C (bottom left)) and the bilinear-like $G^2$ five-patch Screw knob star domain (Domain D (bottom right)).}
    \label{fig:domains}
\end{figure}

Below, we will present now the numerical examples in which we will choose $s=m-1$.

\begin{ex} \label{ex:Biharmonic}
We begin by solving the polyharmonic equation of order $m=2$ (the biharmonic equation) over the two bilinear multi-patch domains (Domain A and B) presented in Fig.~\ref{fig:domains} (top row), and over the bilinear-like $G^1$ thirteen-patch controller domain, shown in Fig.~\ref{fig:domains} (bottom left). The computations are performed by using the underlying spline space $\mathcal{S}_h^{\ab{3},\ab{1}}([0,1]^2)$ with 
mesh sizes $h= h_{0}/2^j$, $j=0,\ldots,3$, where $h_{0}=1/5$ in case of Domain A and B, and $h_0=1/6$ in case of Domain~C. We study the resulting relative errors~\eqref{eq:eqiuv2seminorms} with respect to the $H^2$-seminorm, and observe for all three domains an {optimal} convergence of order~$\mathcal{O}(h^{2})$, cf. Fig.~\ref{fig:biharmonic} (right) for Domain A and B, and Fig.~\ref{fig:biharmonicBL} (right) for Domain C.
\begin{figure}[h!]
    \centering
     \includegraphics[scale=0.14]{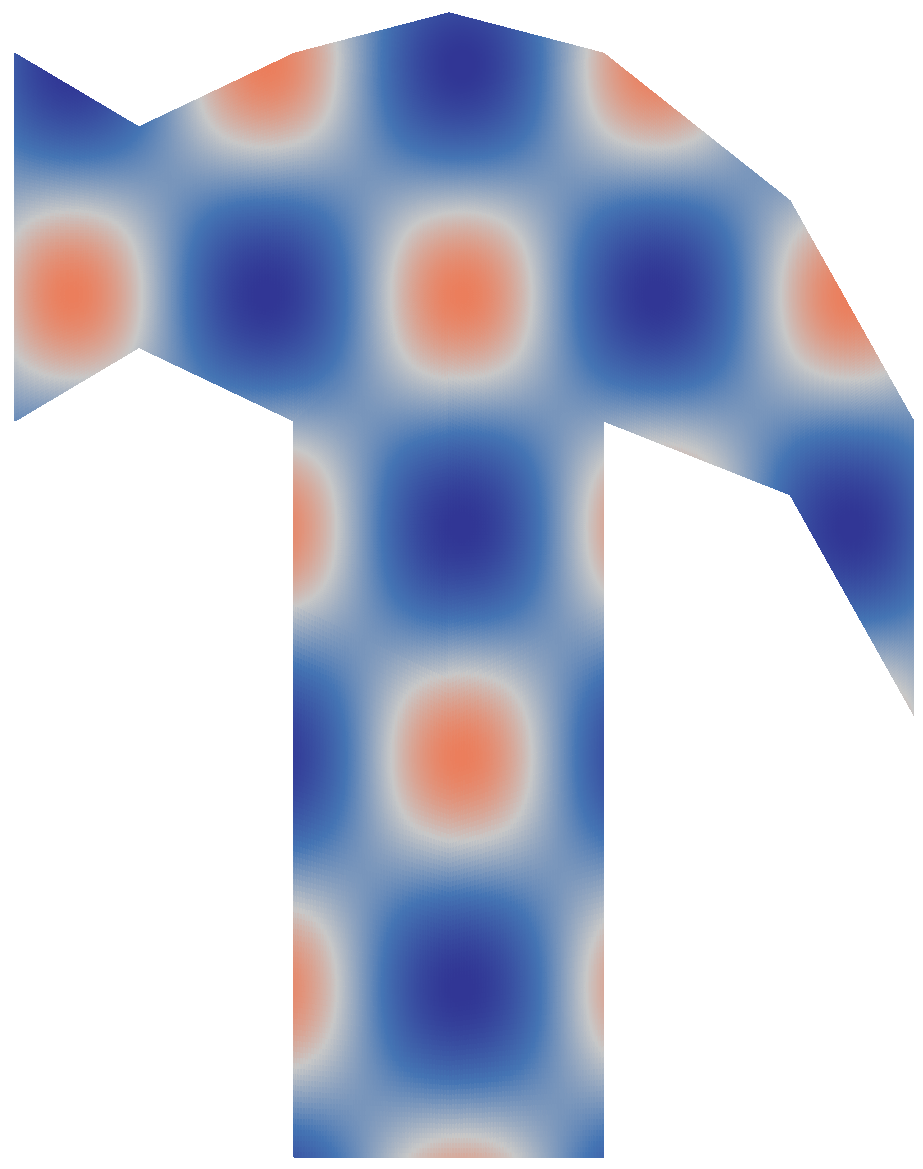} \hskip4.5em
  \includegraphics[scale=0.14]{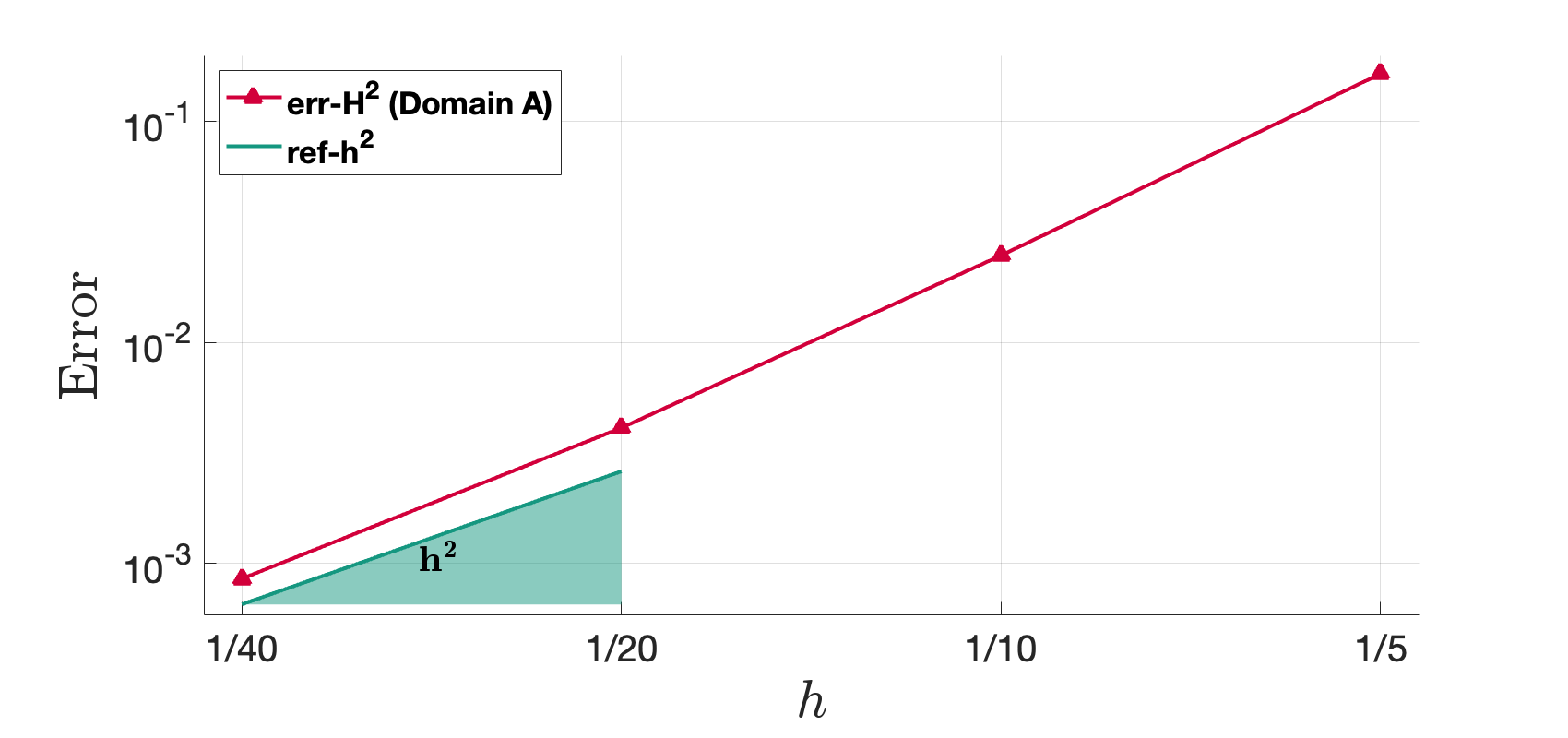}
  \vskip1.5em
  \includegraphics[scale=0.12]{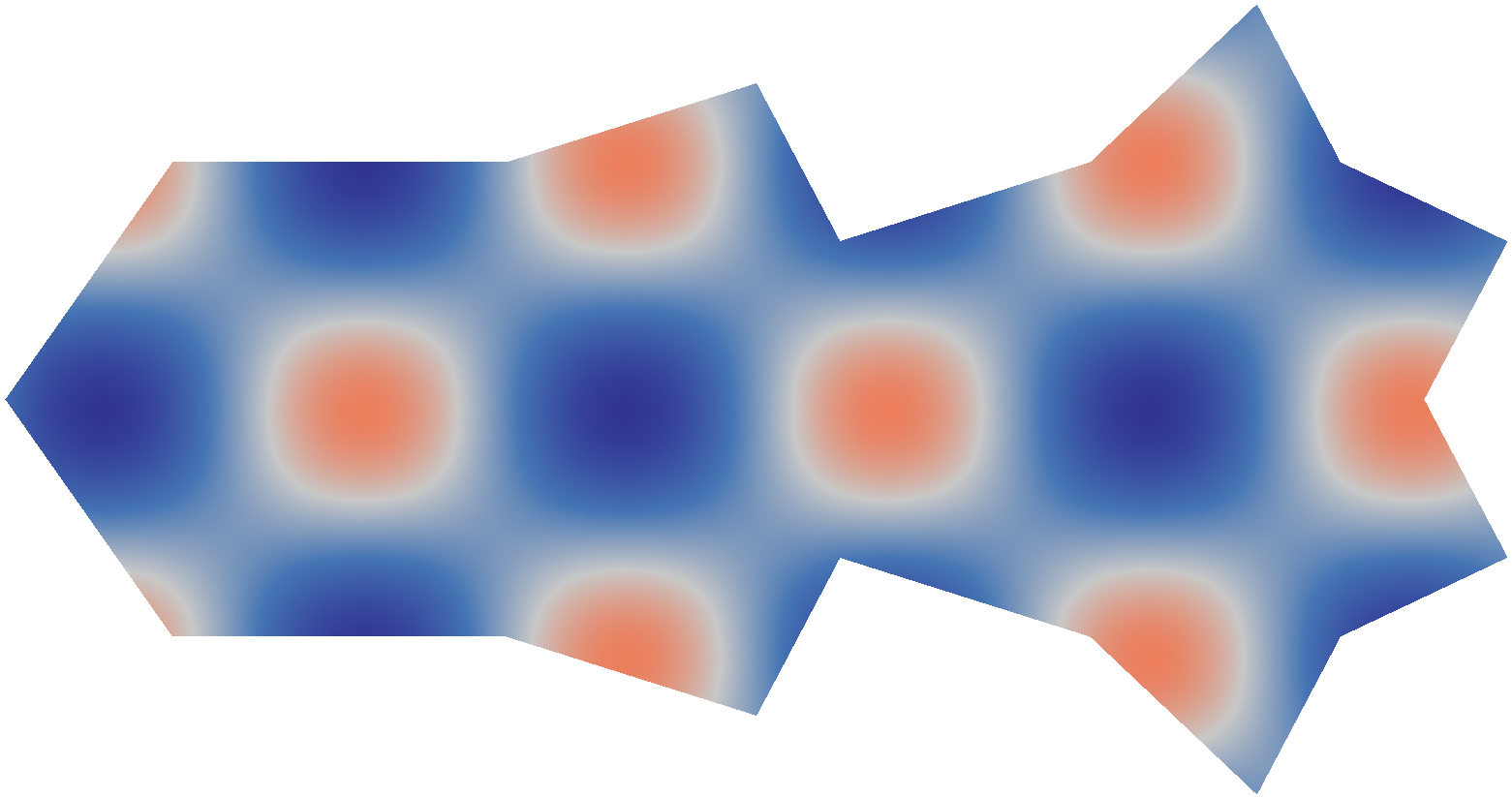} \hskip0.em
  \includegraphics[scale=0.14]{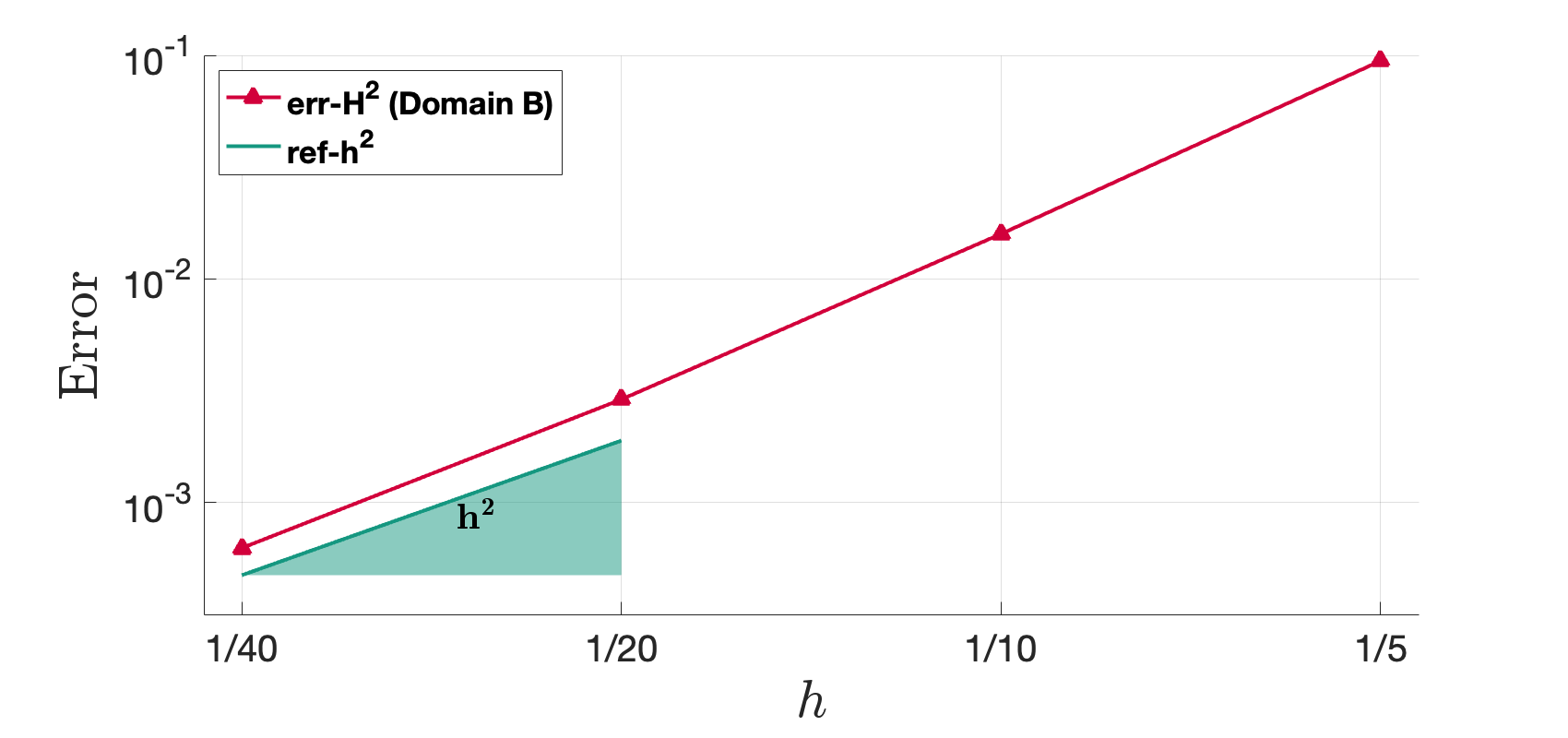}
\caption{Example~\ref{ex:Biharmonic}. The considered exact solution~\eqref{eq:exactSolution} (left) over the two bilinear multi-patch domains (Domain A and B) and the resulting relative errors~\eqref{eq:eqiuv2seminorms} with respect to the $H^2$-seminorm (right) for solving the biharmonic equation. The first row shows the results for Domain A and the second one for Domain B.}
\label{fig:biharmonic}
\end{figure}

\begin{figure}[h!]
    \centering
    \includegraphics[scale=0.13]{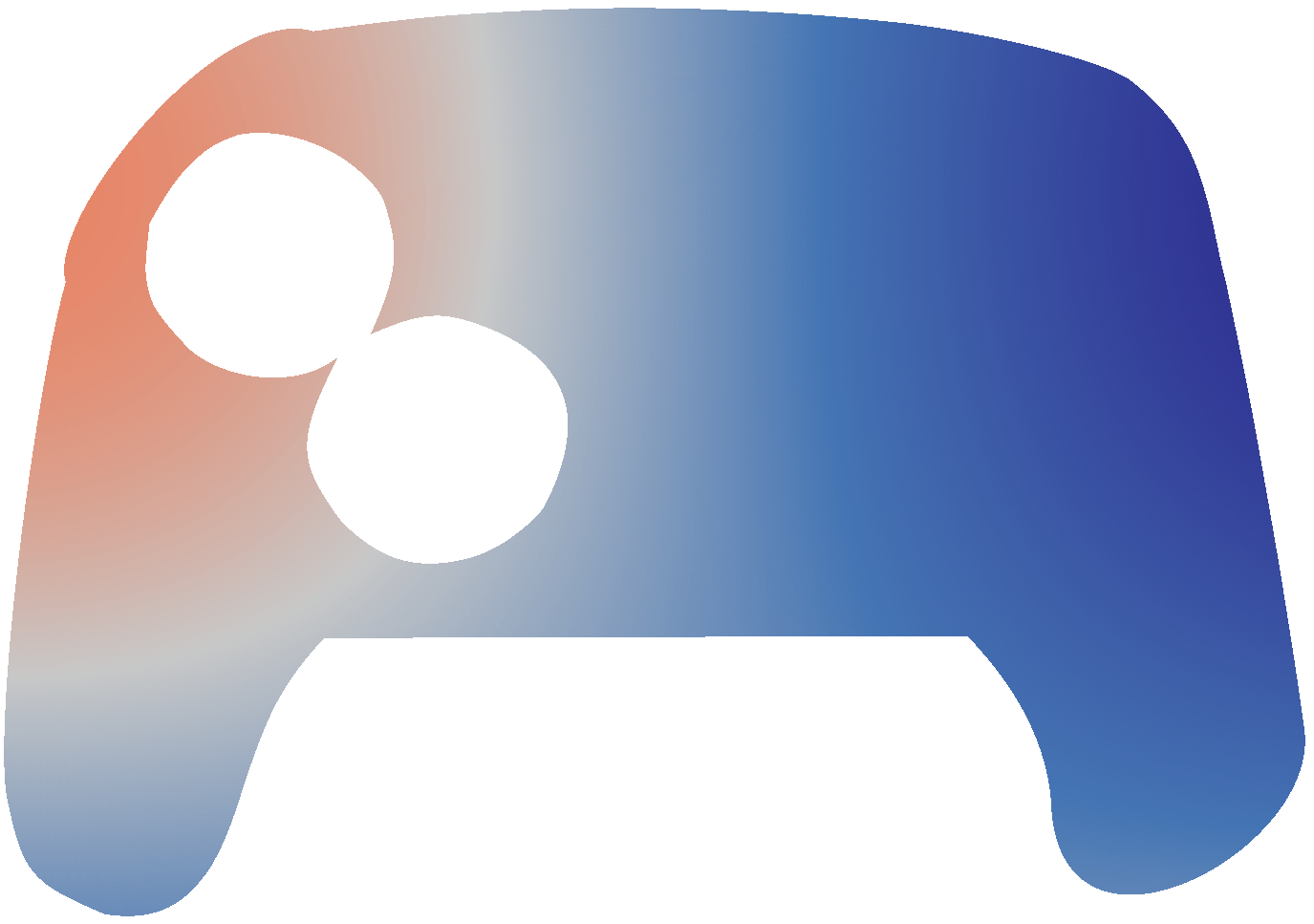}
    \hskip1em
    \includegraphics[scale=0.14]{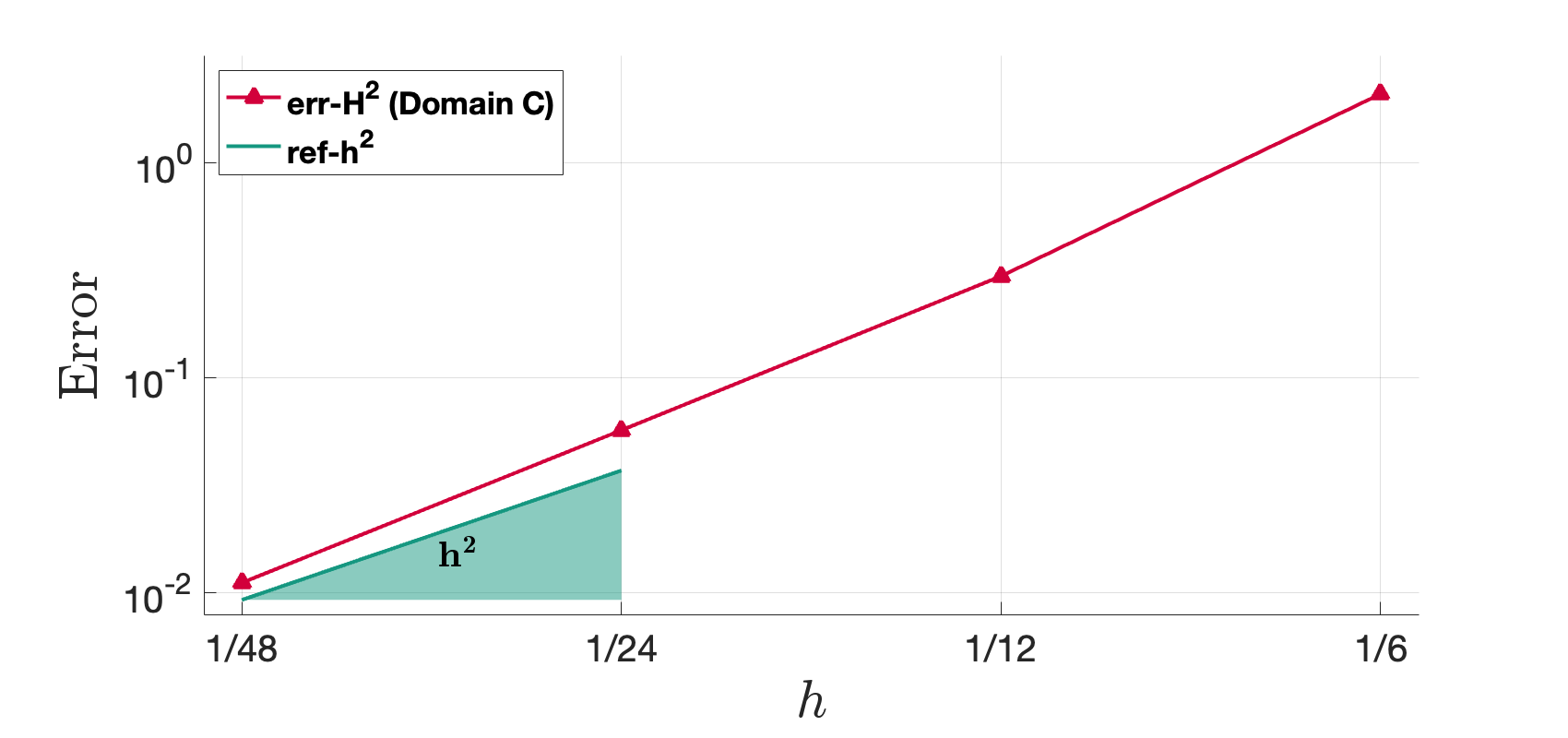}
\caption{Example~\ref{ex:Biharmonic}. The considered exact solution~\eqref{eq:exactSolution} (left) and the resulting relative errors~\eqref{eq:eqiuv2seminorms} with respect to the $H^2$-seminorm (right) for solving the biharmonic equation over the bilinear-like $G^1$ thirteen-patch Domain C.}
\label{fig:biharmonicBL}
\end{figure}
\end{ex}

\begin{ex} \label{ex:Triharmonic}
We continue by solving the polyharmonic equation of order $m=3$ (the triharmonic equation) on the one hand again over the two bilinear multi-patch domains (Domain A and Domain B) presented in Fig.~\ref{fig:domains} (top row), and on the other hand instead of over Domain C, which is just a bilinear-like $G^1$ multi-patch domain, now over the bilinear-like $G^2$ five-patch Screw knob star domain shown in Fig.~\ref{fig:domains} (bottom right).  To apply our method, we use the underlying spline space $\mathcal{S}_h^{\ab{5},\ab{2}}([0,1]^2)$ with mesh sizes $h= h_{0}/2^j$, $j=0,\ldots,3$, where $h_{0}=1/5$ in case of Domain~A and B, and $h_{0}=1/8$ in case of Domain~D. By analyzing the relative errors~\eqref{eq:eqiuv2seminorms} with respect to the $H^3$-seminorm, we observe for all three multi-patch domains an optimal convergence of order~$\mathcal{O}(h^{3})$, cf. Fig.~\ref{fig:triharmonic} (right) for Domain A and B, and Fig.~\ref{fig:triharmonicBL} (right) for Domain D.
\begin{figure}[h!]
    \centering
      \includegraphics[scale=0.14]{pics/solutionA.png} \hskip5em
  \includegraphics[scale=0.14]{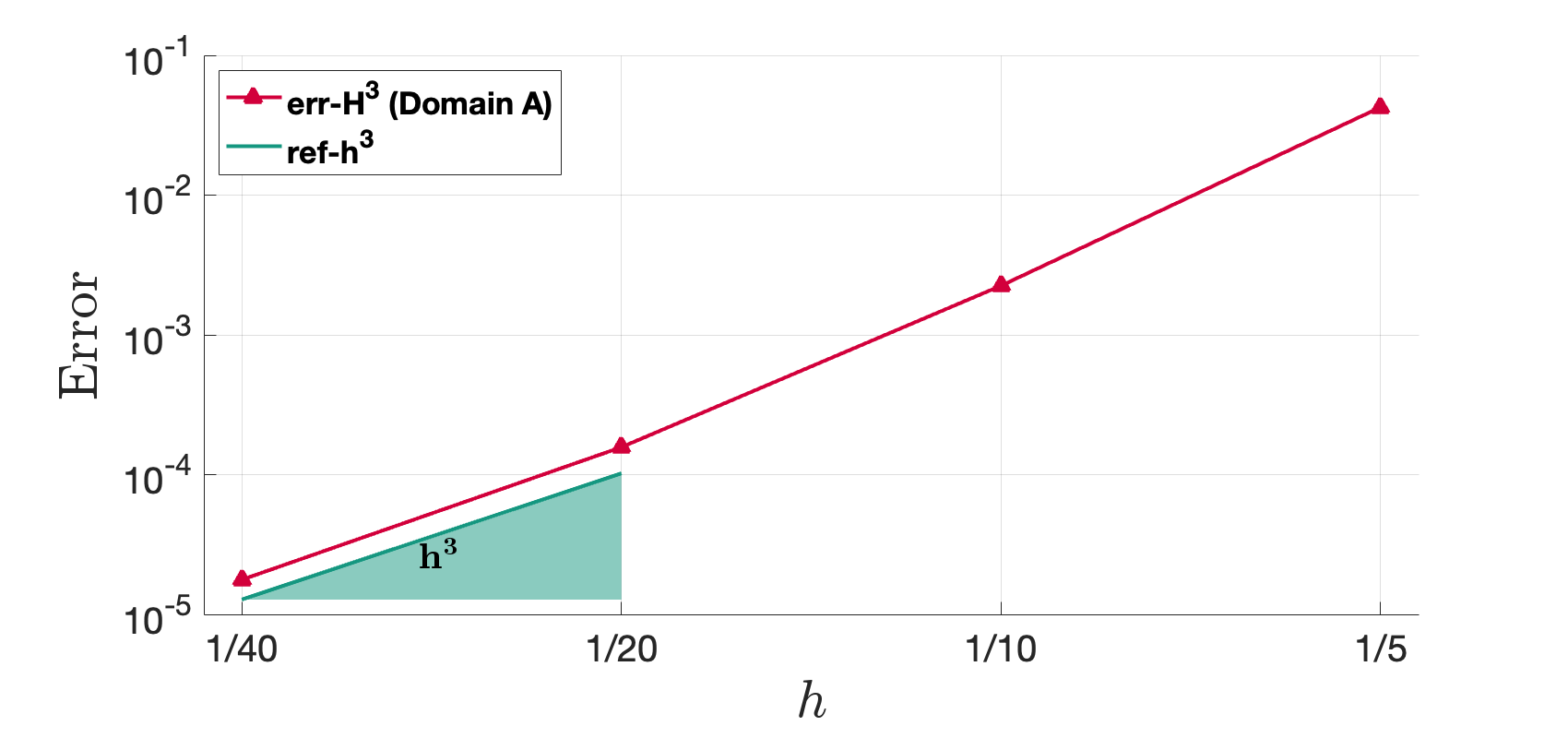}
  \vskip1.5em
  \includegraphics[scale=0.12]{pics/solutionB.png} \hskip0.em
  \includegraphics[scale=0.14]{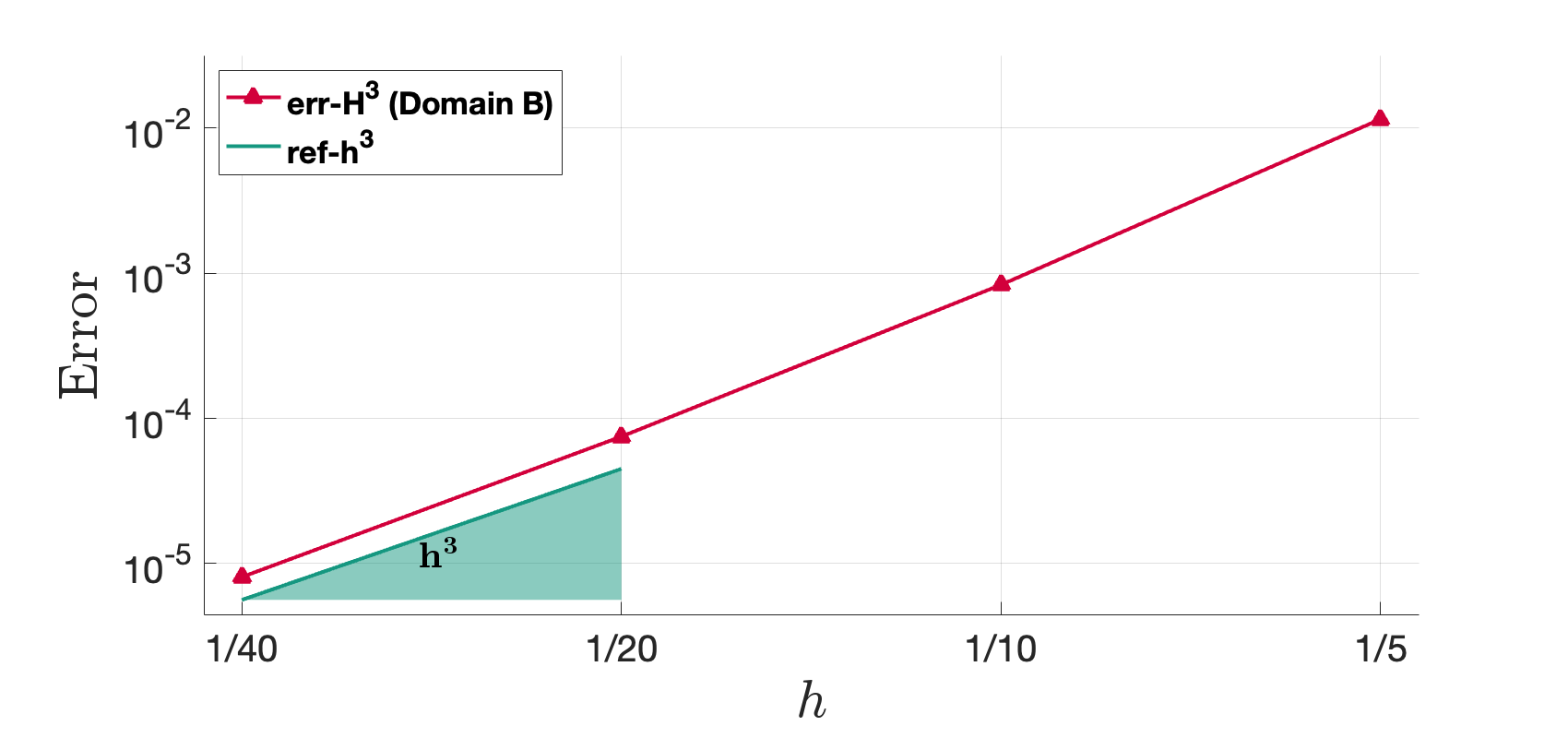}
\caption{Example~\ref{ex:Triharmonic}.
 The considered exact solution~\eqref{eq:exactSolution} (left) over the two bilinear multi-patch domains (Domain A and B) and the resulting relative errors~\eqref{eq:eqiuv2seminorms} with respect to the $H^3$-seminorm (right) for solving the triharmonic equation. The first row shows the results for Domain A and the second one for Domain B. 
 }
    \label{fig:triharmonic}
\end{figure}
\begin{figure}[h!]
    \centering
    \includegraphics[scale=0.13]{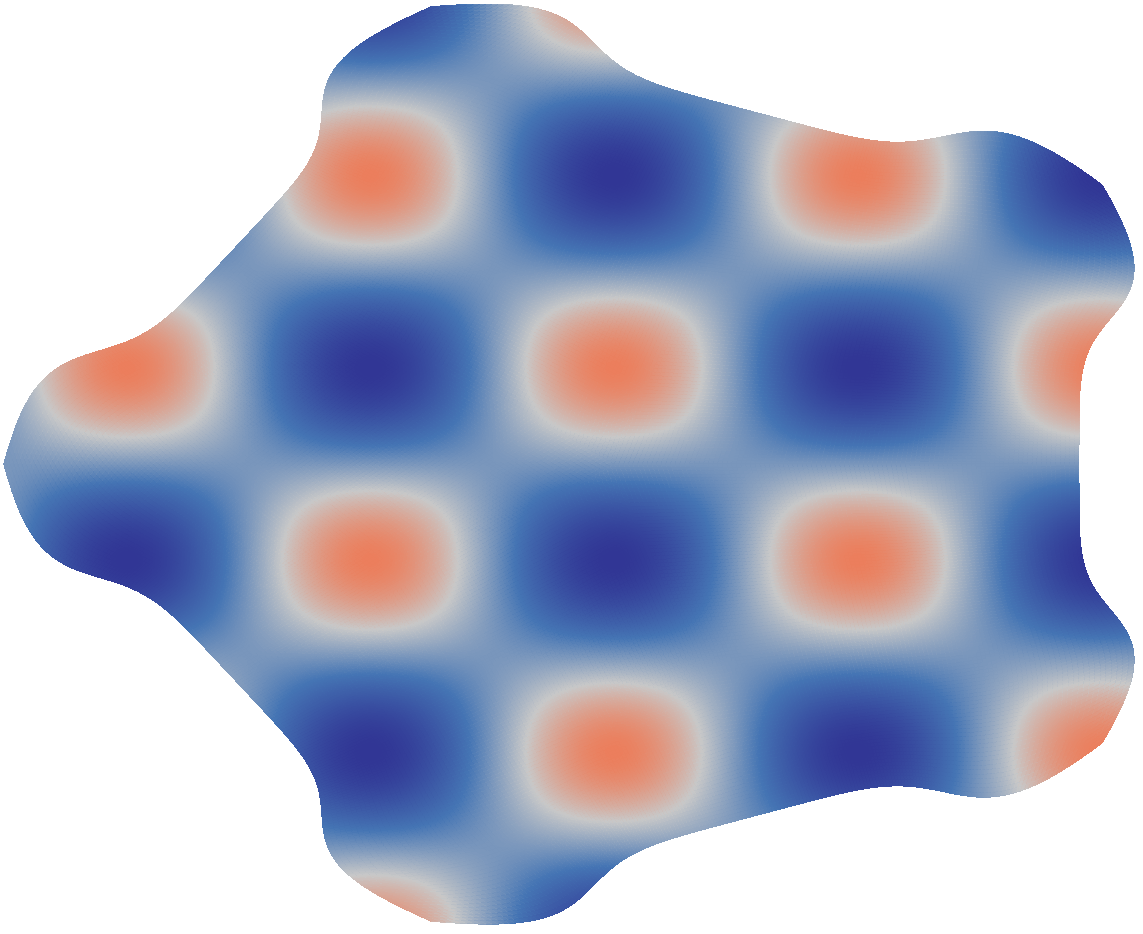}
    \hskip3em
    \includegraphics[scale=0.14]{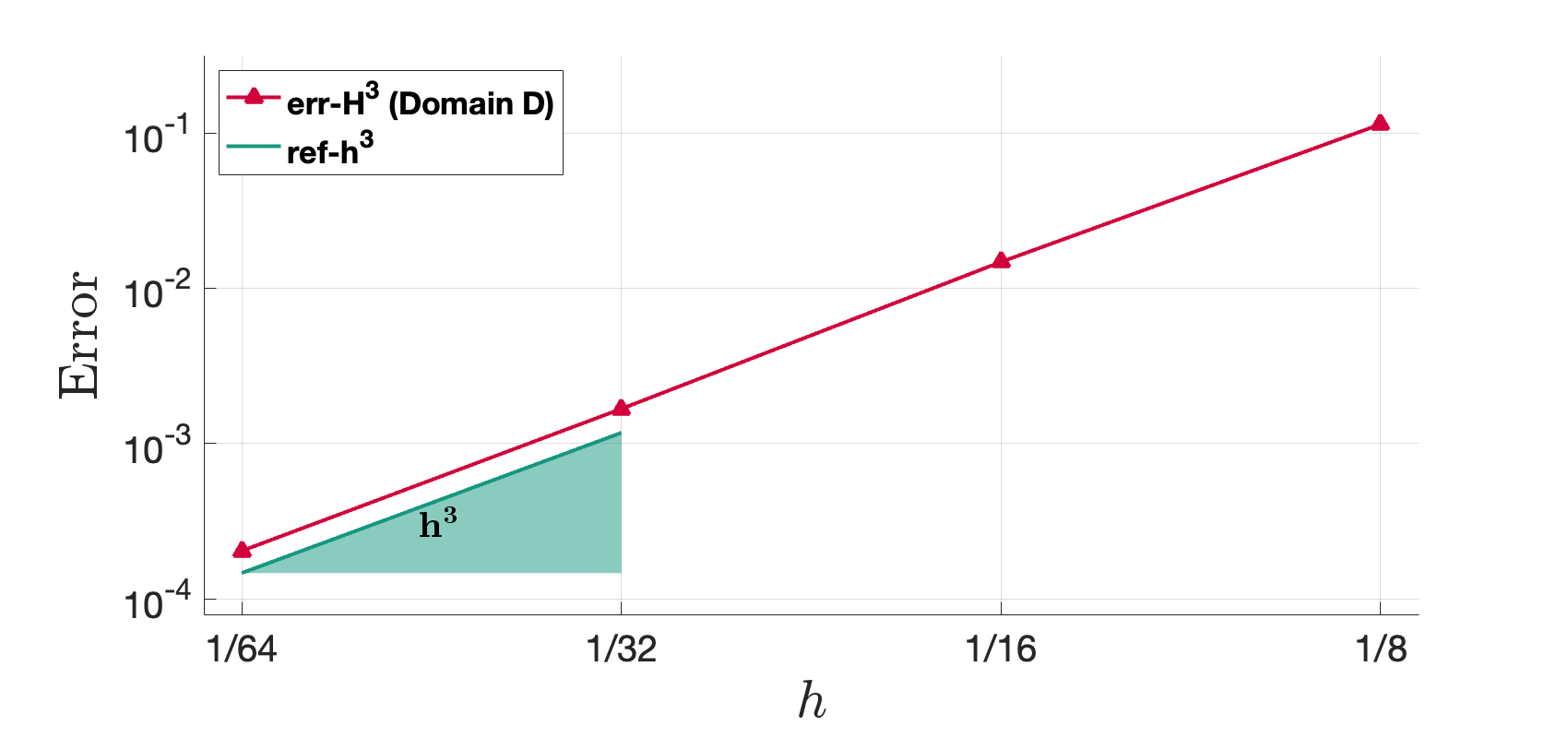}
\caption{Example~\ref{ex:Triharmonic}. 
The considered exact solution~\eqref{eq:exactSolution} (left) and the resulting relative errors~\eqref{eq:eqiuv2seminorms} with respect to the $H^3$-seminorm (right) for solving the triharmonic equation over the bilinear-like $G^2$ five-patch Domain D.
}
\label{fig:triharmonicBL}
\end{figure}
\end{ex}

\begin{ex} \label{ex:MixedSpace}
In Examples~\ref{ex:Biharmonic} and \ref{ex:Triharmonic} we have utilized underlying spline spaces $\mathcal{S}_{h}^{\ab{p},\ab{r}}([0,1]^2)$ with degree $\ab{p}=(p,p)=(2m-1,2m-1)$ and regularity $\ab{r}=(r,r)=(m-1,m-1)$. The chosen degree~$p=2m-1$ represents the lowest possible degree while the regularity~$r=m-1$ is the highest possible to ensure an $h$-refinable isogeometric spline space~$\W_h^s(\overline{\Omega})$, cf. Section~\ref{subsec:Cs_space}. By relaxing the condition to use underlying spline spaces with the same regularity~$r=m-1$ everywhere, it is possible to use underlying spline spaces with this regularity just in the vicinity of $\partial [0,1]^2$, where it is needed to connect the solution of the polyharmonic equation with $C^s$-smoothness across the single patches~$\Omega^{(i)}$, $i \in \mathcal{I}_{\Omega}$, and with the increased regularity~$r_1=p-1=2m-2$ within the interior of $[0,1]^2$. Such mixed regularity underlying spline spaces were introduced in \cite{KaKoVi24b} to construct mixed regularity isogeometric multi-patch spline spaces that possess the lower regularity $r=m-1$ just in a small neighborhood around the edges and vertices of the multi-patch domain, and have the higher regularity $r_1=p-1=2m-2$ in all other parts of the domain. We denote these mixed regularity underlying spline spaces by $\mathcal{S}_h^{\ab{p},(\ab{r}_1,\ab{r})}([0,1]^2)$ with degree $\ab{p}=(p,p)$, lower regularity $\ab{r}=(r,r)$ and higher regularity~$\ab{r}_1=(r_1,r_1)$, and refer for further details on these spaces to~\cite{KaKoVi24b}. Since the basis functions of the mixed regularity underlying spline spaces~$\mathcal{S}_h^{\ab{p},(\ab{r}_1,\ab{r})}([0,1]^2)$ coincide in the vicinity of $\partial [0,1]^2$ with those of the underlying spline spaces~$\mathcal{S}_{h}^{\ab{p},\ab{r}}([0,1]^2)$ with the same regularity everywhere, and differ only within the interior of $[0,1]^2$, the mixed regularity underlying spline spaces~$\mathcal{S}_h^{\ab{p},(\ab{r}_1,\ab{r})}([0,1]^2)$ can be seamlessly integrated into our IETI-based method. Thereby, the integration requires only an adaption in the assembly of the local stiffness matrices~$\ab{K}^{(i)}$.

We now solve the polyharmonic equations of order $m=2$ (the biharmonic equation) and $m=3$ (the triharmonic equation) over the two bilinear multi-patch domains (Domain A and Domain B), presented in Fig.~\ref{fig:domains} (top row). For this purpose, we use on the one hand the underlying spline spaces $\mathcal{S}_h^{\ab{p},\ab{r}}([0,1]^2) = \mathcal{S}_h^{\ab{2m-1},\ab{m-1}}([0,1]^2)$ with the same regularity (SR) everywhere, and on the other hand the mixed regularity (MR) underlying spline spaces $ \mathcal{S}_h^{\ab{p},(\ab{r}_1,\ab{r})}([0,1]^2) = \mathcal{S}_h^{\ab{2m-1},(\ab{2m-2},\ab{m-1})}([0,1]^2)$, and compare the obtained results. For the comparison, we study the convergence behavior under $h$-refinement with respect to relative errors~\eqref{eq:eqiuv2seminorms} in the $H^m$-seminorm by considering mesh sizes~$h= h_{0}/2^j$, $j=0,\ldots,3$, with $h_{0}=1/5$. The obtained convergence results with respect to the number of degrees of freedom (DOF) are visualized in Fig.~\ref{fig:mixed} and show for both underlying spline spaces optimal convergence behavior, which means a convergence of order~$\mathcal{O}(h^{2})$ in the $H^2$-seminorm when solving the biharmonic equation ($m=2$) and a convergence of order~$\mathcal{O}(h^{3})$ in the $H^3$-seminorm when solving the triharmonic equation ($m=3$). However, the magnitude of the relative errors with respect to the number of degrees of freedom (DOF) is for the mixed regularity (MR) underlying spline spaces much lower than for the underlying spline spaces with the same regularity (SR) everywhere. 

\begin{figure}[h!]
    \centering
    \includegraphics[scale=0.128]{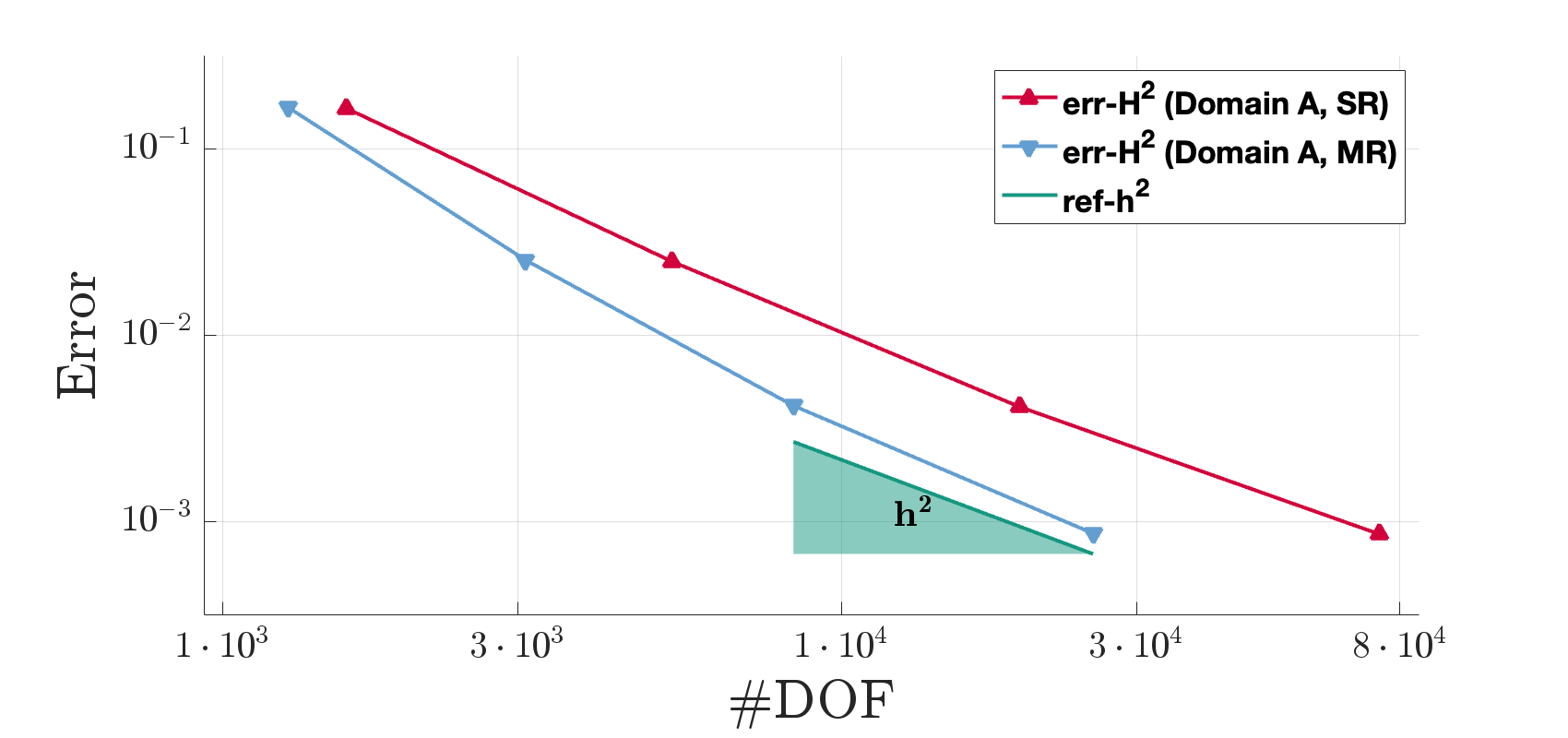}
    \includegraphics[scale=0.128]{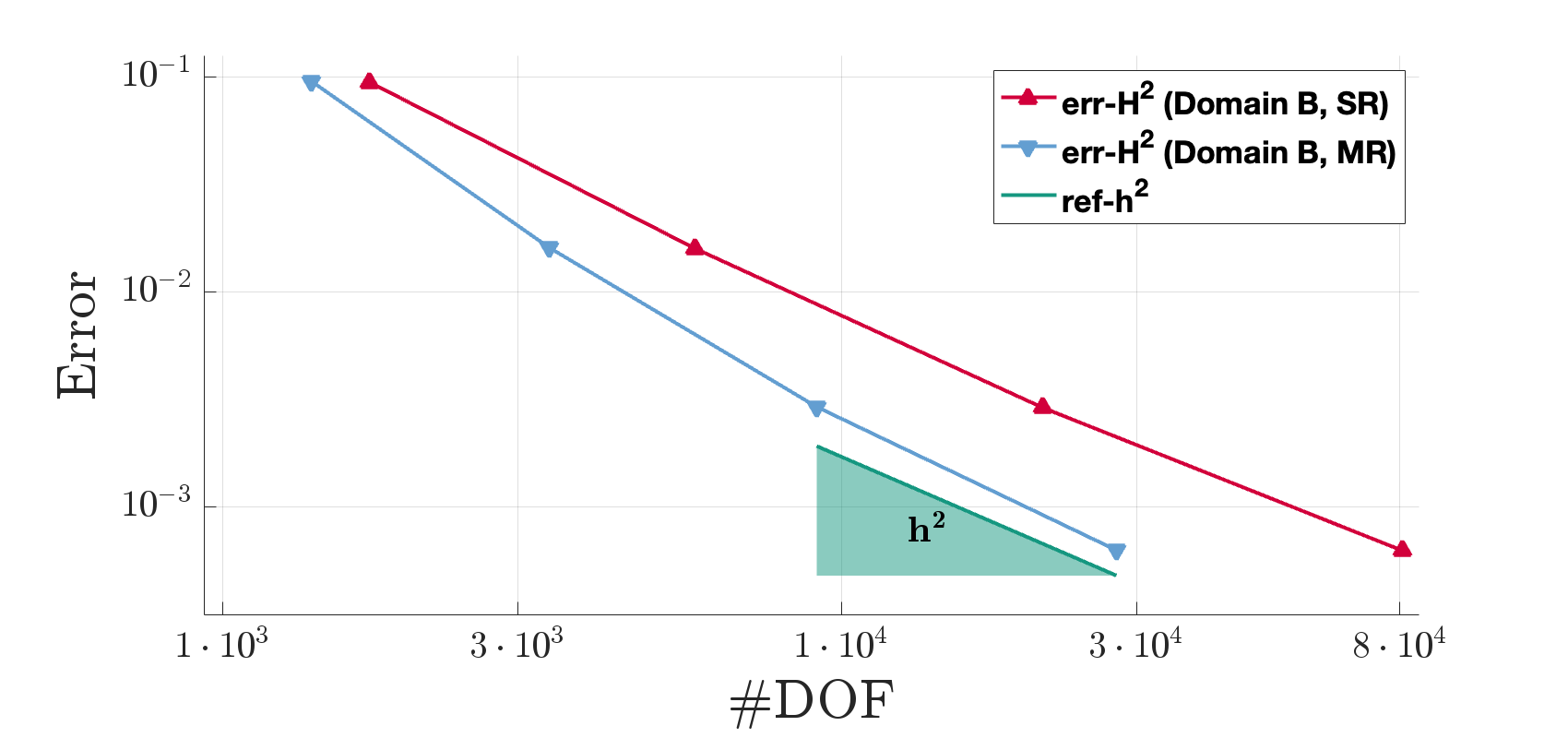}
    \includegraphics[scale=0.128]{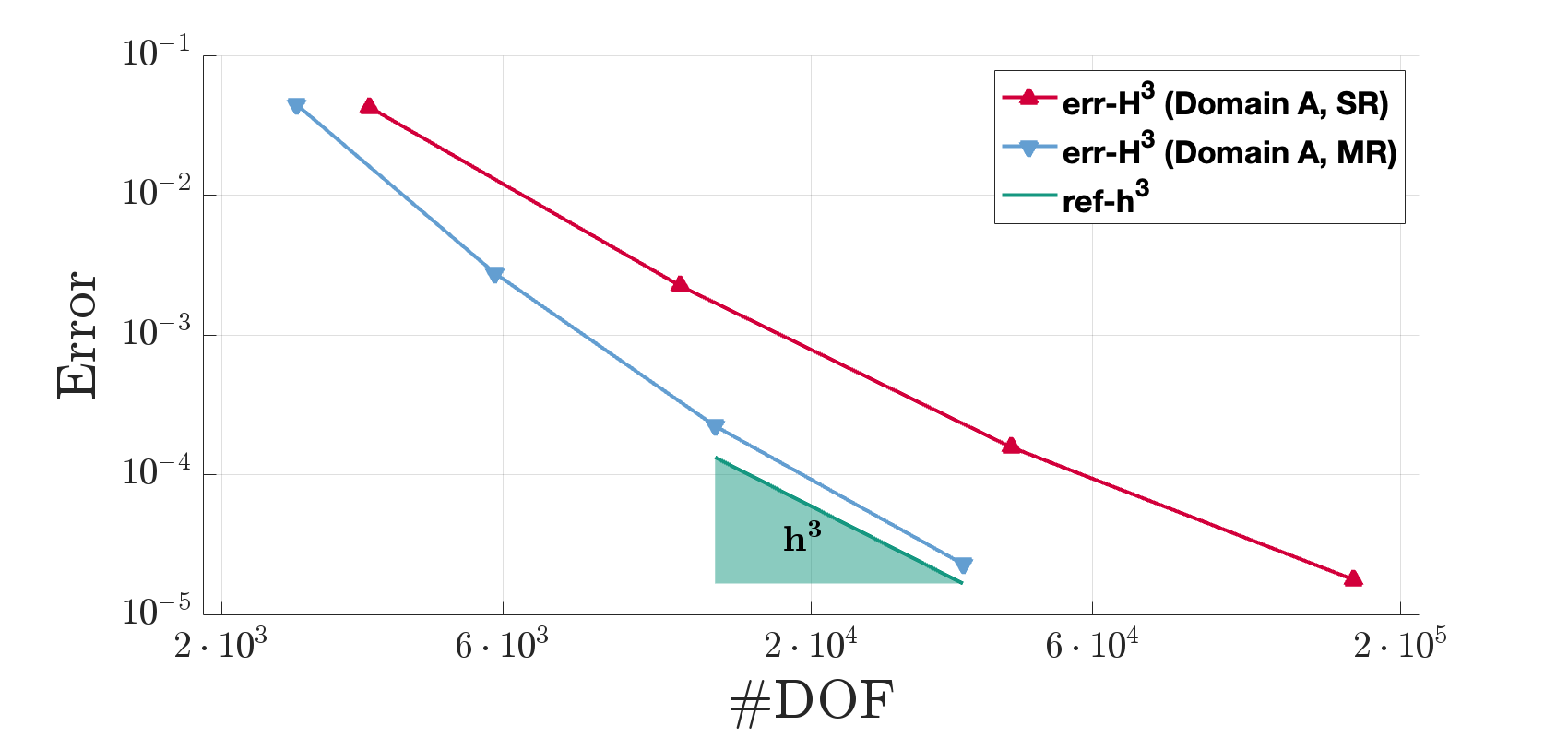}
    \includegraphics[scale=0.128]{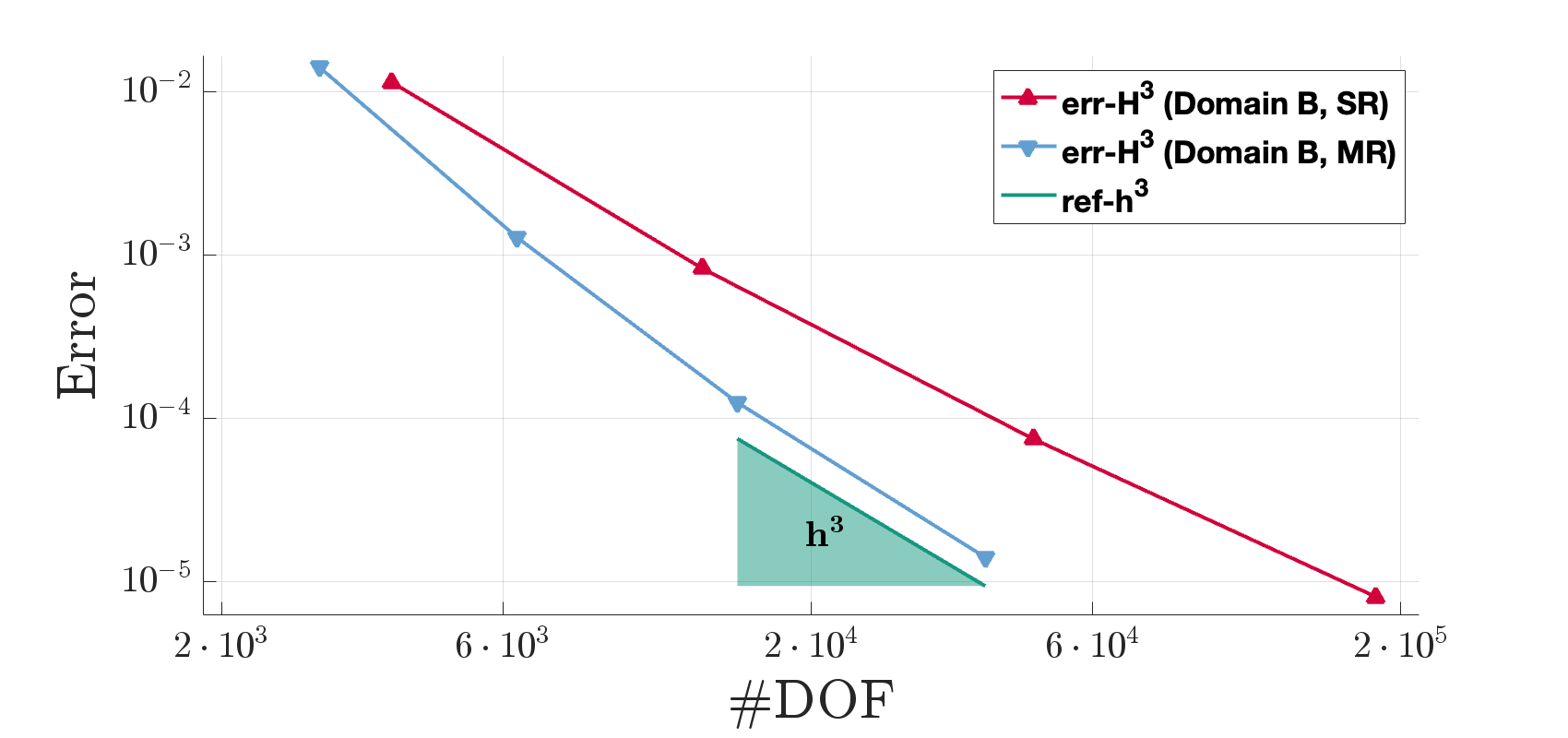}
\caption{Example~\ref{ex:MixedSpace}.
Relative errors~\eqref{eq:eqiuv2seminorms} with respect to the number of degrees of freedom (DOF) for the $H^2$-seminorm when solving the biharmonic equation ($m=2$), see top row, and for the $H^3$-seminorm when solving the triharmonic equation ($m=3$), see bottom row, over the two bilinear multi-patch domains (Domain A and B) shown in Fig.~\ref{fig:domains} (top row). 
The red color corresponds to the convergence results using the underlying spline spaces $\mathcal{S}_h^{\ab{2m-1},\ab{m-1}}([0,1]^2)$ with the same regularity (SR) everywhere, while the blue color presents the results using the mixed regularity (MR) underlying spline spaces $\mathcal{S}_h^{\ab{2m-1},(\ab{2m-2},\ab{m-1})}([0,1]^2)$.
}
\label{fig:mixed}
\end{figure}
\end{ex}

\section{Conclusion} \label{sec:Conclusion}

We presented a novel IETI-based method for solving the polyharmonic equation of order~$m \geq 1$, which is a particular linear elliptic PDE of order $2m$, over the class of so-called bilinear-like $G^s$ multi-patch geometries, $s \geq m-1$. Our technique is based on the solving of a dual-primal formulation of a saddle point problem, where the $C^s$-smoothness of the numerical solution across the patches is enforced by means of Lagrange multipliers. More precisely, our proposed algorithm efficiently computes first a subset of the Lagrange multipliers via a small linear system, followed by the remaining Lagrange multipliers as well as the degrees of freedom of the numerical solution via local, parallelizable linear systems for the single patches. The developed IETI-based method was tested on several numerical examples to solve the biharmonic and triharmonic equation, i.e. the polyharmonic equation of order $m=2$ and $m=3$, respectively, over bilinear-like $G^s$ multi-patch geometries. Thereby, the resulting optimal convergence behavior in the studied $H^m$-seminorm with respect to $h$-refinement demonstrated the great potential of our technique for solving higher order PDEs over planar multi-patch geometries with possibly extraordinary vertices.    

Possible topics for future research could be the extension of our method to multi-patch surfaces or multi-patch volumes on the one hand and to other high order PDEs on the other hand such as fourth order problems like the Kirchhoff-Love shell problem~\cite{kiendl-bletzinger-linhard-09} and the Cahn-Hilliard equation~\cite{gomez2008isogeometric} or sixth order problems such as the Phase-field crystal equation~\cite{Gomez2012}, the Kirchhoff plate model based on the Mindlin’s gradient elasticity theory~\cite{Niiranen2016} and the gradient-enhanced continuum damage models~\cite{GradientDamageModels}. Another interesting topic could be the development of an adaptive IETI method based on the use of hierarchical splines~\cite{Kr97} or truncated hierarchical splines~\cite{GiJuSp2012} to represent the numerical solution of the high order PDEs over the individual patches.

\paragraph*{\bf Acknowledgment}

V.~Vitrih has been partially supported by the 
Slovenian Research and Innovation Agency (research program P1-0404 and research project N1-0296). A.~Kosma\v c has been partially supported by the Slovenian Research and Innovation Agency (research program P1-0404, research project N1-0296 and Young Researchers Grant). This support is gratefully acknowledged.


\appendix

\section{Rank deficiencies of stiffness matrices} \label{sec:AppendixA}

It is known that the stiffness matrices $\ab{K}^{(i)}$ (and consequently the matrix $\ab{K}$) are not invertible. Their rank deficiency depends on the used spline space $\mathcal{S}^{\ab{p},\ab{r}}_h$ as numerically experimented for solving the biharmonic ($m=2$) and triharmonic equation ($m=3$) for $m \leq p \leq 9$ over several bilinear multi-patch domains. The obtained upper bounds for their rank deficiency are presented in Tables~\ref{tab:rankDef_bih} and \ref{tab:rankDef_trih} and indicate to be $2p+m$.

    \begin{table}[htb!]
    \centering
    \begin{tabular}{c||c|c|c|c|c|c|c|c|}
            $p\backslash r$  & $1$ & $2$ & $3$ & $4$ & $5$ & $6$ & $7$ & $8$\\
              \hline
       $2$  & 5& & & & & & &\\
       $3$  & 8 & 8 & & & & & &\\
       $4$  & 9 & 9 & 9 & & & & &\\
       $5$  & 12 & 12 & 12 & 12 & & & &\\
       $6$  & 13 & 13 & 13 & 13 & 13 & & &\\
       $7$  & 16 & 16 & 16 & 16 & 16 & 16 & &\\
       $8$  & 17 & 17 & 17 & 17 & 17 & 17 & 17 &\\
       $9$  & 20 & 20 & 20 & 20 & 20 & 20 & 20 & 20\\
    \end{tabular}
    \caption{Upper bounds for the rank deficiencies of the matrices $\ab{K}^{(i)}$ depending on ${p}$ and ${r}$ for the biharmonic equation.}
    \label{tab:rankDef_bih}
\end{table}
\begin{table}[htb!]
    \centering
         \begin{tabular}{c||c|c|c|c|c|c|c|}
            $p\backslash r$ & $2$ & $3$ & $4$ & $5$ & $6$ & $7$ & $8$\\
              \hline
       $3$  & 9 &  & & & & &\\
       $4$  & 10  & 10 &  & & & &\\
       $5$  & 13 & 13 & 13 &  & & &\\
       $6$  & 14 & 14 & 14 & 14 &  & &\\
       $7$  & 17 & 17 & 17 & 17 & 17 &  & \\
       $8$  & 18 & 18 & 18 & 18 & 18 & 18 &  \\
       $9$  & 21 & 21 & 21 & 21 & 21 & 21 & 21 \\
    \end{tabular}
    \caption{Upper bounds for the rank deficiencies of the matrices $\ab{K}^{(i)}$ depending on ${p}$ and ${r}$ for the triharmonic equation.}
    \label{tab:rankDef_trih}
\end{table}

\section{An efficient computation of the inverse of the matrix $\ab{S}$} \label{sec:AppendixB}

Recall the structure of the matrix $\ab{S}$ from  \eqref{eq:Schur_complement} and \eqref{eq:decomposed_problem_Schur_withS_Notation}, which 
can now be written as
$$
\ab{S} = \left( \begin{array}{cc}
\ab{S}_P &  {\ab{C}}_{\Xi}^T \\
{\ab{C}}_{\Xi} & \ab{0} 
\end{array} \right) = \left( \begin{array}{cc}
\ab{I} &  \ab{0} \\
\ab{0} & \widetilde{\ab{S}} 
\end{array} \right), \quad 
\widetilde{\ab{S}} = \left( \begin{array}{cc}
\ab{S}_F &  {\ab{C}}_{\Xi,F}^T \\
{\ab{C}}_{\Xi,F} & \ab{0} 
\end{array} \right).
$$
In order to compute the inverse of $\ab{S}$ efficiently, we have to compute the inverse of $\widetilde{\ab{S}}$ efficiently. Recall further that $
\ab{S}_{F} = \Diag (\{ \ab{S}_{F}^{(i)}\}_{i \in \mathcal{I}_{\Omega}})
$ and that $\ab{S}_{F}^{(i)}, \, i\in \mathcal{I}_\Omega^I$, is in general not invertible. Moreover, the structure of the matrix ${\ab{C}}_{\Xi,F}$ is presented in \eqref{eq:C_XiF} and \eqref{eq:C_XiF_blocks}.
The goal now is to slightly extend and restructure the matrix $\widetilde{\ab{S}}$ into matrix $\overline{\ab{S}}$, whose inverse will be possible to compute efficiently in a parallel way, and then show how to compute $\widetilde{\ab{S}}^{-1}$ back from $\overline{\ab{S}}^{-1}$. 

Let us first extend the matrix ${\ab{C}}_{\Xi,F}$ to the matrix
$$
{\ab{C}}_{\Xi,F}^J := \left( \begin{array}{cc}
{\widehat{\ab{C}}}_{\Xi,F} &  {\widehat{\ab{J}}} 
\end{array} \right), \;\;
    \widehat{\ab{C}}_{\Xi,F} =\left[
\widehat{{\ab{C}}}_{\Xi,F}^{(i)} \right]_{i \in \mathcal{I}_{\Xi}^{I}}, \;\;
\widehat{{\ab{J}}} = \left[
\widehat{{\ab{J}}}^{(i)}
\right]_{i \in \mathcal{I}_{\Xi}^{I}},
$$
where $\widehat{{\ab{C}}}_{\Xi,F}^{(i)} \in \R^{\nu_i \binom{2s+2}{2} \times |\ab{u}_F|}$ and its $\ell$-th block row consists now of only one nonzero matrix, given in \eqref{eq:C_XiF_blocks}, namely 
$\ab{C}^{(i,i_{\ell-1})}_{\Xi,F} $, $\ell=1,\ldots,\nu_i$. Moreover, the matrix $\widehat{\ab{J}}^{(i)}$ is of the form
$$
\widehat{\ab{J}}^{(i)} = \left( \begin{array}{ccccccc}
\ab{0} &  \cdots & \ab{0} & \ab{I} & \ab{0} & \cdots & \ab{0} \\
\vdots & \vdots & \vdots & \vdots & \vdots & \vdots & \vdots \\
\ab{0} &  \cdots & \ab{0} & \ab{I} & \ab{0}& \cdots & \ab{0}
\end{array}\right) \in \R^{\nu_i \binom{2s+2}{2} \times 
|\mathcal{I}_\Xi^I|\binom{2s+2}{2}},
$$
where the identity blocks correspond to the $i$-th block column. Therefrom we get the extension of the matrix $\widetilde{\ab{S}}$ to the matrix
\begin{equation*}
    \widehat{\ab{S}} = \left( \begin{array}{ccc}
\ab{S}_F & \ab{0} & \widehat{\ab{C}}_{\Xi,F}^T \\
\ab{0} & \ab{0} & \widehat{\ab{J}}^T \\
\widehat{\ab{C}}_{\Xi,F} & \widehat{\ab{J}} & \ab{0}
\end{array} \right).
\end{equation*}
Recall that $\widetilde{\ab{S}}$ acts on $  
\left(\ab{u}_F^T,  
\ab{\lambda}_\Xi ^T
\right)^T
$, which implies that the matrix $\widehat{\ab{S}}$
acts on  
$
\left(
\ab{u}_F^T, 
\ab{\mu}^T,  
\widehat{\ab{\lambda}}_\Xi^T
\right)^T,
$
where
$
\ab{\mu} \in \R^{|\mathcal{I}_\Xi^I| \binom{2s+2}{2}}
$ are some additional Lagrange multipliers, and where the Lagrange multipliers 
${\ab{\lambda}}_\Xi = \left( {\ab{\lambda}}_{\Xi^{(i)}}\right)_{i=1}^{|\mathcal{I}_\Xi^I|}, \; {\ab{\lambda}}_{\Xi^{(i)}} \in \R^{(\nu_i-1)\binom{2s+2}{2}}$, and 
$\widehat{\ab{\lambda}}_\Xi = \big( {\widehat{\ab{\lambda}}}_{\Xi^{(i)}}\big)_{i=1}^{|\mathcal{I}_\Xi^I|}, \; {\widehat{\ab{\lambda}}}_{\Xi^{(i)}} \in \R^{\nu_i\binom{2s+2}{2}}$, 
are related as
\begin{equation}  \label{eq:matricesDi}
{\ab{\lambda}}_{\Xi^{(i)}} = \ab{L}^{(i)}\, \widehat{\ab{\lambda}}_{\Xi^{(i)}}, \quad
\ab{L}^{(i)} = 
\left( \begin{array}{ccccc}
\ab{I} &  \ab{0} & & &   \\
\ab{I} &  \ab{I} & \ab{0} & &   \\
\vdots &  \vdots & \ddots & \ddots &   \\
\ab{I} &  \ab{I} & \cdots &  \ab{I} & \ab{0} \\
\end{array} \right) \in \R^{(\nu_i-1) \binom{2s+2}{2} \, \times \, \nu_i \binom{2s+2}{2}}.
\end{equation}
Let us now permute rows and columns of the matrix $\widehat{\ab{S}}$ with a permutation matrix $\ab{P}$ to get the matrix
$$  
\overline{\ab{S}} = \ab{P} \widehat{\ab{S}}  \ab{P}^T, 
$$
such that 
\begin{equation*}
  \overline{\ab{S}} = \left( \begin{array}{ccccc|c}
\ab{S}_F^{(1)} & {\ab{C}^{(1)}}^T & & & & \ab{0}\\
{\ab{C}^{(1)}} & \ab{0} & & & & {\ab{J}^{(1)}}^T \\
 &  & \ddots & & & \vdots \\
 & & & \ab{S}_F^{(|\mathcal{I}_\Omega|)} & {\ab{C}^{(|\mathcal{I}_\Omega|)}}^T & \ab{0}\\
 & &  & {\ab{C}^{(|\mathcal{I}_\Omega|)}} & \ab{0} & {\ab{J}^{(|\mathcal{I}_\Omega|)}}^T \\
 \hline\\[-0.4cm]
\ab{0} & \ab{J}^{(1)} & \cdots & \ab{0} & \ab{J}^{(|\mathcal{I}_\Omega|)} & \ab{0}
\end{array} \right) =: \left(\begin{array}{cc} 
\ab{T} & \ab{J}^T \\
\ab{J} & \ab{0}
\end{array} \right), 
\end{equation*}
where 
matrix $\ab{C}^{(j)}$, $j\in \mathcal{I}_\Omega$, consists of all matrices $\ab{C}^{(i,j)}_{\Xi,F}$, $i\in \mathcal{I}_\Xi^I$, in $\widehat{\ab{C}}_{\Xi,F}$, and matrix $\ab{J}^{(j)}$ is a block matrix with $|\mathcal{I}_\Xi^I|$ block rows, while the number of block columns equals the number of inner vertices of $\Omega^{(j)}$. Most blocks are zero blocks, except in each block column we have one identity block (in the block rows that correspond to inner vertices of $\Omega^{(j)}$).

Note that $\ab{T}$ is blockwise diagonal, i.e.
$$
\ab{T} = \Diag (\{ \ab{T}^{(j)}\}_{j \in \mathcal{I}_{\Omega}}), \quad 
\ab{T}^{(j)} = \left(\begin{array}{cc}
\ab{S}_F^{(j)} & {\ab{C}^{(j)}}^T  \\
\ab{C}^{(j)} &  \ab{0} 
\end{array} \right),
$$
and that all its diagonal blocks $\ab{T}^{(j)}$ are invertible. Therefore we can compute its inverse $\ab{T}^{-1}$ efficiently in a parallel way by solving local problems on each patch $\Omega^{(i)}$ by
(sparse) LU factorization.  
We can then compute
\begin{equation}  \label{eq:overlineSinv}
\overline{\ab{S}}^{\,-1} = 
\left( \begin{array}{cc}
\ab{T}^{-1} - \left( \ab{J}\ab{T}^{-1}  \right)^T \left( \ab{J}\ab{T}^{-1} \ab{J}^T\right)^{-1} \left( \ab{J}\ab{T}^{-1}  \right) \;\; & \;\;   \left( \ab{J}\ab{T}^{-1}  \right)^T \left( \ab{J}\ab{T}^{-1} \ab{J}^T\right)^{-1} \\
\left( \ab{J}\ab{T}^{-1} \ab{J}^T\right)^{-1} \left( \ab{J}\ab{T}^{-1}  \right) \; & \; -\left( \ab{J}\ab{T}^{-1} \ab{J}^T\right)^{-1} 
\end{array} \right),
\end{equation}
and finally 
\begin{equation}  \label{eq:tildeSinv}
 \widetilde{\ab{S}}^{-1} = \left(\ab{R} \ab{P}^T\right) \overline{\ab{S}}^{\,-1}
\left(\ab{R} \ab{P}^T\right)^T,   
\end{equation}
with
$$
\ab{R}= \left( \begin{array}{ccc}
\ab{I} & \ab{0} &   \\
 &  & \ab{L}
\end{array} \right), \quad
\ab{L}= \Diag (\{ \ab{L}^{(i)}\}_{i \in \mathcal{I}_{\Xi}^I}),
%
\quad \ab{I} \in R^{|\ab{u}_F| \times |\ab{u}_F|}, \quad \ab{0} \in \R^{|\ab{u}_F|\times |\ab{\mu}|},
$$
where the matrices $\ab{L}^{(i)}$ are given in \eqref{eq:matricesDi}. 
Note that the computations in \eqref{eq:overlineSinv} and \eqref{eq:tildeSinv} can be done very efficiently once $\ab{T}^{-1}$ is known. 

\section{Symmetric positive definiteness of the dual system \eqref{eq:solution_LambdaGamma} for $\ab{\lambda}_\Gamma$} \label{sec:AppendixC}
 
The aim is to show that the matrix $\widetilde{\ab{C}}_{\Gamma} \ab{S}^{-1} \widetilde{\ab{C}}_{\Gamma}^T$ is symmetric positive definite.  
Since the matrix is invertible (saddle point problem has a unique solution) it is enough to show symmetric positive semidefiniteness. 
Recall \eqref{eq:decomposed_problem_Schur_withS_Notation}, i.e.,
$$
\ab{S} = \left( \begin{array}{cc}
\ab{S}_P &  {\ab{C}}_{\Xi}^T \\
{\ab{C}}_{\Xi} & \ab{0} 
\end{array} \right) \quad {\rm and} \quad 
\widetilde{\ab{C}}_{\Gamma} = \left( \begin{array}{cc}
\overline{\ab{C}}_{\Gamma} \; \, \ab{0} \end{array}\right),
$$
where $\ab{S}_P$ is symmetric positive semidefinite matrix, ${\ab{C}}_{\Xi}$ has a full row rank, and
$\ker(\ab{S}_P) \cap \ker({\ab{C}}_{\Xi}) = \{\ab{0}\}$.
Note that $\ab{S}$ is in general not symmetric positive semidefinite, so also $\ab{S}^{-1}$ is not. 
Writing
$$
\ab{S}^{-1} =: \ab{\Sigma} = \left( \begin{array}{cc}
\ab{\Sigma}_{11} &  \ab{\Sigma}_{12} \\
\ab{\Sigma}_{21} & \ab{\Sigma}_{22}
\end{array} \right),
$$
the matrix $\widetilde{\ab{C}}_{\Gamma} \ab{S}^{-1} \widetilde{\ab{C}}_{\Gamma}^T$ simplifies to $\overline{\ab{C}}_{\Gamma} \ab{\Sigma}_{11} \overline{\ab{C}}_{\Gamma}^T$ and it suffices to show symmetric positive semidefiniteness for the matrix $\ab{\Sigma}_{11}$. 
Since $\ab{S} \ab{\Sigma} = \ab{I}$, we get the equations
\begin{equation} \label{eq:twoEquations_Sigma}
\ab{S}_P \ab{\Sigma}_{11} + {\ab{C}}_{\Xi}^T \ab{\Sigma}_{21} = \ab{I} \quad {\rm and} \quad 
{\ab{C}}_{\Xi} \ab{\Sigma}_{11} = \ab{0}. 
\end{equation}
The first one implies
\begin{equation} \label{eq:Sigma11}
   \ab{\Sigma}_{11} = \ab{S}_P^{+} \left(\ab{I} -  {\ab{C}}_{\Xi}^T \ab{\Sigma}_{21}\right) + \ab{N} \ab{\Psi},  
\end{equation}
where $\ab{S}_P^{+}$ is the Moore-Penrose pseudoinverse of matrix $\ab{S}_P$, the columns of the matrix $\ab{N}$ form a basis of the null space of $\ab{S}_P$ and each column in $\ab{\Psi}$ represents certain linear combinations of the columns in $\ab{N}$. 
The second equation in \eqref{eq:twoEquations_Sigma} now gives
$
{\ab{C}}_{\Xi} \ab{S}_P^{+} \left(\ab{I} -  {\ab{C}}_{\Xi}^T \ab{\Sigma}_{21}\right) + {\ab{C}}_{\Xi} \ab{N} \ab{\Psi} = \ab{0},
$
which implies 
\begin{equation} \label{eq:Sigma21}
\ab{\Sigma}_{21} = \left( {\ab{C}}_{\Xi} \ab{S}_P^{+} {\ab{C}}_{\Xi}^T \right)^{-1} \left(
{\ab{C}}_{\Xi} \ab{S}_P^{+} + {\ab{C}}_{\Xi} \ab{N} \ab{\Psi} \right). 
\end{equation}
Note that under the assumptions on $\ab{S}_P$ and ${\ab{C}}_{\Xi}$, the matrix $\ab{V} := \left({\ab{C}}_{\Xi} \ab{S}_P^{+} {\ab{C}}_{\Xi}^T\right)^{-1}$ is invertible and symmetric positive definite. 
Inserting \eqref{eq:Sigma21} into \eqref{eq:Sigma11} gives further
\begin{equation} \label{eq:Sigma11_2}
\ab{\Sigma}_{11} = \ab{S}_P^{+} \left( \ab{I} - {\ab{C}}_{\Xi}^T \ab{V} {\ab{C}}_{\Xi} \ab{S}_P^+\right) +  \left( \ab{I} - \ab{S}_P^+ {\ab{C}}_{\Xi}^T \ab{V} {\ab{C}}_{\Xi} \right) \ab{N} \ab{\Psi}.  
\end{equation}
We can now add the orthogonality condition (cf.~\cite{Farhat1991}) 
$$
\ab{N}^T \left( \ab{I} - {\ab{C}}_{\Xi}^T \ab{\Sigma}_{21}\right) = \ab{0}.
$$
Inserting \eqref{eq:Sigma21} implies further
\begin{equation} \label{eq:Psi}
 \ab{\Psi} = \left( \left({\ab{C}}_{\Xi} \ab{N}\right)^T \ab{V} \left({\ab{C}}_{\Xi} \ab{N}\right) \right)^{-1} 
  \left( 
  \ab{N}^T \left( 
  \ab{I} - {\ab{C}}_{\Xi}^T \ab{V} {\ab{C}}_{\Xi} \ab{S}_P^+ \right)
  \right).
\end{equation}
Applying now \eqref{eq:Psi} in \eqref{eq:Sigma11_2} finally gives
\begin{align}  \label{eq:finalSigma11}
\ab{\Sigma}_{11} = & \left( \ab{S}_P^+ - \left({\ab{C}}_{\Xi} \ab{S}_P^+\right)^T  \ab{V} \left({\ab{C}}_{\Xi} \ab{S}_P^+\right) \right) + \nonumber \\[-0.7cm]
\,\\
& \left(  \left( \ab{I} - \left({\ab{C}}_{\Xi} \ab{S}_P^+\right)^T \ab{V} {\ab{C}}_{\Xi} \right) \ab{N} \left( \left({\ab{C}}_{\Xi} \ab{N}\right)^T \ab{V} \left({\ab{C}}_{\Xi} \ab{N}\right) \right)^{-1} \ab{N}^T \left( \ab{I} - \left({\ab{C}}_{\Xi} \ab{S}_P^+\right)^T \ab{V} {\ab{C}}_{\Xi} \right)^T \right). \nonumber
\end{align}
One can observe that the first part of \eqref{eq:finalSigma11} is symmetric positive semidefinite. 
Moreoever, the second part is also clearly symmetric positive semidefinite, since $\ab{V}$ is symmetric positive definite. The sum of two symmetric positive semidefinite matrices is also a symmetric positive semidefinite matrix and this completes the proof.


\begin{thebibliography}{10}

\bibitem{ArReKlSi23}
J.~Arf, M.~Reichle, S.~Klinkel, and B.~Simeon.
\newblock Scaled boundary isogeometric analysis with ${C}^1$ coupling for
  {K}irchhoff plate theory.
\newblock {\em Comput. Methods Appl. Mech. Engrg.}, 415:116198, 2023.

\bibitem{ANU:9260759}
L.~{Beir{\~a}o da Veiga}, A.~Buffa, G.~Sangalli, and R.~V\'azquez.
\newblock Mathematical analysis of variational isogeometric methods.
\newblock {\em Acta Numerica}, 23:157--287, 5 2014.

\bibitem{BeLoSaTa23}
A.~Benvenuti, G.~Loli, G.~Sangalli, and T.~Takacs.
\newblock Isogeometric multi-patch ${C}^1$-mortar coupling for the biharmonic
  equation.
\newblock {\em arXiv}, 2303.07255, 2023.

\bibitem{BeEvMcTa21}
J.~Benzaken, J.~A. Evans, S.~F. McCormick, and R.~Tamstorf.
\newblock {N}itsche’s method for linear {K}irchhoff–{L}ove shells:
  {F}ormulation, error analysis, and verification.
\newblock {\em Comput. Methods Appl. Mech. Engrg.}, 374:113544, 2021.

\bibitem{Bouclier2017}
R.~Bouclier, J.~C. Passieux, and M.~Sala{\"{u}}n.
\newblock {Development of a new, more regular, mortar method for the coupling
  of NURBS subdomains within a NURBS patch: Application to a non-intrusive
  local enrichment of NURBS patches}.
\newblock {\em Comput. Methods Appl. Mech. Engrg.}, 316:123--150, 2017.

\bibitem{IETI_LowRank2024}
A.~B\"unger, T.-C. Riemer, and M.~Stoll.
\newblock {IETI}-based low-rank method for {PDE}-constrained optimization.
\newblock {\em arXiv}, 2405.06458, 2024.

\bibitem{ChAnRa18}
C.~L. Chan, C.~Anitescu, and T.~Rabczuk.
\newblock Isogeometric analysis with strong multipatch {C}$^1$-coupling.
\newblock {\em Comput. Aided Geom. Design}, 62:294--310, 2018.

\bibitem{ChAnRa19}
C.~L. Chan, C.~Anitescu, and T.~Rabczuk.
\newblock Strong multipatch {C}$^1$-coupling for isogeometric analysis on {2D}
  and {3D} domains.
\newblock {\em Comput. Methods Appl. Mech. Engrg.}, 357:112599, 2019.

\bibitem{CIOrSch00}
F.~Cirak, M.~Ortiz, and P.~Schr{\"o}der.
\newblock Subdivision surfaces: {A} new paradigm for thin-shell finite element
  analysis.
\newblock {\em Int. J. Numer. Methods Eng.}, 47(12):2039--2072, 2000.

\bibitem{CoSaTa16}
A.~Collin, G.~Sangalli, and T.~Takacs.
\newblock Analysis-suitable {G}$^1$ multi-patch parametrizations for {C}$^1$
  isogeometric spaces.
\newblock {\em Comput. Aided Geom. Des.}, 47:93 -- 113, 2016.

\bibitem{CoKiBu21}
L.~Coradello, J.~Kiendl, and A.~Buffa.
\newblock Coupling of non-conforming trimmed isogeometric {Kirchhoff}-{Love}
  shells via a projected super-penalty approach.
\newblock {\em Comput. Methods Appl. Mech. Eng.}, 387:114187, 2021.

\bibitem{CoLoBu21}
L.~Coradello, G.~Loli, and A.~Buffa.
\newblock A projected super-penalty method for the $c^1$-coupling of
  multi-patch isogeometric {Kirchhoff} plates.
\newblock {\em Comput. Mech.}, 67(4):1133--1153, 2021.

\bibitem{CottrellBook}
J.~A. Cottrell, T.~J.~R. Hughes, and Y.~Bazilevs.
\newblock {\em Isogeometric Analysis: Toward Integration of {CAD} and {FEA}}.
\newblock John Wiley \& Sons, Chichester, England, 2009.

\bibitem{DiSchWoHe19}
M.~Dittmann, S.~Schu{\ss}, B.~Wohlmuth, and C.~Hesch.
\newblock Weak ${C}^n$ coupling for multipatch isogeometric analysis in solid
  mechanics.
\newblock {\em Int. J. Numer. Methods Eng.}, 118(11):678--699, 2019.

\bibitem{DuRoSa17}
T.~X. Duong, F.~Roohbakhshan, and R.~A. Sauer.
\newblock A new rotation-free isogeometric thin shell formulation and a
  corresponding continuity constraint for patch boundaries.
\newblock {\em Comput. Methods Appl. Mech. Eng.}, 316:43--83, 2017.

\bibitem{FaJuKaTa22}
A.~Farahat, B.~J\"uttler, M.~Kapl, and T.~Takacs.
\newblock Isogeometric analysis with ${C}^1$-smooth functions over multi-patch
  surfaces.
\newblock {\em Comput. Methods Appl. Mech. Engrg.}, 403:115706, 2023.

\bibitem{FaKaKoVi24}
A.~Farahat, M.~Kapl, A.~Kosma\v{c}, and V.~Vitrih.
\newblock A locally based construction of analysis-suitable ${G}^1$ multi-patch
  spline surfaces.
\newblock {\em Comput. Math. Appl.}, 168:46--57, 2024.

\bibitem{FarhatFETI-DP}
C.~Farhat, M.~Lesoinne, P.~LeTallec, K.~Pierson, and D.~Rixen.
\newblock {FETI-DP}: a dual–primal unified {FETI} method—part {I}: {A}
  faster alternative to the two-level {FETI} method.
\newblock {\em Int. J. Numer. Methods Eng.}, 50(7):1523--1544, 2001.

\bibitem{Farhat1991}
C.~Farhat and F.-X. Roux.
\newblock A method of finite element tearing and interconnecting and its
  parallel solution algorithm.
\newblock {\em Int. J. Numer. Methods Eng.}, 32(6):1205--1227, 1991.

\bibitem{Fa97}
G.~Farin.
\newblock {\em Curves and Surfaces for Computer-Aided Geometric Design}.
\newblock Academic Press, 1997.

\bibitem{GiJuSp2012}
C.~Giannelli, B.~Jüttler, and H.~Speleers.
\newblock {THB}-splines: The truncated basis for hierarchical splines.
\newblock {\em Computer Aided Geometric Design}, 29(7):485--498, 2012.
\newblock Geometric Modeling and Processing 2012.

\bibitem{gomez2008isogeometric}
H.~Gomez, V.~M Calo, Y.~Bazilevs, and T.~J.R. Hughes.
\newblock Isogeometric analysis of the {Cahn--Hilliard} phase-field model.
\newblock {\em Comput. Methods Appl. Mech. Engrg.}, 197(49):4333--4352, 2008.

\bibitem{Gomez2012}
H.~Gomez and X.~Nogueira.
\newblock An unconditionally energy-stable method for the phase field crystal
  equation.
\newblock {\em Comput. Methods Appl. Mech. Engrg.}, 249 -- 252:52 -- 61, 2012.

\bibitem{Guo2015881}
Y.~Guo and M.~Ruess.
\newblock Nitsche's method for a coupling of isogeometric thin shells and
  blended shell structures.
\newblock {\em Comp. Methods Appl. Mech. Engrg.}, 284:881--905, 2015.

\bibitem{HoLa17}
C.~Hofer and U.~Langer.
\newblock Dual-primal isogeometric tearing and interconnecting solvers for
  multipatch {dG}-{Iga} equations.
\newblock {\em Comput. Methods Appl. Mech. Eng.}, 316:2--21, 2017.

\bibitem{HoLa93}
J.~Hoschek and D.~Lasser.
\newblock {\em Fundamentals of computer aided geometric design}.
\newblock A K Peters Ltd., Wellesley, MA, 1993.

\bibitem{HuCoBa04}
T.~J.~R. Hughes, J.~A. Cottrell, and Y.~Bazilevs.
\newblock Isogeometric analysis: {CAD}, finite elements, {NURBS}, exact
  geometry and mesh refinement.
\newblock {\em Comput. Methods Appl. Mech. Engrg.}, 194(39-41):4135--4195,
  2005.

\bibitem{KaKoVi24}
M.~Kapl, A.~Kosma\v{c}, and V.~Vitrih.
\newblock Isogeometric collocation for solving the biharmonic equation over
  planar multi-patch domains.
\newblock {\em Computer Methods in Applied Mechanics and Engineering},
  424:116882, 2024.

\bibitem{KaKoVi24c}
M.~Kapl, A.~Kosma\v{c}, and V.~Vitrih.
\newblock Isogeometric collocation with smooth mixed degree splines over planar
  multi-patch domains.
\newblock {\em arXiv}, 2411.13338, 2024.

\bibitem{KaKoVi24b}
M.~Kapl, A.~Kosma\v{c}, and V.~Vitrih.
\newblock A ${C}^s$-smooth mixed degree and regularity isogeometric spline
  space over planar multi-patch domains.
\newblock {\em Journal of Computational and Applied Mathematics}, 473:116836,
  2026.

\bibitem{KaSaTa17b}
M.~Kapl, G.~Sangalli, and T.~Takacs.
\newblock Construction of analysis-suitable {G}$^1$ planar multi-patch
  parameterizations.
\newblock {\em Comput.-Aided Des.}, 97:41--55, 2018.

\bibitem{KaSaTa19b}
M.~Kapl, G.~Sangalli, and T.~Takacs.
\newblock Isogeometric analysis with {C}$^{1}$ functions on unstructured
  quadrilateral meshes.
\newblock {\em The SMAI journal of computational mathematics}, 5:67--86, 2019.

\bibitem{KaVi17b}
M.~Kapl and V.~Vitrih.
\newblock Space of {C}$^2$-smooth geometrically continuous isogeometric
  functions on planar multi-patch geometries: {D}imension and numerical
  experiments.
\newblock {\em Comput. Math. Appl.}, 73(10):2319--2338, 2017.

\bibitem{KaVi17c}
M.~Kapl and V.~Vitrih.
\newblock Dimension and basis construction for ${C}^{2}$-smooth isogeometric
  spline spaces over bilinear-like ${G}^{2}$ two-patch parameterizations.
\newblock {\em J. Comput. Appl. Math.}, 335:289--311, 2018.

\bibitem{KaVi19a}
M.~Kapl and V.~Vitrih.
\newblock Solving the triharmonic equation over multi-patch planar domains
  using isogeometric analysis.
\newblock {\em J. Comput. Appl. Math.}, 358:385--404, 2019.

\bibitem{KaVi20}
M.~Kapl and V.~Vitrih.
\newblock Isogeometric collocation on planar multi-patch domains.
\newblock {\em Comput. Methods Appl. Mech. Engrg.}, 360:112684, 2020.

\bibitem{KaVi20b}
M.~Kapl and V.~Vitrih.
\newblock ${C}^s$-smooth isogeometric spline spaces over planar multi-patch
  parameterizations.
\newblock {\em Advances in Computational Mathematics}, 47:47, 2021.

\bibitem{KaPe17}
K.~Kar{\v c}iauskas and J.~Peters.
\newblock Refinable ${G}^1$ functions on ${G}^1$ free-form surfaces.
\newblock {\em Comput. Aided Geom. Des.}, 54:61--73, 2017.

\bibitem{KaPe18}
K.~Kar{\v c}iauskas and J.~Peters.
\newblock Refinable bi-quartics for design and analysis.
\newblock {\em Comput.-Aided Des.}, pages 204--214, 2018.

\bibitem{kiendl-bletzinger-linhard-09}
J.~Kiendl, K.-U. Bletzinger, J.~Linhard, and R.~W{\"u}chner.
\newblock Isogeometric shell analysis with {K}irchhoff-{L}ove elements.
\newblock {\em Comput. Methods Appl. Mech. Engrg.}, 198(49):3902--3914, 2009.

\bibitem{Klawonn2000}
A.~Klawonn and O.~B. Widlund.
\newblock A domain decomposition method with lagrange multipliers and inexact
  solvers for linear elasticity.
\newblock {\em SIAM Journal on Scientific Computing}, 22(4):1199--1219, 2000.

\bibitem{Klawonn2006}
A.~Klawonn and O.~B. Widlund.
\newblock Dual-primal feti methods for linear elasticity.
\newblock {\em Communications on Pure and Applied Mathematics},
  59(11):1523--1572, 2006.

\bibitem{KlPeSaJu12}
S.~K. Kleiss, C.~Pechstein, B.~J\"uttler, and S.~Tomar.
\newblock {IETI} -- isogeometric tearing and interconnecting.
\newblock {\em Comput. Methods Appl. Mech. Engrg.}, 247-248:201--215, 2012.

\bibitem{Kr97}
R.~Kraft.
\newblock Adaptive and linearly independent multilevel {B}-splines.
\newblock In {\em Surface fitting and multiresolution methods. Vol. 2 of the
  proceedings of the 3rd international conference on Curves and surfaces, held
  in Chamonix-Mont-Blanc, France, June 27--July 3, 1996}, pages 209--218.
  Nashville, TN: Vanderbilt University Press, 1997.

\bibitem{LeLiMaKiReGa20}
L.~Leonetti, F.~S. Liguori, D.~Magisano, J.~Kiendl, A.~Reali, and G.~Garcea.
\newblock A robust penalty coupling of non-matching isogeometric
  {Kirchhoff}-{Love} shell patches in large deformations.
\newblock {\em Comput. Methods Appl. Mech. Eng.}, 371:113289, 2020.

\bibitem{MaMaMo2024}
M.~Marsala, A.~Mantzaflaris, and B.~Mourrain.
\newblock ${G}^1$ spline functions for point cloud fitting.
\newblock {\em Applied Mathematics and Computation}, 460:128279, 2024.

\bibitem{MiZoScBoTh21}
D.~Miao, Z.~Zou, M.~A. Scott, M.~J. Borden, and D.~C. Thomas.
\newblock Isogeometric {B{\'e}zier} dual mortaring: the {Kirchhoff}-{Love}
  shell problem.
\newblock {\em Comput. Methods Appl. Mech. Eng.}, 382:113873, 2021.

\bibitem{MoSaSchTaTa23}
M.~Montardini, G.~Sangalli, R.~Schneckenleitner, S.~Takacs, and M.~Tani.
\newblock A {IETI}-{DP} method for discontinuous {Galerkin} discretizations in
  isogeometric analysis with inexact local solvers.
\newblock {\em Math. Models Methods Appl. Sci.}, 33(10):2085--2111, 2023.

\bibitem{mourrain2015geometrically}
B.~Mourrain, R.~Vidunas, and N.~Villamizar.
\newblock Dimension and bases for geometrically continuous splines on surfaces
  of arbitrary topology.
\newblock {\em Comput. Aided Geom. Des.}, 45:108 -- 133, 2016.

\bibitem{NgPe16}
T.~Nguyen and J.~Peters.
\newblock Refinable ${C}^{1}$ spline elements for irregular quad layout.
\newblock {\em Comput. Aided Geom. Des.}, 43:123--130, 2016.

\bibitem{Nguyen2014}
V.~P. Nguyen, P.~Kerfriden, M.~Brino, S.~P.A. Bordas, and E.~Bonisoli.
\newblock {Nitsche's method for two and three dimensional NURBS patch
  coupling}.
\newblock {\em Computational Mechanics}, 53(6):1163--1182, 2014.

\bibitem{Niiranen2016}
J.~Niiranen, J.~Kiendl, A.~H. Niemi, and A.~Reali.
\newblock Isogeometric analysis for sixth-order boundary value problems of
  gradient-elastic {K}irchhoff plates.
\newblock {\em Comput. Methods Appl. Mech. Engrg.}, 316:328--348, 2017.

\bibitem{ReArSiKl23}
M.~F.~M. Reichle, J.~Arf, B.~Simeon, and S.~Klinkel.
\newblock {S}mooth multi-patch scaled boundary isogeometric analysis for
  {K}irchhoff–{L}ove shells.
\newblock {\em Meccanica}, 58(8):1693--1716, 2023.

\bibitem{RiAuFe16}
A.~Riffnaller-Schiefer, U.~H. Augsd\"orfer, and D.W. Fellner.
\newblock Isogeometric shell analysis with {NURBS} compatible subdivision
  surfaces.
\newblock {\em Applied Mathematics and Computation}, 272:139--147, 2016.

\bibitem{SaJu21}
A.~Sailer and B.~J\"uttler.
\newblock Approximately ${C}^1$-smooth isogeometric functions on two-patch
  domains.
\newblock In {\em Isogeometric Analysis and Applications 2018}, pages 157--175.
  Springer, LNCSE, 2021.

\bibitem{SchDiWoKlHe19}
S.~Schu{\ss}, M.~Dittmann, B.~Wohlmuth, S.~Klinkel, and C.~Hesch.
\newblock Multi-patch isogeometric analysis for {K}irchhoff–{L}ove shell
  elements.
\newblock {\em Comput. Methods Appl. Mech. Engrg.}, 349:91--116, 2019.

\bibitem{SoTa23}
J.~Sogn and S.~Takacs.
\newblock Stable discretizations and {IETI}-{DP} solvers for the {Stokes}
  system in multi-patch {IgA}.
\newblock {\em ESAIM, Math. Model. Numer. Anal.}, 57(2):921--952, 2023.

\bibitem{SoglTakacs_IETI_Elasticity}
J.~Sogn and S.~Takacs.
\newblock Isogeometric tearing and interconnecting solvers for linear
  elasticity in multi-patch isogeometric analysis with theory for two
  dimensional domains.
\newblock {\em Comput. Methods Appl. Mech. Engrg.}, 418:116482, 2024.

\bibitem{ToSpHu17}
D.~Toshniwal, H.~Speleers, and T.~J.~R. Hughes.
\newblock Smooth cubic spline spaces on unstructured quadrilateral meshes with
  particular emphasis on extraordinary points: Geometric design and
  isogeometric analysis considerations.
\newblock {\em Comput. Methods Appl. Mech. Engrg.}, 327:411--458, 2017.

\bibitem{GradientDamageModels}
C.~V. Verhoosel, M.~A. Scott, T.~J.~R. Hughes, and R.~de~Borst.
\newblock An isogeometric analysis approach to gradient damage models.
\newblock {\em Internat. J. Numer. Methods Engrg.}, 86(1):115--134, 2011.

\bibitem{WeLiQiHuZhCa22}
X.~Wei, X~Li, K.~Qian, T.~J.~R Hughes, Y.~J Zhang, and H.~Casquero.
\newblock Analysis-suitable unstructured {T}-splines: {M}ultiple extraordinary
  points per face.
\newblock {\em Comput. Methods Appl. Mech. Engrg.}, 391:114494, 2022.

\bibitem{WeTa21}
P.~Weinm\"uller and T.~Takacs.
\newblock {Construction of approximate ${C}^1$ bases for isogeometric analysis
  on two-patch domains}.
\newblock {\em Comput. Methods Appl. Mech. Engrg.}, 385:114017, 2021.

\bibitem{WeTa22}
P.~Weinm\"uller and T.~Takacs.
\newblock {An approximate ${C}^1$ multi-patch space for isogeometric analysis
  with a comparison to Nitsche’s method}.
\newblock {\em Comp. Methods Appl. Mech. Engrg.}, 401:115592, 2022.

\bibitem{WiZaScPa21}
O.~B. Widlund, S.~Zampini, S.~Scacchi, and L.~F. Pavarino.
\newblock Block {FETI}-{DP}/{BDDC} preconditioners for mixed isogeometric
  discretizations of three-dimensional almost incompressible elasticity.
\newblock {\em Math. Comput.}, 90(330):1773--1797, 2021.

\bibitem{WeFaLiWeCa23}
Wen. Z., Md.~S. Faruque, X.~Li, X.~Wei, and H.~Casquero.
\newblock Isogeometric analysis using {G}-spline surfaces with arbitrary
  unstructured quadrilateral layout.
\newblock {\em Comput. Methods Appl. Mech. Engrg.}, 408:115965, 2023.

\bibitem{ZhSaCi18}
Q.~Zhang, M.~Sabin, and F.~Cirak.
\newblock Subdivision surfaces with isogeometric analysis adapted refinement
  weights.
\newblock {\em Computer-Aided Design}, 102:104--114, 2018.

\end{thebibliography}
%

\end{document}